\documentclass[12pt,draft,twoside]{article}
\usepackage[cp1250]{inputenc}
\usepackage{latexsym,amssymb,amsmath}
\title{On rectifiable spaces and its algebraical equivalents, topological
algebraic systems and Mal'cev algebras}
\date{}
\author{N. I. Sandu}
\begin{document}
\maketitle

\begin{abstract}
We investigate the rectifiable spaces, the Mal'cev algebras, the
almost quasivarieties of topological algebraic systems and their
free systems and others. It  specifies and  corrects the roughest
mistakes,  incorrect statements and nonsense  of the introduced
concepts connected with concepts listed before, which are
available in numerous papers on topological algebraic systems,
basically in papers of Academician Choban M. M. and his disciples.
\smallskip\\

Mal'cev algebra, arity of operation, topological
algebraic system, quasiatomic formula, hereditary class,
multiplicative closed class, unitary class, replica complete
class, variety, quasivariety, almost quasivariety, topological
algebraic system with given defining topological space and given
defining relations, topological free system of almost
quasivariety.
\smallskip\\
\textbf{Mathematics of subject classification}: 08A05; 22A20;
54E35; \break 54A25; 08B05
\end{abstract}

By  \cite{Gul} a topological space $X$ is said to be
\textit{rectifiable}  or \textit{a space with rectifiable
diagonal} (the terminology by \cite{Usp1}, \cite{Usp2}) provided
that there is a  homeomorphism $\Phi : X \times X \rightarrow X
\times X$ of $X \times X$ onto itself and an element $e \in X$
such that $\pi_1 \circ \Phi = \pi_1$ and for every $x \in X$ the
equality $\Phi(x,x) = (x,e)$ is fulfilled, where $\pi_1: X^2
\rightarrow X$ is the projection on the first coordinate.

As  mentioned in \cite{Gul}, the notion of rectifiable space was
introduced at the seminar of Prof. Arhangel'skii at the Moscow
State University (see \cite{Usp1}, \cite{Usp2}, \cite{Gul}). A
more general notion than  rectifiable space is the notion of
centralizable space introduced in \cite{Shap} (see, also,
\cite{Chob1}, \cite{Chob2}): in the definition of centralizable
space is not required for the homeomorphism  $g : X \times X
\rightarrow X \times X$ to be surjective.  These definitions have
a  topological character. Rectifiable spaces turn out to be a good
generalization of topological groups. In general, in papers,
containing problems connected with rectifiable spaces, many known
results from topological groups which are  generalized by
topological methods, are shown to  remain valid for rectifiable
spaces. But these papers have a shortcoming: they do not  shown
the algebraic structure, which  the notion of rectifiable space
corresponds to.

This paper elucidates  the algebraic notions of quasigroup and
loop   and  shows  that a topological space $X$ is rectifiable if
and only if $X$ is a topological right loop (implicitly this
relation is used in \cite{Gul}, \cite{LS}, \cite{LLL}). Thus, in
order to study rectifiable spaces the powerful methods of theory
of topological quasigroups and loops will be applied (see, for
example, \cite{Smith1}), which were  earlier ignored. Then the
results of rectifiable spaces become  part of the theory of
topological loops, and many proofs of these results become more
simple. Unfortunately, in some works related to problems connected
with rectifiable spaces \cite{Chob1}, \cite{Chob2}, \cite{AC2},
\cite{AC1}, \cite{ACM1}, \cite{ACM2} serious errors are admitted.
Using the fact that any rectifiable space is a topological right
loop some of these errors become evident.

Others gross blunders, connected with study of rectifiable spaces,
will be exposed and corrected  below. In particular, a rectifiable
spaces is a Mal'cev algebra. According to \cite{Mal2} a Mal'cev
algebra is characterized as an  algebra with permutable
congruences. Proceeding from this  we underlined on roughest
mistakes and non-sense of some concepts, connected with the notion
of Mal'cev algebra, considered in \cite{AC1}, \cite{Chob1},
\cite{Chob2}, \cite{Chob3}, \cite{Chob4}, \cite{Chob12},
\cite{Chob6}, \cite{AC2}, \cite{ACM1}, \cite{ACM2}, \cite{ACM3},
\cite{Usp1}, \cite{Usp2}, \cite{Usp3} and others.

In Section 2  the notions of topological algebra of given defining
relations (identities) and topological free algebra of a variety,
considered in \cite{Mal3}, are generalized for topological
algebraic systems. It  introduces the notion of almost
quasivariety of topological algebraic systems of given signature
as a class of  systems which is closed with respect to Cartesian
products, subsystems and contains an unitary system. Let $\frak K$
be a class of topological algebraic systems of fixed signature
$\Omega$, which is closed with respect to Tychonoff topology and
subsystems. It  proves that every system  of $\frak K$ can be
given by defining topological space $X$ and by defining relations,
which are  a totality of quasiatomic formulas, if and only if
$\frak K$ is almost quasivariety (Theorem 5).

Let $\frak K$ be a class of topological algebraic systems of given
signature with such a defining topological space that $\frak K$ is
closed with respect to Tychonoff topology and subsystems. In
Section 3 is proved.

The  class  $\frak K$  contains a $\frak K$-free systems when and
only when $\frak K$ is an almost quasivariety. More specific, when
and only when the non-trivial almost quasivarieties contain free
topological algebraic systems $\mathcal{F}_m$ of any given rank $m
\geq 1$ (Theorem 7).

Let $\frak K$ be an almost quasivariety  of topological algebraic
systems of fixed signature $\Sigma$. Then:

1) an algebraic system $\mathcal{A}$, free in relation to class
$\frak K$, is  free in relation to any subclass $\mathcal{L}
\subseteq \mathcal{K}$, and  in relation to closure $HS\prod
\mathcal{K}$;

2) all  free $\frak K$-systems of a given rank $M$ are
topologically  isomorphic among each other and  any topological
algebraic $\frak K$-system topologically generated by a set of
cardinality $m$ is an image of a continuous homomorphism of a free
systems $\mathcal{F}_m(\frak K)$ of rank $m$;

3) a free basis of a free system $\mathcal{A} =
\mathcal{F}_m(\mathcal{A}$) of  some class $\frak K$ of
topological algebraic system is a minimal generating set in
$\mathcal{A}$ (Theorem 8).

Using the listed above results in Section 5  underlines the
roughest mistakes (some of them will be   corrected) in
definitions of classical algebraical notions,the non-sense of the
introduced notions such as continuous signature, free algebras,
various types of varieties, quasivarieties of topological
algebraic system of continuous signature (for example, a class of
algebras is called quasivariety if this class is closed with
respect to Tychonoff products and subalgebras) and others, which
make the basis of the monographes, the dissertations and numerous
papers. Some of these works will be specified during the analysis.

\section{Preliminaries. Topological algebraic systems, quasivarieties}

For the general algebraic theory of quasigroups and loops see the
works \cite{Bruk}, \cite{Bel}, \cite{Smith1} and for general
theory of topology see the works \cite{Kelli}, \cite{Eng}. But
first, following \cite{Mal1} (see, also, \cite{Con}, \cite{Mal2})
some facts from the theory of algebraic systems, particularly,
$\Omega$-algebras, will be given.

Let $n \in N = \{0, 1, 2, \ldots\}$. A mapping $F_{i}$, $i \in I$,
which   maps every  ordered sequence $x_1, \ldots, x_n$ of $n$
elements of the set $A$ to  an unique  defined by  element,
$F_i(x_1, \ldots, x_n)$, is called an \textit{ $n$-ary operation},
defined on $A$, and the integer $n$ is called \textit{an arity} of
operation $F_i$. According to \cite{Mal1}, an \textit{algebraic
system} is called the object $\mathcal{A} = <A, \Omega_F,
\Omega_P>$, consisting of three sets: a non-empty \textit{basic
set} $A$, a set of operations $\Omega_F = \{f_0, \ldots, f_{\xi},
\ldots \}$ defined on set $A$, and a set of predicates  $\Omega_P
= \{p_0, \ldots, p_{\eta}, \ldots \}$, given on set $A$. If
$\Omega_P = \emptyset$ and $\Omega _F = \Omega = \cup_{n \in
N}\Omega_n$, where $\Omega_n$ is the set of all operations of
arity  $n$,  then $\mathcal{A} = <A, \Omega>$, is called
\textit{an algebra of  signature $\Omega$}, or
\textit{$\Omega$-algebra}. The operations $f_{\xi} \in \Omega_F$
and the predicates $p_{\eta} \in \Omega_P$ are called
\textit{basic} \cite{Con}, \cite{Mal1}.

Further, unless otherwise stipulated, we shall examine only an
algebraic systems  with basic operations of finite arity
 and follow, basically, terminology
from \cite{Mal1}.

A mapping of the basic set of algebraic system $\mathcal{A}$ into
the basic set of algebraic system $\mathcal{B}$ is called a
\textit{ mapping} of  algebraic system $\mathcal{A}$ into
algebraic system $\mathcal{B}$. A mapping $\varphi$ of algebraic
system $\mathcal{A} = <A, \{f_i\}, \{p_j\}>$ into algebraic system
$\mathcal{B} = <B, \{g_i\}, \{q_j\}>$ of the same type  is called
a \textit{homomorphism} of algebraic system $\mathcal{A}$ into
algebraic system $\mathcal{B}$ if
$$\varphi f_i(x_1, \ldots, x_{m_i}) = g_i(\varphi x_1, \ldots,
\varphi x_{m_i}),$$ $$p_j(x_1, \ldots, x_{m_j}) \Rightarrow
q_j(\varphi x_1, \ldots, \varphi x_{m_i}),\eqno{(1)}$$ for all
$x_1, \ldots, x_{m_r} \in A$, where $r = i;  j$. The relation
$\theta: x \theta y \Leftrightarrow \varphi x = \varphi y$ is
called \textit{nuclear  equivalence} of homomorphism $\varphi$.

A homomorphism $\varphi$ of system $\mathcal{A} = <A, \{f_i\},
\{p_j\}>$ on system $\mathcal{B} = <B, \{g_i\}, \{q_j\}>$ is
called \textit{strong}, if for all elements $b_1, \ldots, b_{m_j}
\in B$ and every basic predicate $q_j$ from $q_j(b_1, \ldots,
b_{m_j} = \{true\}$ it follows the existence in $A$ of inverse
images $a_1, \ldots, a_{m_j}$ of elements $b_1, \ldots, b_{m_j}$
such  that $p_j(a_1, \ldots, a_{m_j} = \{true\}$

Let $\mathcal{A} = <A, \Omega_F, \Omega_P>$ be an algebraic
system. A relation $p(x_1, \ldots, x_n)$ of set $A$ is called
\textit{stable} on algebraic system $\mathcal{A}$ if for any
$m$-ary operation $f$ and for any elements $a_{i1}, a_{i2},
\ldots, a_{in} \in A$ ($i = 1, 2, \ldots, m$) from validity of
relations $p(a_{i1}, a_{i2}, \ldots, a_{in})$ ($i = 1, 2, \ldots,
m$) the validity of relation $p(f(a_{11}, \ldots, a_{m1}), \ldots,
f(a_{1n},  \ldots, a_{mn}))$ follows. A relation of equivalence
$\theta$ of set $A$ is called \textit{congruence} of algebraic
system $\mathcal{A}$ if $\theta$ is stable with respect to every
basic operation of system $\mathcal{A}$. The product of two
congruences $\theta$, $\xi $ on system $\mathcal{A}$ is a
congruence on $\mathcal{A}$ if and only if $\theta$, $\xi $ are
permutable, i.e. $\theta\xi  = \xi \theta$ (\cite[Corollary
I.2.1]{Mal1}.

Let $\theta$ be a congruence on algebraic system $\mathcal{A} =
<A, \Omega_F, \Omega_P>$. We denote $[a] = \{x \in A \vert x
\theta a\}$ for $a \in A$. On the totality $\mathcal{A}/\theta$ of
all classes of congruence we define the operations
$$f_n^{\star}([a_1], \ldots, [a_n]) = [f_n(a_1, \ldots, a_n)] \eqno{(2)}$$
for $f_n \in \Omega_F$ and the predicates, believing
$$p_k^{\star}([a_1], \ldots, [a_k]) \eqno{(3)}$$
true if there exist such  elements $a_j^{\prime} \theta a_j$ that
$p_k(a_1^{\prime}, a_k^{\prime})$ is true for $p_k \in \Omega_P$.
With respect to the introduced operations and predicates we get
the \textit{quotient system} $\mathcal{A}/\theta$. The canonical
mapping $\varphi: a \rightarrow [a]$ of system $\mathcal{A}$ on
quotient system $\mathcal{A}/\theta$ is a strong homomorphism. The
$\Omega$-algebras have no  predicates, so for $\Omega$-algebras
the notions of homomorphism and strong homomorphism coincide.
Homomorphisms can exist for models, but are not strong.

Further up to  Proposition 1, not breaking a generality, for an
algebraic system $\mathcal{A} = <A, \Omega_F, \Omega_P>$ let's
consider that $\Omega_P = \emptyset$ and let's consider the
algebraic system $\mathcal{A} = <A, \Omega_F, \Omega_P>$  as
$\Omega$-algebra  $\mathcal{A} = <A, \Omega>$.

For a set $A$, let $\mathcal{B}(A)$ denote the set of all its
subsets. Any subset of $\mathcal{B}(A)$ will be called
\textit{system} of subsets of set $A$. A system $\mathcal{L}$ of
subsets of set $A$ is called \textit{system of closures} if
$\cap\mathcal{D} \in \mathcal{L}$ for any subsystem $\mathcal{D}
\subseteq \mathcal{L}$. If $\mathcal{D} = \emptyset$ then $A \in
\mathcal{L}$. Let $X$ be a topological space and let $\mathcal{T}$
be the system of all closed subsets. Then $\mathcal{T}$ will be a
system of closures with property: $A \cup B \in \mathcal{T}$. Such
a system is called \textit{topological}. Conversely, if a
topological system of closures $\mathcal{T}$ on set $X$ is given
and $\emptyset \in \mathcal{T}$, than one can  define on $X$ a
topology with elements of $\mathcal{T}$ as closed sets. In
general, the received topology will not be separated.

\textit{A closure operator} on set $A$ is called a mapping $J$ of
$\mathcal{B}(A)$ in to itself, having the following properties: 1)
if $X \subseteq Y$, then $J(X) \subseteq J(Y)$; 2) $X \subseteq
J(X)$; 3) $JJ(X) = J(X)$. Every system of closures $\mathcal{L}$
on  $A$ defines a closure operator $J$ on $A$ by the rule $J(X) =
\cap\{Y \in \mathcal{L} \vert Y \supseteq X\}$. Conversely, every
closure operator  $J$ on $A$ defines a system of closures
$\mathcal{L} = \{X \subseteq A\vert J(X) = X\}$ and the so defined
mapping  is one-to-one \cite[Theorem II. 1.1]{Con}.

Let $\mathcal{A} = <A, \Omega>$ be an $\Omega$-algebra and let
$(\mathcal{A}_{\lambda})_{\lambda \in \Lambda}$ be a family of
subalgebras of $\mathcal{A}$. For every basic operation $\omega
\in \Omega$ any algebra $\mathcal{A}_{\lambda}$ is closed under
$\omega$, hence the intersection $\cap\mathcal{A}_{\lambda}$ also
is closed under $\omega$. Denote the corresponding closure
operator by $J_{\Omega}$. Consequently, $J_{\Omega}(X)$ will be
the intersection $\cap\mathcal{A}_{\lambda}$  of all subalgebras
$\mathcal{A}_{\lambda}$ of $\mathcal{A}$ which contain the set $X
\subseteq A$. By definition $J_{\Omega}(X) =
\cap\mathcal{A}_{\lambda}$ is called \textit{ a subalgebra of
$\Omega$-algebra $\mathcal{A}$ generated by set $X \subseteq A$}.

A closure operator $J$ on $A$ is called \textit{algebraic} if for
any $X \subseteq A$ and $a \in A$ the following holds:  $a \in
J(X)$ implies $a \in J(X_f)$ for some finite subset $X_f$ of set
$X$. The system of closures  is called \textit{algebraic} if the
corresponding closure operator  is algebraic. A system of closures
is algebraic if and only if it is  inductive, i.e. each chain in
$\mathcal{L}$ has a least upper bound in $\mathcal{L}$
\cite[Theorem II.1.2]{Con}.

A characteristic property of algebras considered widely in
literature (see, for example, \cite{Mal1}, \cite{Mal2},
\cite{Mal3}, \cite{Con}), is that  they are defined by $(n,
1)$-operations, where $n$ is a finite integer. Such a $(n,
1)$-operations is called \textit{finite-place operations}. Let's
name them \textit{operations of finite arity}. Further we
investigate only such algebraic system $\mathcal{A} = <A,
\Omega_F, \Omega_P>$ for which the basic operations in $\Omega_F$
are finite-place. The following is of crucial importance for such
a  requirement.

Let $\mathcal{A} = <A, \Omega>$ be an $\Omega$-algebra with basic
operations of finite arity and let $X \subseteq A$ be a subset of
set $A$. By induction on $k$ we define the subset $X_k$ of
$\Omega$-algebra $\mathcal{A}$: $X_0 = X$; $X_{k+1} = \{x \in
\mathcal{A}$ \text{or} $x = \omega(a)$ \text{for some} $a \in
X_k^n$ \text{and} $\omega \in \Omega_n\}$. According to
\cite[Proposition II.5.1]{Con}, \cite[Theorem I.1.2]{Mal1} and
\cite[pag. 94]{Con}  the operator $J_{\Omega}$ is algebraical  and
$\cup_{k=0}^{\infty}X_k = J_{\Omega}(X)$ is the subalgebra of
$\Omega$-algebra $A$, algebraically generated by set $X$. We will
called the set $X$ \textit{set of algebraic generators}.
Conversely, for any given algebraic system of closures
$\mathcal{L}$ of a set $A$ it is possible to define  an
$\Omega$-algebra $\mathcal{A} = <A, \Omega>$ with basic operations
in $\Omega$ of finite arity such that $\mathcal{B}_{\Omega}(A) =
\mathcal{L}$, where $\mathcal{B}_{\Omega}(A)$ denote the set of
all subalgebras of $\Omega$-algebra $\mathcal{A} = <A, \Omega>$
(\cite[Theorem II.5.2]{Con}).

In \cite[Theorem II.5.2]{Con} it is possible to choose the
structure of $\Omega$-algebra by various ways. In particular,
 we shall state the statement \cite[Theorem II.5.6]{Con}
which  shows the special role of unary and $0$-ary operations for
research of algebraic closure operator.

Let's note that in \cite[pag. 56]{Mal1} the notion of algebraic
system of  closures, introduced in \cite{Con}, coincides with
notion of local totality, introduced in \cite{Mal1}. Remind the
last notion.  Let $\mathcal{A} = <A, \Omega>$ be an
$\Omega$-algebra and let $X \subseteq A$ be a subset of the set
$A$. A totality $\frak B = \{A_{\lambda} \vert \lambda \in
\Lambda\}$ of subsets $A_{\lambda}$ of set $A$ is called
\textit{locally in $X$} if any finite subset from $X$  contains in
some set $A_{\lambda}$ from $\frak B$. In such a case the set $B =
\cup A_{\lambda \in \Lambda}$ is the subalgebra of
$\Omega$-algebra $\mathcal{A} = <A, \Omega>$ algebraically
generated by algebraic generators $\{x \in X \vert X \subseteq
A\}$ (\cite[Theorem 1.2.2]{Mal1}).

It is easy to see that the union $C \cup D$ of subalgebras $C, D$
of an $\Omega$-algebra $\mathcal{A}$ with basic operations of
finite arity can not be a subalgebra. However, the union of
non-empty locally totality $\frak B$ of subalgebras of any
$\Omega$-algebra $\mathcal{A}$ is a subalgebra of $\Omega$-algebra
$\mathcal{A}$ (\cite[Theorem 1.2.2]{Mal1}). Moreover, the
following statement (\cite[Theorem II.5.6]{Con}) holds.

For an algebraic system of closures $\mathcal{L}$ on set $A$ the
following statements are equivalent:

(i) $\mathcal{L}$ is a sublattice of lattice $\mathcal{B}(A)$;

(ii) $\mathcal{L}$ is a topologic system of closures;

(iii) there is such a domain of  operations $\Omega$, whose the
arity of operations  has more unity (i.e. unary or 0-ary) and such
a structure of $\Omega$-algebra on $A$ that
$\mathcal{B}_{\Omega}(A) = \mathcal{L}$.

A \textit{polynomial} of letters $x_1, \ldots, x_t$  is called the
meaningful expression, consisting of these letters, brackets and
symbols of basic operations of the algebraic system. Then from
definition of sets $X_k$ (or $A_{\lambda \in \Lambda}$) it follows
that if the system $\mathcal{A}$ algebraically is generated by set
$X$ then for all $a \in \mathcal{A}$ there exist a basic operation
$\omega \in \Omega$ of arity $n$ and such an elements $b_1,
\ldots, b_n \in X$ that $a = \omega(b_1, \ldots, b_n)$. Hence any
element $a \in \mathcal{A}$ can be presented as polynomial of
letters $x_1, \ldots, x_n \in X$. Then from the above-stated it
follows\vspace*{0.1cm}

\textbf{Proposition 1.}  \textit{I). For an $\Omega$-algebra
$\mathcal{A} = <A, \Omega>$ and subset $X \subseteq A$ the
following statements are equivalent:}

\textit{a) $\mathcal{A}$ is algebraically generated by set $X$,
i.e. $X$ is a set of algebraic generators;}

\textit{b) any system of closures of set $\mathcal{B}(A)$ of all
subsets of set $A$ is algebraic;}

\textit{c) the closure operator  $J_{\Omega}(X) = \cup_{k =
0}^{\infty} X_k$, considered in \cite[Proposition II.5.1]{Con}, is
algebraic;}

\textit{d) the totality $\frak B = \{A_{\lambda} \vert \lambda \in
\Lambda\}$ of subsets $A_{\lambda}$ of set $A$, considered in
\cite[Theorem 1.2.2]{Mal1}, is locally in $X$;}

\textit{f) every element $a \in A$ can be presented with respect
to basic operations from $\Omega$  as polynomial of variables
$x_1, \ldots, x_r \in X$.}

\textit{II). The equivalent statements a) - f) of item I) hold if
any basic operation $\omega \in \Omega$ is a finite arity, i.e.
has a form $\omega(x_1, \ldots, x_n) = x_{n+1}$ (\cite[Proposition
II.5.1]{Con}, \cite[Theorem I.2.2]{Mal1}. Conversely, if for an
algebra with basic set $A$ hold the statements a) - f) of item I)
then on set $A$  it is possible to define such an $\Omega$-algebra
$\mathcal{A} = <A, \Omega>$ with basic operations in $\Omega$ of
finite arity that $\mathcal{B}_{\Omega}(A) = \mathcal{L}$
(\cite[Theorem II.5.2]{Con}).}\vspace*{0.1cm}

From Proposition 1 it follows that the study of algebraical
systems with basic operations of finite arity and algebraical
systems with basic operations of non-finite arity requires
different  approaches. The expression $''$a subalgebra of
$\mathcal{A}$ is algebraically generated by set $X \subseteq
A$$''$  is meaningful if and only if the basic
 operations in $\Omega_F$ are finite arity by item I).
According to item II) only in such a  case  every element $a \in
A$ can be presented with respect to the basic operations from
$\Omega$ as polynomial of variables $x_1, \ldots, x_r \in X$. It
is impossible to apply to the second case the algebraic closed
operators as by item I) this notion characterizes the operations
of finite arity.

Probably, there can not be a general way of researching the
operations of non-finite arity and thus in every concrete case the
achieved simplification would be significant.  \cite[pag. 69]{Con}
stipulates that so far the study of  operations with non-finite
arity did not pay  enough attention to mathematics. As shown
above, for this it is possible to study the operations with
non-finite arity with the help of topological spaces. Inversely,
as in \cite{Mal3},  any topological space could be considered as
discrete algebras with systems of partial operations of infinite
signature. But the properties of such partial operation are not
expressed as identities, but as operations of limiting
transitions. \vspace*{0.1cm}

Let $\mathcal{A} = <A, \Omega_F, \Omega_P>$ be an algebraic system
of  signature $\Omega = \Omega_F \cup \Omega_P$,  let $x_1, x_2,
\ldots, x_k$ be object variables and let $f_m^{(n)}$, $m, n = 0,
1, \ldots$, be  functional variables, where $m$ is the number of
order, $n$ is the arity of variable.  We define: 1) every word of
form $x_i$ or $f_i^{(0)}$ is a term; 2) if $a_1, \ldots, a_n$ are
 terms, then $f^{(n)}_1(a_1, \ldots, a_n)$ is a term; 3) a word
is a term if it is term by items 1), 2). The formulas of form
$$P(f_1(x_1, \ldots, x_k), \ldots, f_n(x_1, \ldots, x_k)),$$
$$f(x_1, \ldots, x_k) = g(x_1, \ldots, x_k), \eqno{(4)}$$
where $f, g, f_1, \ldots, f_n$ are some terms of signature
$\Omega$, $P \in \Omega$, are called \textit{quasiatomic} from
variables $x_1, x_2, \ldots, x_k$.

Let $S_1(x_1, \ldots, x_k), \ldots, S_{k+1}(x_1, \ldots, x_k)$ be
some quasiatomic formulas of signature $\Omega$ (switching
equality) from variables $x_1, \ldots, x_k$. A formula of form
$(\forall x_1, \ldots, x_k)S_1(x_1, \ldots, x_k)$ is called
\textit{identity} and a formula of form \break $(\forall x_1,
\ldots, x_k)(S_1 \vee \ldots \vee S_k  \rightarrow S_{k+1})$ is
called \textit{quasiidentity}.

Particularly, if  classes of algebras are considered, then
$\Omega_P = \emptyset$ and the identities have the form $(\forall
x_1, \ldots, x_k))(f(x_1, \ldots, x_k) = g(x_1, \ldots, x_k))$ and
the quasiidentities have the form $(\forall x_1, \ldots, x_k))(f_1
= g_1, \vee, \ldots, \vee f_r = g_r \rightarrow f_{r+1} =
g_{r+1})$.

A class $\mathcal{K}$ of algebraic  systems  is called \textit{a
variety} (respect. \textit{a quasivariety}) if such   a totality
$T$ of identities (respect. quasiidentities) of signature $\Omega$
exists that $\mathcal{K}$ consists of those and only those
algebraic systems of signature $\Omega$, where all  formulae from
$T$ hold true. The totality $T$ is called \textit{the defining
set} of variety (respect. quasivariety) $\mathcal{K}$ . Since  it
is possible to regard anyidentity as quasiidentity, then any
variety can be regarded  as a quasivariety.

We give  some notions and results from \cite{Mal1}. An algebraic
system is called \textit{unitary}  if it is  from one element and
all its basic predicates have the value  true.
 A class $\mathcal{K}$ of algebraic systems is called
\textit{hereditary} if $\mathcal{K}$  is closed with respect to
the subsystems of its systems, and is called
\textit{multiplicatively closed} if $\mathcal{K}$ is closed with
respect to the Cartesian product.\vspace*{0.1cm}

\textbf{Proposition 2.} (Birkhoff's Theorem). \cite[Theorem
VI.13.1]{Mal1}. \textit{A class $\frak{V}$ of algebraic systems of
a fixed signature $\Omega$ is a variety if and only if:}

\textit{$V_1$)  class $\frak{V}$ is hereditary;}

\textit{$V_2$)  class $\frak{V}$ is multiplicatively closed;}

\textit{$V_3$) every homomorphic image of any algebraic system
from $\frak{V}$ belongs to $\frak{V}$.}\vspace*{0.1cm}

\textbf{Proposition 3} \cite[Corollary V.11.3]{Mal1}. \textit{A
class  $\frak{Q}$ of algebraic systems of fixed signature $\Omega$
is a quasivariety if and only if:}

\textit{$Q_1$)  class  $\frak{Q}$  is  ultraclosed;}

\textit{$Q_2$)  class  $\frak{Q}$  is hereditary;}

\textit{$Q_3$)  class  $\frak{Q}$  is multiplicatively closed;}

\textit{$Q_4$)  class  $\frak{Q}$ contains an unitary
system.}\vspace*{0.1cm}

Remind that a formula is called \textit{closed} if it does not
contain any free subject variables. By \cite[Theorem
III.6.1]{Mal1} it is possible to define  that a class
$\mathcal{K}$ of algebraic systems of signature $\Omega$ is called
\textit{axiomatizable} if and only if  such a totality of closed
formulaes $S$  exists that $\mathcal{K}$ consists of those and
only those systems of signature $\Omega$ for which  the formulae
from $S$ hold true. Every axiomatizable class of systems
$\mathcal{K}$ is ultraclosed \cite[Corollary IV.8.10]{Mal1}.
Moreover, an axiomatizable class of systems $\mathcal{K}$
satisfies the conditions $Q_2)$, $Q_3)$ $Q_4)$ if and only if the
class $\mathcal{K}$ is a quasivariety \cite[Corollary
V.11.7]{Mal1}.

Recall that an algebraic system for which the basic set of
elements is a topological space and the basic operations are
continuous is called \textit{topological algebraic system}. In
paper \cite{ProtSid} a variety $\mathcal{K}$  is defined as class
of topological algebraic systems, which satisfy  some set of
limiting identities $\Sigma$. A variety $\mathcal{K}$ is called
primitive class if $\Sigma$ consists of algebraical identities.
For thus defined notions the analogues of the Offers  of
Propositions 2, 3 are proved. For a class $ \mathcal{K}$ to be a
variety is it is necessary and sufficient for $\mathcal{K}$ to be
closed with respect to:  closed subsystems of system from
$\mathcal{K}$;  images of continuous homomorphisms:  Tychonoff
products. A variety $\mathcal{M}$ of topological groups is a
primitive class  if and only if   $\mathcal{M}$ is closed with
respect to ultraproduct.

Let us note that the condition $Q_1$) is essentially for the
definition of quasivariety.  \cite[pag. 295]{Mal1}  gives a class
of groups $\mathcal{G}$ for which  the conditions $Q_2)$, $Q_3)$
$Q_4)$ hold, but $\mathcal{G}$ is not  a quasivariety. We also
mention that from \cite[Theorem V.11.4]{Mal1} it follows that the
condition $Q_4)$ is necessary in Proposition 3.

According to \cite[Theorem V.11.5]{Mal1}, \textit{a class of
algebraic systems $\mathcal{K}$ satisfies the conditions $Q_2)$,
$Q_3)$ $Q_4)$ if and only if the class $\mathcal{K}$ is replica
complete.} Let us remind this definition. Let $\frak K$ be a class
of algebraic systems of fixed signature $\Omega$ and let
$\mathcal{A}$ be  some system of signature $\Omega$, not
necessarily belonging to a class $\frak K$. A homomorphism
$\alpha$ of system $\mathcal{A}$ on some $\frak K$-system
$\mathcal{A}_1$ is called \textit{$\frak K$-morphism} if for every
homomorphism $\gamma$ of system $\mathcal{A}$ in any system $\frak
K$-system $\mathcal{B}$ such  a homomorphism $\varsigma$ of system
$\mathcal{A}_1$ on system $\mathcal{B}$ exists that $\gamma =
\varsigma\alpha$. Every $\frak K$-morphiń image of system
$\mathcal{A}$ is called \textit{replica} of $\mathcal{A}$ in class
$\frak K$ ($\frak K$-replica) and is denoted by
$\mathcal{A}_{\frak K}$.

A class $\frak K$ of signature $\Omega$ is called \textit{replica
complete}, if  every algebraic system has a replica within it.
According to \cite[Theorem V.11.5]{Mal1} we have\vspace*{0.1cm}

\textbf{Lemma 1.} (Lemma-definition).  \textit{A class $\frak K$
of algebraic systems of fixed signature $\Omega$ is replica
complete, if $\frak K$ satisfies the conditions $Q_2)$, $Q_3)$
$Q_4)$, i.e. is hereditary,  multiplicatively closed and contains
an unitary system. According to Proposition 3 such class of
algebraic systems $\frak K$ with listed properties is called an
almost quasivariety.}\vspace*{0.1cm}

From Proposition 3 it follows.\vspace*{0.1cm}

\textbf{Corollary 1.}  \textit{An almost quasivariety $\frak K$ is
a quasivariety  if and only if the   class $\frak K$ is
axiomatizable.}\vspace*{0.1cm}

 Let $\mathcal{A} = <A,\Omega_F, \Omega_P>$ be an algebraic
system generated by set $X$. Every polynomial of letters  $x_1,
\ldots, x_p \in X$ can be considered as a  $p$-ary operation.
Operations obtained with the help of polynomial are  called
\textit{derived}.

Transformations of the set of algebraic system $\mathcal{A}$,
having the form $x \rightarrow F(x)$, where $F(x)$, is a
polynomial of  $x$, are called  \textit{translations of the
system}. The translation $T$ is called \textit{reversible} if such
a translation $S$ exists  that $ST = TS = E$, where $E$ is the
identical mapping. All reversible translations form the
\textit{group of translations} of given system.\vspace*{0.1cm}

\textbf{Theorem 1}  (Mal'cev) \cite{Mal2}. \textit{All congruences
on every algebraic system of some variety are permutable iff   a
polynomial $\Psi(x,y,z)$ exists, satisfying the identities}
$$\Psi(x,x,z) = z, \quad \Psi(x,z,z) = z \eqno{(5)}$$ \textit{on
all systems of this class.}\vspace*{0.1cm}

\textit{A biternary system} of an algebraic system $\mathcal{A} =
<A, \Omega_F, \Omega_P>$ is a pair of ternary operations $\alpha,
\beta : A \times A \times A \rightarrow A$ such that
$\alpha(x,x,y) = y, \break  \alpha(\beta(x,y,z),y,z) = x,
\beta(\alpha(x,y,z),y,z) = x$ for all $x, y, z \in
A$.\vspace*{0.1cm}

\textbf{Theorem 2}  (Mal'cev) \cite{Mal2}. \textit{The groups of
reversible translations of all systems of some variety is
transitive iff  derived ternary operations exist, relative to
which the systems of a given class are biternary.}\vspace*{0.1cm}

\textbf{Theorem 3.} \textit{The Theorems 1, 2 hold for any almost
quasivariety  class of algebraic systems of a given
signature.}\vspace*{0.1cm}

The proof of Theorem 3 almost literally repeats the proofs of
Theorems 1 and 2. It is necessary only to use the Lemma 1 and the
Proposition 1.

\section{On almost quasivarieties of topological algebraic
systems}

A relation $\theta(x_1, \ldots, x_n)$ defined on topological space
is called \textit{continuous}, if for any sequence  of elements
$x_1, \ldots, x_n$ not satisfying the relation $\theta$ such  a
sequence  of open set $U_i$, $x_i \in U_i$ ($i = 1, \ldots, n$)
exists that for all $x_i^{\prime}$ ($x_i^{\prime} \in U_i; i = 1,
\ldots, n$) the relation $\theta(x_1^{\prime}, \ldots,
x_n^{\prime})$ does not hold \cite{Mal2}.

Let $\theta$ be an equivalence  of a set $A$ and let $N$ be a
subset of $A$. According to \cite{Mal2}, the set $\cup_{n \in
N}\{a \in A \vert a \theta n\}$ is called \textit{a
$\theta$-saturation} of $N$. An equivalence $\theta$ of a
topological space  is called \textit{complete  equivalence} if any
$\theta$-saturation of an open set is  open.

Recall that an operation (mapping) $f(x_1, x_2, \ldots, x_n)$ of a
topological space $G$ is continuous if for every neighborhood $V$
of point $f(x_1, x_2, \ldots, x_n)$ such  neighborhoods $U_1,
\ldots, U_n$ of points $x_1, x_2, \ldots, x_n$ exist that \break
$f(U_1, \ldots, U_n) \subset V$. \textit{A homeomorphism} is a
continuous one-to-one mapping $\varphi$ of some  topological space
$X$ on some topological space $Y$, the inverse mapping
$\varphi^{-1}$ is also  continuous.

An algebraic system for which the basic set of  elements is a
topological space and the basic operations are continuous is
called \textit{topological algebraic system}. The totality of
topological algebraic systems which form a variety (respect.
quasivariety or almost quasivariety, or axiomatizable class) in
algebraic sense is called \textit{a variety} (respect. \textit{
quasivariety or almost quasivariety, or axiomatizable class})
\textit{of topological algebraic system}. \textit{Replica of class
of topological algebraic systems} is defined in natural way: all
mappings in algebraic definition should be  continuous.

Let $\mathcal{A} = <A, \{f_i\}, \{p_j\}>$,  $\mathcal{B} = <B,
\{g_i\}, \{q_j\}>$ be topological algebraic systems of the  same
type. A homomorphism $\varphi$ of  $\mathcal{A}$ into
$\mathcal{B}$ is defined by (1). If the mapping $\varphi$ of
topological space $A$ into topological space $B$ is continuous,
then  $\varphi$ is called \textit{continuous homomorphism} of
topological algebraic systems.

If $\theta$ is a congruence on topological algebraic system
$\mathcal{A} = <A, \{f_i\}, \{p_j\}>$, then by (2), (3) the
quotient set $A/\theta$ turns into an algebraic system. We define
the images of open sets of topological space $A$ under
homomorphism $A \rightarrow A/\theta$ as open sets in $A/\theta$.
Then the algebraic system $\mathcal{A}/\theta$ turns into
topological algebraic system and $\mathcal{A}\rightarrow
\mathcal{A}/\theta$ will be a continuous homomorphism. If $\theta$
is a complete congruence, then $\mathcal{A}\rightarrow
\mathcal{A}/\theta$ will be an open continuous homomorphism.

If $\theta$ is a complete and continuous congruence of topological
algebraic system $\mathcal{A}$ then the quotient system
$\mathcal{A}/\theta$ satisfies the axiom $T_2$, i.e. is a
Hausdorff space. Conversely, let $\varphi$ be an open and
continuous homomorphism of topological algebraic system
$\mathcal{A}$ on topological algebraic system $\mathcal{B}$, which
satisfies the axiom $T_2$. Then the nuclear congruence of
$\varphi$ on $\mathcal{A}$ is complete and continuous and the
mapping $\mathcal{A}/\theta \leftrightarrow \mathcal{B}$ is an
open and continuous isomorphism (see, \cite{Mal2}).

Further all topological systems are assumed to have  axiom
$T_2$.\vspace*{0.1cm}

 \textbf{Proposition 4.} \textit{A nuclear equivalence $\theta$
of every continuous  homomorphism $\varphi$ of topological
algebraic system $\mathcal{A}$ on  topological algebraic system
$\mathcal{B}$ of the same type is a congruence on $\mathcal{A}$
and the canonical mapping $\tau: \mathcal{A}/\theta \rightarrow
\mathcal{B}$ is a continuous homomorphism. If the homomorphism
$\varphi$ is open and continuous  then the canonical mapping
$\tau: \mathcal{A}/\theta \rightarrow \mathcal{B}$ is an open
continuous homomorphism. If the homomorphism $\varphi$ is strong,
open and continuous then the canonical mapping $\tau:
\mathcal{A}/\theta \rightarrow \mathcal{B}$ is an open continuous
isomorphism.}\vspace*{0.1cm}

\textbf{Proof.} For algebraic systems the Proposition 4 coincides
with Theorem I.2.1 from \cite{Mal1}. We consider the topological
case and let the homomorphism $\varphi$ of $\mathcal{A}$ on
$\mathcal{B}$ be continuous. As  shown above, the homomorphism
$\alpha$ of $\mathcal{A}$ on $\mathcal{A}/\theta$ is continuous.

Let $z$ be an element in $\varphi \mathcal{A} = \mathcal{B}$. Then
such  an elements $y \in \mathcal{A}/\theta$, $x \in \mathcal{A}$
exist that $\tau(y) = z$, $\alpha(x) = y$, $\varphi(x) = z$. Let
$W$ be a neighborhood  of $z$ in $\mathcal{B}$. From the
definition of topology of $\mathcal{A}/\theta$ it follows that
such a neighborhood $U$ of $x$ in $\mathcal{A}$ exists  that
$\varphi(U) \subseteq W$. The homomorphism $\alpha$ is open. Then
$\alpha(U) = V$ is a neighborhood  of $y$ in $\mathcal{A}/\theta$
and $\tau(V) \subseteq W$. Hence $\tau$ is a continuous
homomorphism.

Let now $\varphi$ be an open continuous homomorphism. Then the
congruence $\theta$ is complete and the continuous homomorphism
$\alpha: \mathcal{A} \rightarrow \mathcal{A}/\theta$ is open. Let
$V$ be an open set in $\mathcal{A}/\theta$. We denote $U =
\alpha^{-1}(V)$. As $\theta$ is a complete congruence then
$\alpha(U) = V$ and $U$ is an open set in $\mathcal{A}$. The
homomorphism $\varphi$ is open. Then the set $\varphi(U)$ is open
in $\mathcal{B}$. Hence the set $\tau(V) = \tau\alpha(U) =
\varphi(U)$ is open in $\mathcal{B}$ and the continuous
homomorphism $\tau: \mathcal{A}/\theta \leftrightarrow
\mathcal{B}$  is open. If the homomorphism $\varphi$ is strong,
then by \cite[Theorem I.2.1]{Mal1}, the homomorphism of algebraic
systems $\tau: \mathcal{A}/\theta \leftrightarrow \mathcal{B}$ is
an isomorphism.
 This completes the proof of Proposition 4.\vspace*{0.1cm}

The Proposition 4 shows that the totality of all strong, open and
continuous homomorphic images of a given algebraic system
$\mathcal{A}$ up  to homeomorphism is settled by totality of all
quotient systems on its various complete congruences. If
$\mathcal{A}$ is an $\Omega$-algebra, then a similar  statement
holds for open and continuous homomorphic images of $\mathcal{A}$
and complete congruences of $\mathcal{A}$.

The following Theorem 4 is proved in \cite[Theorem 10]{Mal2} for a
varieties of topological $\Omega$-algebras, i.e. for a varieties
of topological algebraic systems $(A, \Omega_F, \Omega_P)$ with
$\Omega_P = \emptyset$. The proof of Theorem 4 literally repeats
the proof of Theorems 10 and 11 from \cite{Mal2}: for this purpose
it is sufficient to use the Theorem 3 instead of Theorem
1.\vspace*{0.1cm}

\textbf{Theorem 4.} \textit{If the congruences are permutable on
all systems of an almost quasivariety $\frak{A}$  then all
congruences are complete  on topological   systems of class
$\frak{A}$ and  those and only those congruences are continuous on
topological systems of class $\frak{A}$, where adjacent classes
are closed.}\vspace*{0.1cm}

\textbf{Corollary 2.} \textit{Let $\frak{K}$ be an almost
quasivariety of topological algebraic systems with permutable
congruences and let $\mathcal{A} \in \frak{K}$. Then any strong
continuous homomorphism of $\mathcal{A}$ is an open,  strong and
continuous homomorphism of $\mathcal{A}$. If $\mathcal{A}$ is an
$\Omega$-algebra, then any continuous homomorphism of
$\mathcal{A}$ is an open  continuous homomorphism of
$\mathcal{A}$.}\vspace*{0.1cm}

\textbf{Remark 1.} Let $X$ and $Y$ be  topological spaces, $X
\times Y$ be their Cartesian product, and $f$ be a mapping of $X
\times Y$ into  a third topological space. The mapping $f$ is
called \textit{continuous on $x$} (respect. on $y$) if and only if
for every $y \in Y$ (respect. $x \in X$) the function $f(,y)$, the
value of which at  point $x$  equals $f(x,y)$, is continuous. If
the  function $f$ is continuous, then it is continuous with
respect to every variable. The opposite  converts statement is not
always true: there exist  non-continuous functions, with respect
to every variable \cite[pag. 143] {Kelli}.\vspace*{0.1cm}

According to (4), let
$$P_{\lambda}(f_{\lambda_1}(x_{\lambda_{11}}, \ldots, x_{\lambda_{1k}}),
\ldots, f_{\lambda_n}(x_{\lambda_{1n}}, \ldots,
x_{\lambda_{kn}}),$$
$$f_{\varsigma}(x_{\varsigma 1}, \ldots, x_{\varsigma k}) = g_{\varsigma}(x_{\varsigma 1},
\ldots, x_{\varsigma k}), \eqno{(6)}$$ where $\varsigma, \lambda
\in \Lambda$, $f_{\varsigma}, g_{\varsigma}, f_{\lambda_1},
\ldots, f_{\lambda_n}$ are some terms of signature $\Omega$,
$P_{\lambda} \in \Omega_P$, be a totality $S$ of quasiatomic
formulas  of variables $x_i \in X$ ($i \in I$). Let us note that
every term of any topological algebraic system is a continuous
operation.

By analogy with \cite{Mal3} we give.\vspace*{0.1cm}

\textbf{Definition 1.}  \textit{Let $\frak K$ be a class of
topological algebraic systems  of given signature $\Omega =
\Omega_F \cup \Omega_P$, let  $X$ be a topological space
 and let $S$ be a totality of  formulas $S$ of form (6)  of
signature $\Omega$ of variables $x_i \in X$ ($i \in I$). A
topological algebraic system $\mathcal{A} = <A, \Omega>$ of the
class $\frak K$ with the given continuous mapping $\sigma$ of
topological space $X$ in topological space $A$ will called
\textit{defined in $\frak K$ by defining space  $X$ and defining
relations $S$} if:}

$F_1$) \textit{in  $\mathcal{A}$ the formulas of $S$} hold
$$P_{\lambda}(f_{\lambda_1}(\sigma x_{\lambda_{11}}, \ldots, \sigma x_{\lambda_{1k}}),
\ldots, f_{\lambda_n}(\sigma x_{\lambda_{1n}}, \ldots, \sigma
x_{\lambda_{kn}}),$$
$$f_{\varsigma}(\sigma x_{\varsigma 1}, \ldots, \sigma x_{\varsigma k}) =
g_{\varsigma}(\sigma x_{\varsigma 1}, \ldots, \sigma x_{\varsigma
k});$$ \textit{are fulfilled;}

$F_2$) \textit{the system $\mathcal{A}$ is topologically generated
by images of elements of $X$, i.e.  algebra $< A, \Omega_F>$ does
not contain a closed subalgebra, not equal to  $<A, \Omega_F>$ and
containing $\sigma X$;}

$F_3$) \textit{for every continuous mapping  $\gamma$ of the space
$X$ into any topological algebraic system $\mathcal{C}$ of the
class $\frak{K}$ under which  the relations}
$$P_{\lambda}(f_{\lambda_1}(\gamma x_{\lambda_{11}}, \ldots, \gamma x_{\lambda_{1k}}),
\ldots, f_{\lambda_n}(\gamma x_{\lambda_{1n}}, \ldots, \gamma
x_{\lambda_{kn}}),$$
$$f_{\varsigma}(\gamma x_{\varsigma 1}, \ldots, \gamma x_{\varsigma k}) =
g_{\varsigma}(\gamma x_{\varsigma 1}, \ldots, \gamma x_{\varsigma
k});$$ \textit{are satisfied there exists  a continuous
homomorphism $\alpha$ of system $\mathcal{A}$ in to
$\mathcal{C}$, compatible with mappings $\sigma, \gamma$, i.e.
such that $\alpha\sigma(x) = \gamma(x)$ for all $x \in
X$.}\vspace*{0.1cm}

The proof of following Lemma 2 is similar to \cite[Theorem
2]{Mal3}. It is necessary only to use the Proposition 1, which
describes the subsystem of algebraic system $<A, \Omega>$,
algebraically generated by a subset of set $A$. Let us show the
proof of Lemma 2 to ease the reading.\vspace*{0.1cm}

\textbf{Lemma 2.}  \textit{The condition $F_2$) of Definition 1 is
equivalent to condition:}

$F^{\prime}_2)$ \textit{the algebraic system $\mathcal{A}$ is
algebraically generated  by the images of elements of $X$, i.e.
the system $\mathcal{A}$ is a totality of elements of set $A$,
expressed as  finite polynomials of images of elements from $X$
with respect to operations $\Omega_F$.}\vspace*{0.1cm}

\textbf{Proof.} We denote by $A^{\star}$ the totality of elements
of $A$, which are expressed as finite $\Omega_F$-polynomials of
imagines of elements of $X$. Leaving in  $A^{\star}$ the topology,
which it has a  subset of topological space $A$, we convert
$A^{\star}$ into a  topological algebraic system
$\mathcal{A}^{\star}$, where  the space $X$ is reflected
continuously in $\mathcal{A}^{\star}$ with the help of the same
mapping $\sigma$ as in the given  system $\mathcal{A}$. By $F_3$),
there is a continuous homomorphism $\tau$ of system $\mathcal{A}$
into system $\mathcal{A}^{\star}$, coordinated  with mappings $X
\rightarrow A$, $X \rightarrow \mathcal{A}^{\star}$. As
$\mathcal{A}^{\star}$ algebraically generate by images of elements
of  $X$ (Proposition 1), then $\tau$ will be the reflection
$\mathcal{A}$ on all system $\mathcal{A}^{\star}$, moreover, for
$a \in \mathcal{A}^{\star}$ we have $\tau(a) = a$. As
$\mathcal{A}^{\star}$ is dense â $\mathcal{A}$, then for each $a
\in \mathcal{A}$ there is  a directed set $\{a_{\lambda}\}$ of
elements of $\mathcal{A}^{\star}$, converging to  $a$. From
continuously of  $\tau$ it follows that $\tau(a) = \lim
\{\tau(a_{\lambda})\} = a$.

According to \cite[Proposition 1.6.7]{Eng} a topological space $X$
is Hausdorff if and only if any directed set in $X$ has no more
than one limit point. Hence $\mathcal{A} = \mathcal{A}^{\star}$,
as required.\vspace*{0.1cm}

\textbf{Corollary 3.} \textit{If in  some class  $\frak{K}$ of
topological algebraic systems of given signature, there exist
systems $\mathcal{A}$, $\mathcal{B}$, having the same  generalized
topological space $X$ and the same totality of  quasiatomic
formulas as their defining relations, then the systems
$\mathcal{A}$, $\mathcal{B}$ are  images of continuous mappings of
$X$ and  are topologically isomorphic.}\vspace*{0.1cm}

\textbf{Proof.} According to Definition 1 and Lemma 2 such
continuous mappings of defining topological space $X$
$$x_i \rightarrow a_i \quad (a_i \in \mathcal{A}, i \in I), \quad
x_i \rightarrow b_i \quad (b_i \in \mathcal{B}, i \in I)
\eqno{(7)}$$ exist that the elements $a_i$ algebraically generate
$\mathcal{A}$, the elements $b_i$ algebraically  generate
$\mathcal{B}$ and in suitable continuous homomorphisms
$$\varphi: \mathcal{A} \rightarrow \mathcal{B}, \quad
\psi: \mathcal{B} \rightarrow \mathcal{A}$$ we have $\varphi(a_i)
= b_i$, $\psi(b_i) = a_i$, hence
$$\psi\varphi(a_i) = a_i, \quad  \varphi\psi(b_i) = b_i (i \in I).
\eqno{(8)}$$

Let $f$ be a polynomial of signature $\Omega_F$. As the mappings
$\varphi\psi$ and $\psi\varphi$ are  continuous homomorphisms of
$\mathcal{A}$ and $\mathcal{B}$ into themselves, then from (8) it
follows that
$$\psi\varphi f(a_{\lambda}, \ldots, a_{\gamma}) = f(a_{\lambda}, \ldots,
a_{\gamma}), \quad \varphi\psi f(b_{\lambda}, \ldots, b_{\gamma})
= f(b_{\lambda}, \ldots, b_{\gamma}),$$ i.e. $\varphi\psi$,
$\psi\varphi$ are  identical mappings, $\psi = \varphi^{-1}$ and
$\varphi$ are the continuous isomorphism of $\mathcal{A}$ on
$\mathcal{B}$. Consequently,  the mapping $\varphi$ is a
topological isomorphism of $\mathcal{A}$ on $\mathcal{B}$,
moreover  the image of any  element of $X$ in $\mathcal{A}$ passes
to the image of this element in $\mathcal{B}$. This completes the
Proof of Corollary 3.\vspace*{0.1cm}

From Corollary 3 it follows that if a topological space is a
defining space for the same topological algebraic system then this
system is defined unequivocally within the accuracy of topological
isomorphism. Moreover, the following holds.

\textbf{Corollary 4.} \textit{Let $\frak K$ be a class of
topological algebraic systems of fixed signature and let $X$, $Y$
be the  defining spaces with continuous mappings $\alpha: X
\rightarrow \mathcal{A}$, $\beta: Y \rightarrow \mathcal{B}$  for
the systems $\mathcal{A}, \mathcal{B} \in \frak K$. If there is a
 continuous mapping $\varphi$ of $\mathcal{A}$ on $\mathcal{B}$
then there exists a continuous epimorphism $\xi: \mathcal{A}
\rightarrow \mathcal{B}$, satisfying the condition  $\xi\alpha(x)
= \beta\varphi(x)$, $x \in X$.}

\textbf{Proof.} We have that $\beta\varphi: X \rightarrow
\mathcal{B}$ and $\beta\varphi$ is a continuous mapping as a
product of continuous mappings. Further, by hypothesis $\alpha: X
\rightarrow \mathcal{A}$ and $X$ is a defining space with
continuous mapping for system $\mathcal{A}$. Then from item $F_3)$
of Definition 1  the existence of continuous epimorphism $\xi:
\mathcal{A} \rightarrow \mathcal{B}$ and $\xi\alpha(x) =
\beta\varphi(x)$ follows, $x \in X$, as required.

If in Definition 1 we fix the space $X$ and change the class
$\frak K$ then, in general, the system $\mathcal{A}$ also changes.
The following holds.

\textbf{Corollary 5.} \textit{Let   a topological space $X$ and a
totality of quasiatomic formulas $S$ define a system $\mathcal{A}$
in class $\frak K$ of topological algebraic systems of given
signature, and let us define a system $\mathcal{B}$  in subclass
$\frak L \subseteq \frak K$.   Then $\mathcal{B}$ is a continuous
homomorphic image of $\mathcal{A}$.}\vspace*{0.1cm}

\textbf{Proof.} By Definition 1 there   exists a mapping (7) for
which all formulas from $S$ are true in $\mathcal{A}$ and
$\mathcal{B}$, and the totalities of elements  $\{a_i\}$,
$\{b_i\}$ algebraically generate these systems by  hypothesis.
From $\mathcal{B} \in \frak L$, it follows that $\mathcal{B} \in
\frak K$, hence from statement $F_3$) it follows the existence of
such a continuous homomorphism $\varphi: \mathcal{A} \rightarrow
\mathcal{B}$ follows that $\varphi(a_i) = b_i$ ($i \in I$). As the
image $\varphi(\mathcal{A})$ of system $\mathcal{A}$ into
$\mathcal{B}$ contains the elements $b_i$, which generate
$\mathcal{B}$, then $\varphi$ is a continuous homomorphism of
$\mathcal{A}$ on $\mathcal{B}$, as required.\vspace*{0.1cm}

\textbf{Corollary 6.} \textit{Let the  sets $S_1$, $S_2$ of
defining relations on the  same  topological space $X$ defines
 respective the systems
$\mathcal{A}$ and $\mathcal{B}$ in a class $\frak K$ of
topological algebraic systems of given signature. If all formulas
from $S_1$ are consequences of the formulas of $S_2$ in $\frak K$,
particularly, if $S_1 \subseteq S_2$, then the system
$\mathcal{B}$ is an image of $\mathcal{A}$ at a continuous
homomorphism $\varphi: \mathcal{A}\rightarrow \mathcal{B}$. If the
sets of formulas $S_1$ and $S_2$ are equivalent, then the systems
$\mathcal{A}$ and $\mathcal{B}$ are topologically
isomorphic.}\vspace*{0.1cm}

\textbf{Proof.} According to Definition 1 and Lemma 2, there exist
such  continuous mappings  $x_i \rightarrow a_i (a_i \in
\mathcal{A}, i \in I), \quad x_i \rightarrow b_i (b_i \in
\mathcal{B}, i \in I)$ of space $X$  that the elements $a_i$
generate algebraically $\mathcal{A}$, the elements $b_i$ generate
algebraically $\mathcal{B}$, and the formulas from $S_1$ and $S_2$
for the specified values  $x_i$ are true  in $\mathcal{A}$ and,
respectively $\mathcal{B}$. As the formulas from $S_1$ derive
from those of $S_2$, then in $\mathcal{B}$ all formulas from $S_1$
are true. According to item $F_3)$ from here it follows that such
a continuous homomorphism $\varphi$ of  system $\mathcal{A}$ into
system $\mathcal{B}$ exists, translating $a_i$ in $b_i$ ($i \in
I$). Since the elements $b_i$ generate $\mathcal{B}$, then
$\varphi$ is the required continuous homomorphism of $\mathcal{A}$
on $\mathcal{B}$.

If the sets of formulas $S_1$ and $S_2$ are equivalent, then in
addition to  homomorphism $\varphi: \mathcal{A} \rightarrow
\mathcal{B}$ there is a  homomorphism $\psi: \mathcal{B}
\rightarrow \mathcal{A}$ sending  $b_i$ to $a_i$. Since the sets
$\{a_i \vert i \in I\}$, $\{b_i \vert i \in I\}$ generate the
systems $\mathcal{A}$, $\mathcal{B}$ and $\psi\varphi(a_i) = b_i$,
$\varphi\psi(b_i) = a_i$, then $\psi = \varphi^{-1}$. Hence the
systems $\mathcal{A}$ and $\mathcal{B}$ are topological
isomorphic, as required.\vspace*{0.1cm}

Remind that a formula $r(x_{\lambda}, \ldots, x_{\gamma})$ is
called \textit{topological consequence} of totality of formulas
$S$ in class $\frak K$  if for every continuous mapping $x_i
\rightarrow c_i$ of totality symbols $x_i$ of topological space
$X$ in any $\frak K$-system $\mathcal{C}$, at which in
$\mathcal{C}$ all formulas from $S$ are true, the formula
$r(c_{\lambda}, \ldots, c_{\gamma})$  is also true in
$\mathcal{C}$. Now we generalize \cite[Theorem V.11.1]{Mal1},
proved for algebraic systems.\vspace*{0.1cm}

\textbf{Proposition 5.} \textit{Let $\frak{K}$ be a class of
topological algebraic systems of given signature $\Omega =
\Omega_F \cup \Omega_P$. A system $\mathcal{A} = <A, \Omega> \in
\frak K$ is defined by a topological space $X$, a continuous
mapping $\sigma$ of  $X$ into the basic set $A$ of $\mathcal{A}$
and a totality $S$ of quasiatomic formulas of signature $\Omega$
of the form (6) if and only if  the following conditions are met:}

\textit{$f_1)$ all formulas from $S$ hold in $\mathcal{A}$ at
continuous mapping $\sigma: x_i \rightarrow a_i$ ($x_i \in X$,
$a_i \in \mathcal{A}$, $i \in I$);}

$f_2$) \textit{the system $\mathcal{A}$ is generated topologically
by images of elements of $X$, i.e.  algebra $< A, \Omega_F>$ does
not contain a closed subalgebra, not equal to  $<A, \Omega_F>$ and
containing $\sigma X$;}

\textit{$f_3)$ every  formula from $S$, true in $\mathcal{A}$ at
continuous mapping $\sigma: x_i \rightarrow a_i$ ($i \in I$), is a
topological consequence  of totality of formulas $S$ in  class
$\frak K$} \textit{hold true.}\vspace*{0.1cm}

\textbf{Proof.} We use the Lemma 2 without  reference.
\textit{Sufficiency.} Let the system $\mathcal{A}$ satisfy the
conditions $f_1)$, $f_2)$ and $f_3)$. Clearly,  $F_1)$, $F_2)$
follow from  $f_1)$ $f_2)$.  Let us show that $\mathcal{A}$  also
have property $F_3)$. ( Let us have a continuous mapping $\beta: X
\rightarrow \mathcal{C}$, $x_i \rightarrow c_i$ in some  $\frak
K$-system $\mathcal{C} = <C, \Omega>$, where all formulas from $S$
hold true. Let us consider on pair $<\mathcal{A}, \mathcal{C}>$
the relation $\varphi$, consisting of every possible pairs of form
$$<f(a_{\lambda}, \ldots, a_{\nu}), f(c_{\lambda}, \ldots,
c_{\nu})>, \eqno{(9)}$$ where $f$ is a term, we shall also  show,
that $\varphi$ is a continuous homomorphism of $\mathcal{A}$ in
$\mathcal{C}$. As the elements $a_i$ algebraically generate
$\mathcal{A}$, then it is possible to present any element of
system $\mathcal{A}$ as a polynomial   $f((a_{\lambda}, \ldots,
a_{\nu})$. Consequently, the left part  of relation $\varphi$ is
$\mathcal{A}$. Further, if for a some terms $f, g$ we have
$$f(a_{\lambda}, \ldots, a_{\nu}) = g(a_{\lambda}, \ldots,
a_{\nu}), \eqno{(10)}$$ then by $f_3)$ the formula
$$f(x_{\lambda}, \ldots, x_{\nu}) = g(x_{\lambda}, \ldots,
x_{\nu}) \eqno{(11)}$$ is a topological consequence of totality of
formulas $S$ in class $\frak K$. But in system $\mathcal{C}$ all
formula from $S$ at mapping $x_i \rightarrow c_i$ are true.
Therefore, the  formula (11) also is true $\mathcal{C}$, i.e.
$$f(c_{\lambda}, \ldots, c_{\nu}) = g(c_{\lambda}, \ldots,
c_{\nu}). \eqno{(12)}$$ holds in $\mathcal{C}$.

It means that $\varphi$ is a mapping from $\mathcal{A}$ on
$\mathcal{C}$. Any polynomial of topological algebraic system is a
continuous function. As $\mathcal{A}$ and $\mathcal{C}$ are a
topological algebraic systems then from \cite[Proposition
1.4.9]{Eng} it follows that $\varphi$ is a continuous mapping.

It is similarly proved that if  some relation of the form
$$p(f(a_{\lambda}, \ldots, a_{\nu}), \ldots, g(a_{\lambda}, \ldots,
a_{\nu})), \eqno{(13)}$$ hold true in $\mathcal{A}$, where $p$ is
a signature predicate, then  the relation
$$p(f(c_{\lambda}, \ldots, c_{\nu}), \ldots, g(c_{\lambda}, \ldots,
c_{\nu})). \eqno{(14)}$$ holds true in $\mathcal{C}$.

From (12) and (14) it follows that the mapping $\varphi$ is a
homomorphism. Remind, that a mapping of a topological algebraic
system on another one  is called continuous homomorphism if it is
continuous and is homomorphism in algebraic forms \cite[pag.
12]{Mal2}. Hence $\varphi: \mathcal{A} \rightarrow \mathcal{C}$ is
a continuous homomorphism.

As the pairs of form (9) belong to $\varphi$, than $<a_i, c_i> \in
\varphi$,  i.e. $\varphi(a_i) = c_i$ ($i \in I$), as
required.\vspace*{0.1cm}

\textit{Necessity.} Let the system $\mathcal{A}$ satisfy the
conditions $F_1)$ $F_2)$, $F_3)$, and let for  some terms $f, g,
\ldots, h$  the relation (10) or respectively (13)) hold  true in
$\mathcal{A}$. Then it is necessary to show that the formula (12)
(or respectively (14)) in class $\frak K$ is a topological
consequence of the totality of formulas $S$, i.e. that in every
$\frak K$-system $\mathcal{C}$ for any mapping $x_i \rightarrow
c_i$ ($c_i \in \mathcal{C}, i \in I$), where  all formulas from
$S$ hold true, the formula (12) (or respectively (14)) also holds.
At the specified assumptions, there exists a continuous
homomorphism of system $\mathcal {A}$ in the system $\mathcal {C}$
for which $\varphi(a_i) = c_i$ ($i \in I$) by $F_3)$.  But under
homomorphisms the quasiatomic formulas keep their validity, and
consequently from the validity of (10), or (13) there follows the
validity follow for (12), or respectively the formula (14). This
completes the proof of Proposition 5.

We denote by $\frak{K}$ a class of algebraic systems of signature
$\Omega$ with defining relations $S$ and let  $A_{\alpha}$
($\alpha \in \Sigma$) be a system of topological
 algebraic systems of  class  $\frak{K}$. Remind that the totality of all
functions $h(\alpha)$ defined on $\Sigma$ with  value in $\cup
A_{\alpha}$, satisfying the requirements $h(\alpha) \in
A_{\alpha}$, is \textit{the Cartesian product} $H = \prod
A_{\alpha}$.

In $H$  \textit{the Tychonoff topology} is introduced, announcing
open the Cartesian product of open subsets, chosen in a finite
number of multiplied spaces, multiplied by others spaces, as well
as unions  of these products. Further, according to \cite[section
IV.7.5]{Mal1} quasiatomic formulas multiplicative are steady. We
consider (6) and we assume by definition that $\varphi_i(h_1,
\ldots, h_{n_i}) = h$ ($h_i, h \in H$), where $h(\alpha) =
\varphi_i(h_1(\alpha), \ldots, h_{n_i}(\alpha))$, $\alpha \in
\Sigma$. Thus,  $H$ turns into a topological algebraic system of
the  same class $\frak{K}$ as the given algebraic systems
$A_{\alpha}$, which are called  \textit{(direct) topological
product of topological algebraic systems $A_{\alpha}$}
(\cite{Mal3}).

To research the  almost quasivarieties of topological algebraic
systems it will be necessary to use the following notions  from
\cite{Mal3}. Let $\frak{K}$ be a class of topological algebraic
systems of fixed signature $\Omega$ ($\frak{K}$-systems) defined
by a totality  of quasiatomic formulas of  signature $\Omega$  of
the form (6). We will examine classes of topological algebraic
systems, the topological spaces of which satisfy  the conditions:

$K_1)$ the class $\frak{K}$ is closed under  topological product;

$K_2)$ every subsystem of $\frak{K}$-system is an
$\frak{K}$-system.

Such classes, satisfying $K_1)$, $K_2)$ are the classes of
Hausdorf spaces, of regular spaces, of completely regular spaces,
with functionally separated  pairs of points etc. Moreover,  to
construct a class with the specified properties it is possible to
take a certain set of topological spaces and to consider all those
spaces, which can  be obtained from operations, taken with the
help of $K_1)$, $K_2)$. Particularly, it is possible to take as
initial spaces topological spaces with discrete topology. ++++ The
received topological spaces are regular and completely
discontinuous in the sense that for any finite set of  points
$x_1, \ldots, x_n $ of such space $X$ there is a splitting $X$  in
$n$  pairwise disjoint closed sets, containing  one of the
specified points.

To describe the  topological algebra $\mathcal{A}$, given by
defining space  $X$ and defining relations $S$ in class $\frak K$
in \cite{Mal3} the following questions are considered:

A) Is  $\mathcal{A}$ isomorphic to the abstract algebra
$\overline{\mathcal{A}}$, given in class $\frak K$ by the same
defining relations  $S$ and generating set $X$, as algebra
$\mathcal{A}$ itself?

B) Does  $\mathcal{A}$ contain a subspace  $X$ topologically, i.e.
is the mapping $\sigma$ of space $X$ in $\mathcal{A}$ a
homeomorphism?

Sufficient conditions for a positive answer to questions A) and B)
are contained in Remarks 1, 2. These remarks are transferred
literally for topological algebraic systems. It follows directly
from them that for topological algebraic systems generated by
regular and completely discontinuous space $X$ the problem A) is
solved positively for any ő $\frak K$ and $S$, and problem B) is
solved positive under the condition that the abstract  algebraic
system $\overline{\mathcal{A}}$ contain $X$. The last  assertions
hold true and for  Hausdorff space $X$ as well.

The following Theorem 5 characterizes  classes of topological
algebraic systems, any system of which can be given by defining
spaces and defining relations. Some part of proof of Theorem 5 is
similar in form  to the proofs of \cite[Theorem 1]{Mal3},
\cite[Theorem V.11.5]{Mal1}. Given the importance of Theorem 5,
and also to simplify the reading, the proof of Theorem 5 is given
fully.\vspace*{0.1cm}

\textbf{Theorem 5.} \textit{Let $\frak K$ be a class of
topological algebraic systems of fixed signature $\Omega$. Then
any system $\mathcal{A}$ of $\frak K$ can be given by defining
topological space $X$ and by defining relations $S$, i.e. by a
totality $S$ of quasiatomic formulas  of form (6), if and only if
$\frak K$ is an  almost quasivariety. i.e. $\frak K$  satisfies
the conditions:}

\textit{$Q_2)$  class  $\frak K$  is topologically  hereditary,
i.e. if $\mathcal{A} \in \frak K$ and $\mathcal{B}$ is a subsystem
of $\mathcal{A}$, then $\mathcal{B} \in \frak K$;}

\textit{$Q_3)$  class  $\frak K$ is  topologically
multiplicatively closed, i.e. the direct topological  product of
$\frak K$-systems is a topological $\frak K$-system;}

\textit{$Q_4)$ class  $\frak K$ contains an unitary
system.}\vspace*{0.1cm}

\textbf{Proof.} \textit{Necessity.} Let any  system $\mathcal{A}$
of a class  $\frak K$  be given by defining topological space and
defining relations. If we consider the topological algebraic
system $\mathcal{A} \in \frak K$ as algebraic system, then by
Lemma 2 the algebraic system $\mathcal{A}$ is  with defining
relations in the sense of \cite[Theorem V.11.2.1]{Mal1}. Then from
\cite[Lemma V.11.3.4]{Mal1} it follows that the class $\frak K$ of
algebraic systems is replica complete. Hence for a topological
algebraic system $\mathcal{A}$  considered as class of algebraic
system we use the notions of $\frak K$-replica $\mathcal{A}_{\frak
K}$ and $\frak K$-morphism.

$Q_2)$. Let $\mathcal{A}$ be a topological algebraic system of
class $\frak K$ defined by defining topological space $X$, by
continuous mapping $\sigma: X \rightarrow \mathcal{A}$ and by set
$S$ of defining relations of form (6). Let $\mathcal{B}$ be a
subsystem of $\mathcal{A}$. We denote $Y = \{x \in X \vert
\sigma(x) \in \mathcal{B}\}$. Leaving in  $\mathcal{B}$ the
topology which it has as a subset of the topological space
$\mathcal{A}$ we convert  $\mathcal{B}$ into a topological system,
where  $Y$ is reflected continuously in  $\mathcal{B}$ with the
help of narrowing  $\tau$ of $\sigma$ on $\mathcal{B}$. According
to Lemma 2 the set $\{\sigma(x) \vert x \in X\}$ algebraically
generate the system $\mathcal{A}$, then the set $\{\tau(y) \vert y
\in Y\}$ algebraically generate the system $\mathcal{B}$. The
continuity  of $\sigma: X \rightarrow \mathcal{A}$ induces the
continuity of $\tau: Y \rightarrow \mathcal{B}$ by \cite[pag.
112]{Eng}. Moreover, any element in $\mathcal{A}$ is expressed in
the form of a finite polynomial from  images elements $X$ by item
$F_2)$ of Definition 1 and Lemma 2. Then from item $F_1)$ it
follows that the formulas from $S$ hold true in system
$\mathcal{B}$.

Let now $\mathcal{B}_{\frak K}$ be the replica of system
$\mathcal{B}$ and let $\beta: \mathcal{B} \rightarrow
\mathcal{B}_{\frak K}$ be the corresponding $\frak K$-morphism. On
algebraic system $\mathcal{B}_{\frak K}$ we introduce a topology
regarding  the images of open sets of topological algebraic system
$\mathcal{B}$ as open sets. Then $\beta$ is a continuous
homomorphism, hence and $\beta\tau: Y \rightarrow
\mathcal{B}_{\frak K}$ is a continuous mapping.

From definition of replica $\mathcal{B}_{\frak K}$ it follows that
for identical embedding $\varepsilon: \mathcal{B} \rightarrow
\mathcal{A}$  there should  exist a continuous homomorphism
$\varsigma: \mathcal{B}_{\frak K} \rightarrow \mathcal{A}$,
satisfying the condition $\varepsilon = \varsigma\beta$. As
$\beta$ is an epimorphism and $\varepsilon$ is the identical
mapping of $\mathcal{B}$ on itself, then from $\varepsilon =
\varsigma\beta$ it follows that the mappings $\alpha$, $\varsigma$
mutually inverse, $\varsigma = \beta^{-1}$, therefore system
$\mathcal{B}$ is isomorphic to system $\mathcal{B}_{\frak K}$,
belonging to  class $\frak K$.

Let us use  the following statement: the quasiatomic formulas
maintain their validity in case of homomorphisms and in case of
converting into subsystem  of algebraic systems \cite[pag.
278]{Mal1}. From here it follows that if a system $\mathcal{L}$ of
class $\frak K$ of topological algebraic systems of fixed
signature satisfies the condition $f_3)$ of Proposition 5 then the
$\frak K$-morphism image of system $\mathcal{L}$, i.e. the replica
$\mathcal{L}_{\frak K}$ satisfies this condition $f_3)$. Moreover,
using the definition of replica complete class and Lemma 2 we
verify immediately that if the replica $\mathcal{L}_{\frak K}$
satisfies  the condition $f_3)$, then any subsystem of
$\mathcal{L}_{\frak K}$, which  by \cite[Theorem V.11.3.5]{Mal1}
belongs  to replica complete class $\frak K$, also satisfies
condition $f_3)$.

Now again we consider the topological algebraic system
$\mathcal{A}$ defined by  defining topological space $X$, by
continuous mapping $\sigma: X \rightarrow \mathcal{A}$ and by set
$S$ of defining relations and let $\alpha: \mathcal{A} \rightarrow
\mathcal{A}_{\frak K}$ be the corresponding $\frak K$-morphism. On
algebraic system $\mathcal{A}_{\frak K}$ we define the topology,
induced by topology of topological system $\mathcal{A}$ similarly
for system $\mathcal{B}$. Also we can show similarly  that the
mapping $\alpha\sigma: X \rightarrow \mathcal{A}_{\frak K}$ is
continuous, $\alpha \rightarrow \mathcal{A}_{\frak K}$ is  a
continuous isomorphism and any element of $\mathcal{A}_{\frak K}$
is expressed in form of finite polynomial of variables
$\alpha\sigma(x)$, $x \in X$. From here and the previous paragraph
it follows that the system $\mathcal{A}_{\frak K}$ satisfies the
conditions $f_10$ - $f_3)$ of Proposition 5 from which it follows
that the topological algebraic system $\mathcal{A}_{\frak K}$ is
defined by defining space $X$, by continuous mapping
$\alpha\sigma: X \rightarrow \mathcal{A}_{\frak K}$ and by
defining relations $S$. From Definition 1 it follows that
$\alpha^{-1}: \mathcal{A}_{\frak K} \rightarrow \mathcal{A}$ is a
continuous homomorphism, correlated with continuous mappings
$\alpha\sigma: X \rightarrow \mathcal{A}_{\frak K}$, $\sigma: X
\rightarrow \mathcal{A}$. We have shown above, that $\alpha:
\mathcal{A} \rightarrow \mathcal{A}_{\frak K}$ is a continuous
isomorphism. Hence $\alpha: \mathcal{A} \rightarrow
\mathcal{A}_{\frak K}$ is a topological isomorphism. Now the
topological isomorphism $\alpha$ induces the topological
isomorphism $\beta: \mathcal{B} \rightarrow \mathcal{B}_{\frak K}$
by narrowing the  defining space $X$ for $\mathcal{A}$ to defining
space $Y$ for $\mathcal{B}$.

We have above shown, that the replica $\mathcal{B}_{\frak K}$ is a
topological algebraic system defined by defining space $X$ and by
defining relations $S$. Then $\mathcal{B}_{\frak K}$ satisfies the
condition $f_3)$ of Proposition 5. We have also shown above that
in such a case the topological system $\mathcal{B}_{\frak K}$ also
satisfies the condition $f_3)$, as well as  $\mathcal{B} =
\beta^{-1}(\mathcal{B}_{\frak K})$, since $\beta$ is a topological
isomorphism. We have  shown above, that  $\mathcal{B}$ satisfies
the conditions $f_1)$, $f_2)$ and from Proposition 5 it follows
that the topological system $\mathcal{B}$ will be defined by
defining space $Y$, continuous mapping $\beta: Y \rightarrow
\mathcal{B}$ and defining relations $S$. The item $Q_2)$ is
proved.

$Q_3)$. Like in the previous case we consider that $\frak K$ is a
replica complete class of algebraic systems. Let $\mathcal{A} =
\prod_{i \in I}\mathcal{A}_i$ be a Cartesian product, where
$\mathcal{A}_i \in \frak K$. By \cite[Theorem V.11.3.5]{Mal1}
$\mathcal{A} \in \frak K$. By hypothesis $\mathcal{A}_i$ is a
topological algebraic system defined by defining space $X_i$,
continuous mapping $\sigma: X_i \rightarrow \mathcal{A}_i$ and
defining relations $S_i$. The system $\mathcal{A}$ together with
Tychonoff topology induced by topologies of $\mathcal{A}_i$ is
also a topological system.

Let us prove  that $\mathcal{A}$ is a topological algebraic system
given by defining space $X = \prod_{i \in I}X_i$ (Cartesian
product of spaces $X_i$), by defining relations $S = \cap_{i \in
I}S_i$ and by mapping $\sigma = \prod_{i \in I}\sigma_i$, using
the items $f_1)$ - $f_3)$ of Proposition 5. The Cartesian product
$\prod_{i \in I}\sigma_i$ of mappings $\sigma_i$ is a mapping
$\sigma: X \rightarrow \mathcal{A}$ of space $X$ in space
$\mathcal{A}$ whose each  point $x = \{x_i\} \in \prod_{i \in
I}X_i$ corresponds the point $\sigma(x) = \{\sigma_i(x_i)\} \in
\prod_{i \in I}\mathcal{A}_i$. Any mapping $\sigma_i$ is
continuous, then by \cite[Corollary 2.3.5]{Eng} the mapping
$\sigma: X \rightarrow \mathcal{A}$ is also continuous.

According to item $f_2)$ any set $\sigma_i(X_i)$ topologically
generate the system $\mathcal{A}_i$, i.e. the space
$\sigma_i(X_i)$ is dense in space $\mathcal{A}_i$. Then from
\cite[Proposition 2.3.6]{Eng} it follows that the space
$\sigma(X)$ is dense in space $\mathcal{A}$, i.e. $\sigma_i(X)$
topologically generate the system $\mathcal{A}$. Hence the
condition $f_2)$ for system $\mathcal{A}$ holds true.

We have shown above  that the algebraic system $\mathcal{A}$
belongs to replica complete class $\frak K$. Then similarly to
system $\mathcal{B}$ from item $Q_2)$ we prove that the
topological system $\mathcal{A}$ satisfies the condition $f_3)$.
Clearly,  $\mathcal{A}$ satisfies the relations $S = \cap_{i \in
I}S_i$. Hence $\mathcal{A}$ satisfies the condition $f_1)$ for
$S$. Consequently, by Proposition 4 the topological system
$\mathcal{A}$ satisfies the condition $Q_3)$.

$Q_4)$. Let us consider the unitary system $\mathcal{A}_e$ as
$\mathcal{A}$.  The unitary system satisfy any quasiatomic formula
of form (6) and any homomorphic image of unitary system can  be
only an unitary system. Then $\mathcal{A}_e \in \frak K$. Hence
the unitary system $\mathcal{A}_e$ is given by defining
topological space $X$, consisting of one point, and continuous
mapping $\sigma: X \rightarrow \mathcal{A}_e$, $\sigma(X) = e$.
The condition $Q_4)$ holds.

The necessity of Theorem 5 is proved. The  sufficiency follows
from the next statement.\vspace*{0.1cm}

\textbf{Theorem 6.} \textit{Let $\frak K$ be  an almost
quasivariety  of topological algebraic systems of given signature
$\Omega$. Then for any topological space $X$, linked with $\frak
K$ by conditions $K_1)$, $K_2)$, and for any totality $S$ of
quasiatomic formulas of signature $\Omega$ of form (6) there
exists a topological algebraic system $\mathcal{A}$ with
properties $F_1$, $F_2$, $F_3$ and is defined by these properties
univocally with accuracy of topological isomorphism over the image
of space $X$ in it.}\vspace*{0.1cm}

\textbf{Proof.} Let us consider  the property $F^{\prime}_2)$
instead of $F_2)$ by  Lemma 2. Let $X = \{x_i \vert i \in I\}$ and
let  $m = \max (|\Omega|, |I|, \aleph_0)$. According to \cite[pag.
163]{Mal1} any algebraic system of class $\frak K$, which  is
generated by a set of cardinality $\leq |I|$ has itself the
cardinality $\leq m$.

We denote by $M$ the totality of every possible pairwise
topological non-isomorphic $\frak K$-systems of cardinality $\leq
m$. According to \cite[page 280]{Mal1}, the cardinality of $M$ is
no more than a number  depending only on $m$ and $|\Omega|$.

We consider the totality of all $\frak K$-systems from $M$ which
satisfy  conditions $F_1)$, $F_2)$. We denote by
$\sigma_{\lambda}$ that continuous mapping of $X$ into
$\mathcal{A}_{\lambda}$,  which is considered  in $F_1)$, $F_2)$.
Then with respect to mappings
$$\sigma_{\lambda}: x_i \rightarrow a_i^{\lambda} \quad
(a_i^{\lambda}\in \mathcal{A}_{\lambda}, \mathcal{A}_{\lambda} \in
M) \eqno{(15)}$$  all formulas from $S$ hold in
$\mathcal{A}_{\lambda}$. The system of mappings $\sigma_{\lambda}$
is non-empty, since the almost quasivariety $\frak K$ contains the
unitary system $\mathcal{A}_e = <e, \Omega>$ (the condition
$Q_4)$) in which  not only formulas from $S$ hold at mapping $x_i
\rightarrow e$, but also  all quasiatomic formulas in general.

We consider the Cartesian  product $\mathcal{B} =
\prod\mathcal{A}_{\lambda}$ ($\lambda \in J$). As all its factors
belong to the  class $\frak K$, which is multiplicatively closed
by  condition $Q_3)$, then $\mathcal{B} \in \frak K$ Moreover,
$\mathcal{B}$ satisfies the condition $K_1)$. Hence, $\mathcal{B}$
is a topological space with topology induced by topology of space
$X$ under continuous mapping $\sigma: x_i \rightarrow
\mathcal{B}$.

By  \cite[Proposition 2.3.6]{Eng} a continuous homomorphism
$\sigma: x_i \rightarrow \mathcal{B}$ ($i \in I$) exists for the
totality of continuous homomorphisms (15). According to
\cite[Section III.7.5]{Mal1} the quasiatomic formulas are
multiplicatively stable. The formulas from $S$ are quasiatomic and
hold true in factors $\mathcal{A}_{\lambda}$ at continuous
mappings $\sigma_{\lambda}$. Hence, the formulas from $S$ hold at
continuous mappings $\sigma$.

We put $\sigma(x_i) = a_i$ and denote by $\mathcal{A}$ the
subsystem of system $\mathcal{B}$, generated in $\mathcal{B}$ by
elements $a_i$ ($i \in I$). As the formulas from $S$   do not
contain a quantifier, then from their validity  in system
$\mathcal{B}$ at mappings $\sigma: x_i \rightarrow a_i$  follows
their validity in subsystem $\mathcal{A}$. We want to show that
the generating symbols $x_i$ and the formulas from $S$ generate
the very system $\mathcal{A}$ in class $\frak K$. Indeed,
$\mathcal{A}$ is a subsystem of $\mathcal{B}$, hence by $Q_2)$
belongs to $\frak K$. The elements $a_i$ generate $\mathcal{A}$
and under continuous mapping $\sigma: x_i \rightarrow a_i$ the
formulas from $S$ hold in $\mathcal{A}$.

It remains to show  that the mapping $\sigma$ satisfies the
condition $F_3)$ from Definition 1. Indeed, let  a topological
algebraic $\frak K$-system $\mathcal{C}$  and the mapping $x_i
\rightarrow c_i$ ($c_i \in \mathcal{C}$, $i \in I$) be given. We
denote by $\mathcal{D}$ the subsystem generated in $\mathcal{C}$
by elements $c_i$. The system $\mathcal{D}$ belongs to $\frak K$
and under the mapping $x_i \rightarrow c_i$ the formulas from $S$
hold in $\mathcal{D}$. It is necessary to find a continuous
homomorphism $\gamma: \mathcal{A} \rightarrow \mathcal{D}$ at
which $\gamma(a_i) = c_i$ ($i \in I$).

The elements $c_i$ generate $\mathcal{D}$ and the cardinality of
$\mathcal{D}$ does not exceed  $|I|$. Hence with the accuracy of a
topological isomorphism the system $\mathcal{D}$ coincides with
some system $\mathcal{A}_j$ and the mapping $x_i \rightarrow c_i$
coincides with mapping  $\sigma_j: x_i \rightarrow a_i^j$, and it
is sufficient  to find a continuous homomorphism $\gamma:
\mathcal{A} \rightarrow \mathcal{A}_j$ with condition $\gamma(a_i)
= a_i^j$ ($i \in I$).

The continuous homomorphism $\sigma: x_i \rightarrow \mathcal{B}$
($i \in I$) is induced by mappings (15), hence for projection
$\pi_j$ we have $\pi_j: \mathcal{A} \rightarrow \mathcal{A}_j$,
$\pi_j\sigma_i = \sigma_i^{j}$. Hence $\pi_j$ is the required
continuous homomorphism of $\mathcal{A}$ into $\mathcal{A}_j$.

Uniqueness of system $\mathcal{A}$ follows from Corollary 3. This
completes the proof of Theorem 6.\vspace*{0.1cm}

Let $\frak K$ be a class of topological algebraic systems defined
by a totality of quasiatomic formula. In general, the class $\frak
K$ is not closed with respect to continuous homomorphic images.
However it takes place\vspace*{0.1cm}

\textbf{Corollary 7.} \textit{Let $X$ be such a topological space
 that a class $\frak K$  of topological algebraic systems of
fixed signature over $X$ satisfies  the conditions $K_1)$, $K_2)$.
Then any set of generating symbols $\{x_i \in X \vert i \in I\}$
and any totality of quasiatomic formulas from these symbols
determine the appropriate  $\frak K$-system if and only if  $\frak
K$ is  a  almost quasivariety.}

\textbf{Proof.} Sufficiency follows from Theorem 6, and necessity
follows from necessity of Theorem 5.\vspace*{0.1cm}

\textbf{Corollary 8.} \textit{An axiomatizable class $\frak K$ of
topological algebraic systems is an  almost quasivariety if and
only if  $\frak K$ is a quasivariety.}\vspace*{0.1cm}

\textbf{Proof.} Every axiomatizable class of systems $\mathcal{K}$
is ultraclosed \cite[Corollary IV.8.10]{Mal1}. Then the proof of
this  statement follows from Proposition 3 and Lemma
1.\vspace*{0.1cm}

 For the class $\frak K$ we  denote  by
symbol $\prod \frak K$ the class of all homeomorphic copies of
topological  products of topological $\frak K$-systems, by $S\frak
K$ the class of all subsystems of $\frak K$-systems and by $\frak
K_e$  the class, obtained by connection to $\frak K$ of an unitary
system.\vspace*{0.1cm}

\textbf{Proposition 6.} \textit{For a class of topological
algebraic systems $\frak K$ the class $S\prod\frak K_e$ is a
minimal  almost quasivariety, containing in itself the class
$\frak K$.}\vspace*{0.1cm}

\textbf{Proof.} The minimality  of   $S\prod\frak K_e$ is obvious,
since the almost quasivariety contains the unitary system, the
Cartesian products and the subsystems of $\frak K$-systems.
Therefore it is necessary to show  that $S\prod\frak K_e$ is an
almost quasivariety. This class contains the unitary system. If
$\mathcal{A}$ is a subsystem of the product $\prod \mathcal{A}_i$
($i \in I, \mathcal{A}_i \in \frak K_e$), then each of its
subsystems is a subsystem of the same product. Hence the class
$S\prod\frak K_e$ is hereditary. At last,  let  a sequence
$\{\mathcal{A}_i \vert i \in I\}$ of systems from $S\prod
\mathcal{A}_e$ be given. By  condition,
$$\mathcal{A}_i \subseteq \prod \mathcal{B}_{\lambda} \quad
(\lambda \in K_i) \quad (\mathcal{B}_{\lambda} \in \frak K_e; K_i
\cap K_j = \emptyset, i\neq j).$$

But in this case, according to \cite[section 1.2.5]{Mal1}, there
is an embedding
$$\prod\mathcal{A}_i \subseteq \prod_{\nu}\prod_{\lambda}
\mathcal{B}_{\lambda} \cong \prod \mathcal{B}_{\nu} \quad (\nu \in
\cup K_i),$$ from which it follows that $\prod\mathcal{A}_i \in
S\prod \frak K_e$.

From Proposition 6 it follows that the class $S\prod\frak K$ is a
minimal topological almost quasivariety which contains the class
$\frak K_e$. Let us call   $S\prod\frak K_e$ the \textit{closure
of class $\frak K$ with respect to topological almost
quasivarieties}.\vspace*{0.1cm}

\textbf{Corollary 9.} \textit{Suppose that up to topological
isomorphism the class $\frak K$ consists only of one topological
algebraic system $\mathcal{A}$. Then the minimal almost
quasivariety of topological algebraic systems, containing
$\mathcal{A}$, consists of unitary system and topological
isomorphic copies of subsystems of Cartesian degrees of system
$\mathcal{A}$.}\vspace*{0.1cm}

\textbf{Corollary 10.} \textit{If a class $\frak K$ of topological
algebraic system is axiomatizable, then the class $S\prod \frak
K_e$ is a  quasivariety.}\vspace*{0.1cm}

\textbf{Proof.} If the unitary system is joined with the
axiomatizable class $\frak K$  the axiomatizable class $\frak K_e$
is obtained. By Corollary 8 $S\prod \frak K_e$ is a quasivariety,
as required.

The class $(\mathcal{A})$ of topological isomorphic copies of
finite topological algebraic system $\mathcal{A}$ is
axiomatizable. Then from Corollary 8 it follows\vspace*{0.1cm}

\textbf{Corollary 11.} \textit{Let $\mathcal{A}$ be a finite
topological algebraic system. Then \break $S\prod (\mathcal{A})_e$
is a   quasivariety.}\vspace*{0.1cm}

Really, the class  $(\mathcal{A})$ of topological isomorphic
copies  of finite system $\mathcal{A}$ is axiomatizable. Then
Corollary 11 follows from Proposition 6.\vspace*{0.1cm}

\textbf{Remark 2.} The Definition 1,  Lemma 2,  Corollaries 3 --
6, Proposition 5 have sense according to Theorem 5 only for
 almost quasivarieties of topological algebraic systems
$\frak K$, defined by  sets of quasiatomic formulas. Since
$\mathcal{A}, \mathcal{C} \in \frak K$ in item $F_3)$ of
Definition 1, then its continuous homomorphism $\gamma:
\mathcal{A} \rightarrow \mathcal{C}$ is strong as is mentioned in
the beginning of Section 1. Moreover, all continuous homomorphisms
considered in proofs of listed above statements, and also in
Theorems 5 and 6, Proposition 6, Corollaries 7 -- 11 are
strong.\vspace*{0.1cm}

Let now $\varphi$ be a strong continuous homomorphism of a
topological algebraic system $\mathcal{A}$ of a topological almost
quasivarieties of topological algebraic systems $\frak K$  and let
$\theta$ be the congruence on $\mathcal{A}$ induced by $\varphi$.
The congruence $\theta$ will be complete and continuous. In
general, the quotient system $\mathcal{A}/\theta$ will not be a
topological algebraic system. But from Corollary 2 it
follows\vspace*{0.1cm}

\textbf{Proposition 7.} \textit{Let $\frak{K}$ be an almost
quasivariety of topological algebraic systems with permutable
congruences and let $\mathcal{A} \in \frak{K}$. Then any strong
continuous homomorphism of $\mathcal{A}$ is  open, strong and
continuous and the  quotient system $\mathcal{A}/\theta$ is a
topological algebraic system.}\vspace*{0.1cm}

\textbf{Remark 3.} Till now we examined any Hausdorff topological
spaces and we use a characteristic property of such spaces only to
prove  Lemma 2. On page 176 of \cite{Mal1} it is mentioned that
the Lemma 2 holds true for such topological space that  the
conditions $K_1)$ and $K_2)$ hold. Hence all statements of Section
2 will stay in force, if instead of Hausdorff spaces to examine a
topological spaces with properties of  capture subsystems of
topological systems and closed under direct topological
product.\vspace*{0.1cm}

\textbf{Corollary 12.} \textit{If a topological space $X$ is
defining space for some topological algebraic system, then it is
necessary for $X$ to be completely regular (Tychonoff) space.}

\section{On free systems of almost quasivariaties of topological
algebraic systems.}

This section is the continuation of Section 2. Similarly to
Section 2, in this section any considered homomorphism of
topological algebraic system is strong according to Remark 2 and
according to Remark 3 the considered topological spaces are or
Hausdorff, or such spaces that satisfy the conditions $K_1)$,
$K_2)$. In both cases these spaces are completely regular
(Tychonoff) by Corollary 12.

We shall follow the concepts from \cite[section V.12.2]{Mal1}. Let
$\frak K$ be a class of topological algebraic systems of fixed
signature $\Omega$  and let $\mathcal{A} = <A, \Omega>$ be a
topological $\frak K$-system with such a continuous mapping
$\sigma: X \rightarrow A$ that $\sigma(X)$ topologically generate
$\mathcal{A}$. A non-empty totality of elements $S$ from $A$ is
called \textit{independent in $\mathcal{A}$ with respect to class
$\frak K$} (\textit{$\frak K$-independent}) if every continuous
mapping of $S$ into any $\frak K$-system $\mathcal{B}$ can be
continued up to continuous homomorphism of
$\overline{\mathcal{S}}$ in $\mathcal{B}$, where
$\overline{\mathcal{S}}$ denotes the subsystem of $\mathcal{A}$
generated topologically by elements  $S \subseteq \sigma(X)$. The
class $\frak K$ satisfies the condition $K_2)$. Then by Lemma 2
the topological $\frak K$-system $\overline{\mathcal{S}}$ is
generated  algebraically by set $S$.\vspace*{0.1cm}

\textbf{Proposition 8.} \textit{If  pairwise distinct elements
$a_1, \ldots, a_n$ of topological algebraic system $\mathcal{A}$
are $\frak K$-independent and  $\mathcal{A}$ satisfies  some
quasiatomic relation
$$p(f_1(a_1, \ldots, a_n), \ldots, f_s(a_1, \ldots, a_n)) = \{true\},
\eqno{(16)}$$ where $p \in \{\Omega, =\}$, and $f_i(x_1, \ldots,
x_n)$ are  some terms of signature $\Omega$ not necessarily
containing all variable $x_1, \ldots, x_n \in X$, than  the
identity}
$$(\forall x_1, \ldots, x_n) p(f_1(x_1, \ldots, x_n), \ldots, f_s(x_1, \ldots, x_n)) =
\{true\} \eqno{(17)}$$  \textit{holds true in class $\frak K$.
Conversely, let for a certain set $S$ of elements of topological
algebraic system $\mathcal{A}$ from the validity of a quasiatomic
relation of form (16) for  some pairwise distinct elements $a_1,
\ldots, a_n \in S$  follows the validity  of the identity (17) in
$\frak K$. Then the set $S$ is $\frak
K$-independent.}\vspace*{0.1cm}

\textbf{Proof.} Let $a_1, \ldots, a_n \in A$ be distinct  $\frak
K$-independent elements, $S_1 = \{a_1, \ldots, a_n\}$ and let
$b_1, \ldots, b_n$ be some (not necessarily distinct) elements of
a $\frak K$-system $\mathcal{B} = <B, \Omega>$ with continuous
mapping $\eta: X \rightarrow B$. Let $b_1 = \eta(x_1), \ldots, b_n
= \eta(x_n)$. From the $\frak K$-independence of elements $a_1,
\ldots, a_n$ it follows that the mapping $\alpha: a_i \rightarrow
b_i$ ($i = 1, \ldots, n$) can be continued up to a homomorphism
$\beta: \overline{S_1} \rightarrow \mathcal{B}$. From here and
(16) it follows that
$$\{true\} = p(f_1(\beta(a_1), \ldots, \beta(a_n)), \ldots, f_s(\beta(a_1), \ldots,
\beta(a_n)) =$$ $$p(f_1(b_1, \ldots, b_n), \ldots, f_s(b_1,
\ldots, b_n)) =$$ $$p(f_1(\eta(x_1, \ldots, \eta(x_n)), \ldots,
f_s(\eta(x_1, \ldots, \eta(x_n))).$$ The elements $b_1 =
\eta(x_1), \ldots, b_n = \eta(x_n) \in B$, where $x_1, \ldots, x_n
\in X$, are arbitrary. From here it follows that (17) is an
identity of class $\frak K$.

Conversely, identity (17) follows in  $\frak K$ from any relation
of form (16) for an arbitrary airwise distinct elements $a_1,
\ldots, a_n$ of a some set $S$ of elements in $A$. We take some
mapping $\alpha$ of  $S$ in an arbitrary $\frak K$-system
$\mathcal{B}$. It is necessary to continue $\alpha$ up to
homomorphism $\overline{S}$ in $\mathcal{B}$. Above we have shown
that each element $a \in \overline {S} $ can be  presented as

$$a = f(a_1, \ldots, a_n) \quad (a_i \in S, a_i \neq a_j; i \neq j,
i, j = 1, \ldots, n),$$ where $f(x_1, \ldots, x_n)$ is a term of
signature $\Omega$. Let us introduce the mapping
$$\varphi: f(a_1, \ldots, a_n) \rightarrow f(\alpha(a_1), \ldots,
\alpha(a_n)).$$

Now show that the mapping $\varphi$ is the unique continuous
homomorphism of $\overline{S}$ into $\mathcal{B}$, i.e. that
 the validity of relation $$p(f_1(\alpha(a_1), \ldots,
\alpha(a_n), \ldots, f_s(\alpha(a_1), \ldots, \alpha(a_n))
\eqno{(18)}$$ follows from the validity of relation of form (16).

But it is obvious, since the relation (18) follows from identity
(17) at mapping  $x_i \rightarrow \alpha(a_i)$ ($i = 1, \ldots,
n$). The continuity of homomorphism $\varphi$ follows from the
definition of  identity (17) (see, for example, item $F_1)$). This
completes the proof of Proposition 8.\vspace*{0.1cm}

\textbf{Corollary 13.} \textit{If each finite subset of an
infinite set $S$ of elements of topological algebraic system
$\mathcal{A}$ is $\frak K$-independent, then  the set $S$ is also
$\frak K$-independent.}\vspace*{0.1cm}

Really, every quasiatomic relation of form (16) between the
elements of  $S$ contains only a finite subset of elements from
$S$. Then identities (17) follow from $\frak K$-independence
$\frak K$, as required.

In the definition of $\frak K$-independence  not the system
$\mathcal{A}$ itself is important, but and only subsystem
$\overline{\mathcal{S}}$, so \textit{the $\frak K$-independence of
set $S$ in system $\mathcal{A}$ is equivalent to $\frak
K$-independence of $S$ in any subsystem of system $\mathcal{A}$.}
Also, \textit{the $\frak K$-independence of any subset of  $S$
follows from $\frak K$-independence of $S$.}

Literally  repeating the proof of Theorem V.12.3 from \cite{Mal1}
(it is necessary to apply only the Proposition 8) we prove the
following Corollary 13 in which  it is underlined that a $\frak
K$-independent set of elements is independent not only with
respect to subclass $\frak L \subseteq \frak K$, but is
independent with respect to an wider class. Remind, that for any
class $\frak L$   symbols $S\frak L$, $\prod\frak L$, $H\frak L$
denote  the classes of all subsystem of $\frak L$-systems,
topological isomorphic copies of topological products and,
respectively, continuous homomorphic images.\vspace*{0.1cm}

\textbf{Corollary 14.} \textit{If a set $S$ of elements in
topological algebraic system $\mathcal A$ is independent with
respect to a class $\frak K$, then $S$ is independent also with
respect to  class $HS\prod\mathcal A$ and with respect to any
 subclass $\frak L \subseteq \frak
K$.}\vspace*{0.1cm}

The topological $\frak K$-system $\mathcal{A}$ will be called
\textit{free with respect to class $\frak K$}, if in $\mathcal{A}$
there exists a set $S$ of elements, independent from and
topologically generating  $\mathcal{A}$. The totality $S$ with
this property is called a \textit{$\frak K$-free basis of system
$\mathcal{A}$}. A system $\mathcal{A}$ is called \textit{free
system of rang $m$ in  class $\frak K$} (denoted by
$\mathcal{F}_m(\mathcal{A}$)), if $\mathcal{A} \in \frak K$ and a
$\frak K$-free basis of cardinality $m$ exists in $\mathcal{A}$.

Let $\frak K$ be a class of topological algebraic systems of fixed
signature $\Omega$. Comparing with each other the definition of
free system with respect to class $\frak K$  and the system
described by Proposition 5 given in class $\frak K$ by defining
relations (Definition 1) we ascertain.\vspace*{0.1cm}

\textbf{Lemma 3.} \textit{Let $\mathcal{A} = <A, \Omega>$ be a
$\frak K$-system, given by topological space $X$, by continuous
mapping $\sigma: X \rightarrow A$ and by defining relations
$\Delta$. Then:}

\textit{1). Any $\frak K$-independent set $S$ of elements of
system $\mathcal{A}$ generate a subsystem $\overline{S}$, free
with respect to class $\frak K$. The set $S$ is a $\frak K$-free
basis for system $\overline S$:}

\textit{2). The system $\mathcal{A}$ is free with $\frak K$-free
basis $\sigma(X)$ when and only when  the set of defining
relations $\Delta$ is empty.}\vspace*{0.1cm}

 It is worth mentioning that the notion of free system was
 introduced  earlier in \cite{Mal3} in form of item 2).

The classes of systems, consisting only from one-element systems
is called \textit{trivial}. In trivial classes only free systems
of rank 1  can exist. Thus the system $\mathcal{A}$ will be free,
if in it each formula of a form $p(a, \ldots, a)$ $(p \in
\Omega_P, a \in \mathcal{A})$ is true then and only then, when it
holds true in each system of this class.\vspace*{0.1cm}

\textbf{Theorem  7.} \textit{A non-trivial class  $\frak K$ of
topological algebraic systems of given signature $\Omega$ contains
a $\frak K$-free systems when and only when $\frak K$ is an almost
quasivariety. More specific,  when and only when the non-trivial
almost quasivarieties contain free topological algebraic systems
$\mathcal{F}_m$ of any given rank $m \geq 1$.}\vspace*{0.1cm}

\textbf{Proof.} The first part follows from item 2) of Lemma 3 and
Theorem 5. Let us prove the second part. We use  Lemma 3 and
Proposition 5. Let $\sigma: X \rightarrow A$ be the continuous
mapping from Proposition 5. Assume that $|\sigma(X)| = m$. By
Proposition 5 the totality of elements $\sigma(x_i)$, $x_i \in X$,
topologically generate the system $\mathcal{A}$ and is $\frak
K$-independent. According to Lemma 3 to complete the proof we will
show that  the cardinality of set of all elements $\sigma x_i$ is
equal to $m$, i.e. that $\sigma x_i \neq \sigma x_j$ for $i \neq
j$. But it is were the case that  $i \neq j$ and $\sigma x_i =
\sigma x_j$, then from the independence of elements $\sigma x_i$,
$\sigma x_j$ it follows that in class $\frak K$ the identity $x =
y$ would be true and class $\frak K$ would be trivial. The Theorem
7 is proved.\vspace*{0.1cm}

\textbf{Corollary 15.} \textit{A non-trivial class $\frak K$ of
topological algebraic systems of given signature $\Omega$ contains
an axiomatizable $\frak K$-free systems when and only when $\frak
K$ is a quasivariety. More specific,  when and only when  the
non-trivial quasivarieties contain an axiomatizable free
topological algebraic systems $\mathcal{F}_m$ of any given  rank
$m \geq 1$.}\vspace*{0.1cm}

Corollary 10 also Follows from Theorem 7.

A generating set of elements of a system $\mathcal{A}$ is called
\textit{minimal}, if any of its  proper subsets does not generate
$\mathcal{A}$.\vspace*{0.1cm}

\textbf{Theorem  8.} \textit{Let $\frak K$ be an almost
quasivariety  of topological algebraic systems of fixed signature
$\Sigma$. Then:}

1) \textit{an algebraic system $\mathcal{A}$, free in relation to
class $\frak K$, is  free in relation to any subclass $\mathcal{L}
\subseteq \mathcal{K}$, and  in relation to closure $HS\prod
\mathcal{K}$;}

\textit{2) all  free $\frak K$-systems of a given rank $M$ are
topologically  isomorphic among each other and  any topological
algebraic $\frak K$-system topologically generated by a set of
cardinality $m$ is an image of a continuous homomorphism of a free
systems $\mathcal{F}_m(\frak K)$ of rank $m$;}

 \textit{3) a free basis of a free system
$\mathcal{A} = \mathcal{F}_m(\mathcal{A}$) of  some class $\frak
K$ of topological algebraic system is a minimal generating set in
$\mathcal{A}$.}\vspace*{0.1cm}

\textbf{Proof.} The item 1) and the first part of item 2) follows
from Lemma 3, Corollaries  3, 14 and established after  Corollary
13 the  properties of $\frak K$-independent sets of elements. We
prove the second part of item 2).

Let $S$ be a $\frak K$-free basis of $\mathcal{A}$ with
cardinality $|S| = m$ and let $\mathcal{B}$ be a topological
algebraic $\frak K$-system topologically generated by a set $U$ of
cardinality $m$. We consider a mapping $\alpha$ of set $S$ on $U$.
By hypothesis $\alpha$ may be continued  to continuous
homomorphism $\varphi: \mathcal{A} \rightarrow \mathcal{B}$. Since
$\varphi(\mathcal{A}) \subseteq U$ and $U$ generate $\mathcal{B}$,
then $\varphi: \mathcal{A} \rightarrow \mathcal{B}$ is a
continuous homomorphism.

3). We admit on the contrary  that there is a proper part $S_1$ of
basis $S$ which generates $\mathcal{A}$. Then by
 Lemma 2 for each element $a \in S \backslash S_1$  a representation of the form
$$a = f(a_1, \ldots, a_n) \eqno{(19)}$$ exists,
where $f(x_1, \ldots, x_n)$ is  some term of signature $\Omega$
and $a_1, \ldots, a_n$ are  different elements from $S_1$. From
(19) it follows that the identity  $x = \break f(x_1, \ldots,
x_n)$ holds
 true in class $\frak K$   by  Proposition 8. As $\mathcal{A} \in \frak K$, then this
identity  also holds true  in $\mathcal{A}$. Putting here $x_1 =
a_1$, ..., $x_n = a_n$, $x = a, a_1$ we obtain $a = a_1$, that is
impossible. This completes the Proof of Theorem 8.\vspace*{0.1cm}

A free system of a class of all systems of signature $\Omega$ is
called \textit{absolutely free}. From Theorems 7, 8 it follows
that a free system of rank $m$ of arbitrary class $\frak K$ is an
image of absolutely free system of rank $m$ with respect to
continuous endomorphism.

Let $\frak K$ be an almost guasivariety of topological algebraic
systems of fixed signature. If $\mathcal{A} \in \frak K$ has an
infinite set of topological generators then all sets of such
generators have the same cardinality (proved similarly to
analogical statement for algebraic systems \cite[pag. 318]{Mal1}).
Moreover, if the set of topological generators of $\mathcal{A}$
has the cardinality $m$ then each of its topological isomorphic
system  has minimal set of topological generators of cardinality
$m$. Comparing these  statements with item 3) of Theorem 8 we get
the analog of Fujiwara Theorem \cite[pag.
318]{Mal1}.\vspace*{0.1cm}

\textbf{Corollary 16.} \textit{If in some almost quasivariety of
topological algebraic systems of fixed signature $\frak K$ a free
systems $\mathcal{F}_{m}$, $\mathcal{F}_{n}$ of various ranks $m$,
$n$ are topological isomorphic, then $m$, $n$ are
finite.}\vspace*{0.1cm}

Corollary 16 is also specified by analogy with  \cite[pag.
320]{Mal1}.\vspace*{0.1cm}

\textbf{Corollary 17.} \textit{If an almost quasivariety of
topological algebraic systems of fixed signature $\frak K$
contains a finite non-unitary system then all free systems of
various ranks of $\frak K$ are topologically
non-isomorphic.}\vspace*{0.1cm}

Let $\frak K$ be a class  of topological algebraic systems of
fixed signature and let $\mathcal{F}_n \mathcal{F}_n(\frak K)$ be
a free system of finite rank $n$. We suppose that $\frak K =
\{\mathcal{F}_n\}$. According to Lemma 3   a set $f_1, \ldots, f_n
\in \mathcal{F}_n$ which topologically generate $\mathcal{F}_n$
will called \textit{free}, if any continuous mapping $f_i
\rightarrow a_i\in \mathcal{F}_n$ ($i = 1, \ldots, n$) can be
continued until the continuous homomorphism of $\mathcal{F}_n$ in
itself. The following holds.\vspace*{0.1cm}

\textbf{Corollary 18.} \textit{In free system $\mathcal{F}_n $ of
finite rank $n$ any set of  $n$ elements which topologically
generate $\mathcal{F}_n $ is free  when and only when any
continuous epimorphism of $\mathcal{F}$ on itself is a continuous
isomorphism.}\vspace*{0.1cm}

\textbf{Proof.} Let $\varphi: \mathcal{F}_n \rightarrow
\mathcal{F}_n$ be a continuous epimorphism which is not an
isomorphism and let $\{f_1, \ldots, f_n\} \subseteq \mathcal{F}_n$
be a free generating set. Then $\{\varphi f_1, \ldots, \varphi
f_n\}$ topologically generate $\mathcal{F}_n$. If the latter was
free then by Corollary 3 there would have existed a continuous
isomorphism $\psi: \mathcal{F}_n \rightarrow \mathcal{F}_n$
turning $\{f_1, \ldots, f_n\}$  $\{\varphi f_1, \ldots, \varphi
f_n\}$ respectively. Since the homomorphisms $\varphi$, $\psi$
coincide on the generating set $\{f_1, \ldots, f_n\}$, then they
should coincide on  $\mathcal{F}_n$ as well, i.e. $\varphi =
\psi$, which contradicts the assumption.

Conversely, let  any set from $n$ elements be free in
$\mathcal{F}_n$ and let $\varphi: \mathcal{F}_n \rightarrow
\mathcal{F}_n$ be a continuous epimorphism. Then the elements
$\varphi f_1, \ldots, \varphi f_n$ topologically generate
$\mathcal{F}_n$ and therefore are free in  $\mathcal{F}_n$. That
is why there should exist a continuous  isomorphism $\psi:
\mathcal{F}_n \rightarrow \mathcal{F}_n$ turning $\{f_1, \ldots,
f_n\}$ into  $\{\varphi f_1, \ldots, \varphi f_n\}$ respectively.
Since the homomorphisms $\varphi$, $\psi$ coincide on the
generating set $\{f_1, \ldots, f_n\}$, the should coincide on
$\mathcal{F}_n$ as well, i.e. $\varphi = \psi$, as
required.\vspace*{0.1cm}

\textbf{Remark 4.} We present next \cite[pag. 178]{Mal3}. All
spaces are considered that satisfy the conditions of Lemma 2 of
present paper. Let the topological algebra $A$ be finitely
generated by a subset $X$ of $A$. We denote by $S$ the totality of
all relations between elements from $X$ in algebra $A$ and
consider the algebra $B$ with generating space $X$ and defining
relations $S$. By \cite[Remark 2]{Mal3} the algebra $B$ contains
topologically $X$ and the continuous homomorphism, which is the
continuation of identical mapping of $X$ on $X$, will be an
algebraical isomorphism between $B$ and $A$. Therefore algebra $B$
can be viewed as the same algebra $A$, only with a different
topology. This topology is called \textit{free} regarding to $X$.
From item 2) of Lemma 3 of present paper it follows that the
 identical mapping of algebra $A$ with free topology regarding
to $X$ on algebra $A$ with other topology will be continuous. This
property is \textit{characteristic} for free topology.

\section{Topological and paratopological \break quasigroups}

The theory of quasigroups and loops emerged in 1920's -- 1930's
(\cite{BM}, \cite{D}, \cite{Sus}, \cite{Mouf}) as a specification
of the notion of groupoid.

A \textit{quasigroup} (respect. \textit{right quasigroup} or
\textit{left quasigroup})  is a non-empty set $Q$ together with a
binary operation $Q \times Q \rightarrow Q; (x, y) \rightarrow xy$
such that the equations $ax = c$ and $yb = c$ (respect. $ax = b$
or $xa = b$)  have unique solutions.

A \textit{loop} (respect. right loop) $L$ is a quasigroup with a
base point, or distinguished element, $e \in L$ (respect. $f \in
L$) satisfying the equations $ea = ae = a$ (respect. $af = a$) for
all $a \in L$.

The notions of quasigroups and loops are defined with the help of
only one  binary groupoid operation. These definitions of
quasigroups and loops will be called \textit{groupoid definitions}
and the quasigroups and loops defined in such a way will be
temporarily called   \textit{$GD$-quasigroups} and respectively,
\textit{$GD$-loops.}

In order to apply the universal algebraic techniques, one must use
the universal algebraic description of quasigroups as algebras
$(Q, \cdot, /, \backslash)$ with  three binary operations,
multiplication $\cdot$, left division $\backslash$, and right
division $/$, satisfying the identities
$$(x \cdot y)/y = x, \quad
(x/y)\cdot y = x, \quad x\backslash(x\cdot y) = y,  \quad
x\cdot(x\backslash y) = y. \eqno{(20)}$$

Quasigroups defined by "equations" in this way are sometimes
referred to as \textit{equasigroups}, to distinguish this point of
view from the one often taken that quasigroups $(Q, \cdot)$ are a
special king of groupoid. If one adds a nullary (=constant)
operation $L \rightarrow L; x \rightarrow e$, $ae = a = ea \quad
\forall a \in L$, to the operations of a equasigroup $(L, \cdot,
/, \backslash)$ then $(L, \cdot, /, \backslash, e)$ is called
\textit{eloop}.

A non-empty set $G$ together with two binary operations $(\cdot)$,
$(\backslash)$ and nullary operation $f$ will be called
\textit{right eloop} with right unit $f \in G$ if
$$x\backslash(x\cdot y) = y,  \quad
x\cdot(x\backslash y) = y, \quad xf = x \eqno{(21)}$$ for all $x,
y \in G$.

Similarly, a non-empty set $G$ together with two binary operations
$(\cdot)$, $(/)$ and nullary operation $l$ will be called
\textit{left eloop} with left unit $l \in G$ if
$$(x\cdot y)/y = x, \quad (x/y)\cdot y = x, \quad lx = x \eqno{(22)}$$
for all $x, y \in G$.\vspace*{0.1cm}

\textbf{Proposition 9.} \textit{Every equasigroup (respect. eloop
either right eloop, or left eloop) $Q$ is a $GD$-quasi\-group
(respect. $GD$-loop either right $GD$-loop, or left
$GD$-loop)}.\vspace*{0.1cm}

\textbf{Proof.} For equasigroup $Q$  we consider the equation $ax
= b$, where $a, b \in Q$. Then by (20) $a \backslash (ax) = a
\backslash b$, $x = a \backslash b$. If $ab = ac$, where $a, b, c
\in Q$, then $a \backslash (ab) = a \backslash (ac)$, $b = c$.
Hence the equation $ax = b$  has an unique solution. Similarly it
is proved that the equation $xa = b$ has an unique solution. This
completes the proof of Proposition 9.\vspace*{0.1cm}

Recall that a \textit{topological group} is a group $G$ with such
a (Hausdorff) topology  that the binary operations $G \times G
\rightarrow G, (x, y) \rightarrow xy, x^{-1}y, xy^{-1}$ are
continuous. \textit{A paratopological group} $G$ is a group $G$
with such a topology  that the product maps of $G \times G$ into
$G$, $(x, y) \rightarrow xy$ are jointly continuous.

Now similarly to the definition of topological algebraic system,
and also similarly to the notions of topological group and
paratopological group we define the following.

A \textit{topological quasigroup} (respect. \textit{topological
right quasigroup} or \textit{topological left quasigroup}) is a
equasigroup (respect. right eloop or left eloop) $Q$ together with
such a topology on $Q$  that the binary operations $Q \times Q
\rightarrow Q;$ $(x, y) \rightarrow xy, x/y, x \backslash y$
(respect. $xy, x \backslash y$ or $xy, /$) are  continuous. A
\textit{topological loop} (respect. \textit{topological right
loop} or \textit{topological left loop}) is a eloop (respect.
right eloop or left eloop) $L$ which, in addition, is a
topological quasigroup (respect. topological right quasigroup or
topological left quasigroup).

 A \textit{paratopological  quasigroup} (respect.
\textit{paratopologi\-cal   loop} either
\textit{parato\-po\-logi\-cal right  quasigroup}, or
\textit{paratopological left quasigroup}, or
\textit{paratopological right  loop}, or \textit{paratopological
left loop})  will be called a $GD$-quasigroup (respect. $GD$-loop
either right $GD$-quasigroup, or left $GD$-quasigroup, or  right
$GD$-loop, or left $GD$-loop) $Q$ together with such a topology on
$Q$  that the binary operation $Q \times Q \rightarrow Q;$ $(x, y)
\rightarrow xy$ is continuous.\vspace*{0.1cm}

\textbf{Remark 5.} The proofs of the following assertions almost
literally repeat the proofs of Theorems 7, 8, 9 from \cite{Mal2}.

a). The space of topological right quasigroup (respect.
topological left quasigroup) is regular (compare with \cite[
Proposition IX.1.15] {Smith1}).

b). If a topological right loop (respect. topological left loop)
$L$ is generated by elements of some   connected neighborhood of
right unit (respect. left unit), then the space $L$ is connected.

c) Any convex topological right loop (respect. topological left
loop) is generated by elements of any neighborhood  of right unit
(respect. left unit).\vspace*{0.1cm}

\textbf{Lemma 4.} \textit{If $A$ is an open set of a topological
right quasigroup $(Q, \cdot, \backslash)$ (respect.   topological
left quasigroup $(Q, \cdot, /)$) and $g \in Q$, then the sets
$gA$, $g \backslash A$  (respect. $Ag$, $A / g$) are
open.}\vspace*{0.1cm}

\textbf{Proof.} We prove only the first assertion, the others are
proved similarly. If $b = ga$, $a \in A$, then $g \backslash b = g
\backslash (ga)$, $a = g \backslash b$ by (20). The set $A$ is
open and the  operation $(\backslash)$ is continuous by definition
of topological right quasigroup. Then there exists  such  a
neighborhood $U$ of element $b$ that $g \backslash U \subseteq A$.
By (20) $g \cdot (g \backslash U) \subseteq g \cdot A$, $U
\subseteq g$. Hence for any element $b \in gA$ the existence of
such an open set $U$  that $b \in U \subseteq gA$. Then from the
definition of open set it follows that the set $gA$ is open, as
required.\vspace*{0.1cm}

If $S$ is a set with a binary operation $S \times S \rightarrow S;
(x, y) \rightarrow x\ast y$, then for $a \in S$, the mappings
$L_{(\ast)}(a), R_{(\ast)}(a): S \rightarrow S$, defined by
$L_{(\ast)}(a)x = a\ast x$ and $R_{(\ast)}(a)x = x\ast a$, are
called the \textit{left translations}, respectively, \textit{right
translations } of $S$.

Let $(Q,\cdot,\backslash)$ be a topological right quasigroup, $a,
x \in Q$ and  $E$ denote the identical mapping on $Q$. From the
definitions (20) (or (21))  we get \break
$L_{(\cdot)}(a)L_{(\backslash)}(a)x = x$,
$L_{(\backslash)}(a)L_{(\cdot)}(a)x = x$,
$L_{(\cdot)}(a)L_{(\backslash)}(a) = E$, \break
$L_{(\backslash)}(a)L_{(\cdot)}(a) = E$. Hence the translations
$L_{(\cdot)}(a)$, $L_{(\backslash)}(a)$ are reversible. From
Proposition 9 it follows that the translations $L_{(\cdot)}(a)$,
$L_{(\backslash)}(a)$ are one-to-one. Write the last equalities
as: $L_{(\cdot)}^{(-1)}(a) = L_{(\backslash)}(a)$,
$L_{(\backslash)}^{(-1)}(a) = L_{(\cdot)}(a)$. The operations
$(\cdot)$, $(\backslash)$, which define a topological right
quasigroup, are continuous. Every translation of an algebraic
systems is a continuous  mapping \cite{Mal2}.
 Then from the last equalities it
follows that the translations $L_{(\cdot)}(a)$,
$L_{(\backslash)}(a)$ are homeomorphic. From Lemma 4 it follows
that the translations $L_{(\cdot)}(a)$, $L_{(\backslash)}(a)$ are
open.

According to definition (\cite{Smith1}), let $LM(Q)$ denote the
left multiplication group of the right quasigroup
$(Q,\cdot,\backslash)$, i.e. the group generated by all
translations $L_{(\cdot)}(a)$ of $(Q,\cdot,\backslash)$. By
 equality $L_{(\cdot)}^{(-1)}(a) = L_{(\backslash)}(a)$
the group $LM(Q)$ coincides with the semigroup generated by all
translations $L_{(\cdot)}(a)$, $L_{(\backslash)}(a)$. According to
\cite{Mal2} an element of semigroup $LM(Q)$ will be called
\textit{elementary left translation} of right quasigroup
$(Q,\cdot,\backslash)$. The product of homogeneous or open
mappings is a homogeneous or open mapping. Hence we
proved\vspace*{0.1cm}

\textbf{Theorem 9.} \textit{The left multiplication group $LM(Q)$
of a right equasigroup  $(Q,\cdot,\backslash)$ is transitive on
$Q$ and any elementary left translation in \break $LM(Q)$ is
homeomorphic and open.}\vspace*{0.1cm}

Using  definition (20) or (22) it can be proved\vspace*{0.1cm}

\textbf{Theorem 10.} \textit{The right multiplication group
$RM(Q)$ of a left equasigroup  $(Q,\cdot,/)$ is transitive on $Q$
and any elementary right translation in $RM(Q)$ is homeomorphic
and open.}\vspace*{0.1cm}

A topological space is called \textit{homogeneous} iff its group
of homeomorphisms operates transitively, i.e. iff for any pair
$(x,y)$ of points there is a homeomorphism $\varphi$ of the space
with $\varphi(x) = y$ \cite{Smith1}. Then from Theorems 9, 10 it
follows\vspace*{0.1cm}

\textbf{Corollary 19.} \textit{The underling space of any
topological right quasigroup (respect. topological left
quasigroup) is homogeneous.}\vspace*{0.1cm}

From Theorems 2, 9 and 10 it follows\vspace*{0.1cm}

\textbf{Corollary 20.} \textit{There  exists  derived ternary
operations $\alpha(x, y, z)$, \break $\beta(x, y, z)$ of a
equasigroup $(Q, \cdot, /, \backslash)$ with respect to which the
set $Q$ is a biternary system.}\vspace*{0.1cm}

We consider a biternary system with ternary operations $\alpha,
\beta$. From definition of ternary system it follows that $x =
\beta(\alpha(x,x,y),x,y) = \beta(y,x,y)$, $x = \beta(y,x,y)$. We
put $\Psi(x,y,z) = \beta((\alpha(x, y, a), z, a)$. Using the
definition of ternary system and $x = \beta(y,x,y)$ it is easy to
see that the polynomial $\Psi(x,y,z)$ satisfies the identities
(5).

We also  mention that for equasigroups $(Q, \cdot, /, \backslash)$
the identities (5) are satisfied by polynomial $W(x,y,z) =
(x(a/y))\backslash (a/z)$, where $a$ is a fixed element in $Q$,
(Mal'cev),  and by polynomial $W(x,y,z) = x \cdot (y \backslash
z)$ (respect. $W(x,y,z) = (x / y) \cdot z$) for right eloops $(L,
\cdot, \backslash, f)$ (repect. left eloops $(L, \cdot, /, l)$)
(it is verified directly using the definitions of right and left
eloops).\vspace*{0.1cm}

Then from Theorem 1 it follows\vspace*{0.1cm}

\textbf{Proposition 10.} \textit{Let $(Q, \cdot, \backslash, /)$
(respect. $(Q, \cdot, \backslash, f)$ or $(Q, \cdot, /, l)$) be a
equasigroup (respect. right eloop or left eloop). Then the
polynomial \break $W(x,y,z)  = (x(a/y))\backslash (a/z)$, where
$a$ is a fixed element in $Q$ (respect. \break $W(x,y,z) = x \cdot
(y \backslash z)$ or $W(x,y,z) = (x / y) \cdot z$) of equasigroup
$(Q, \cdot, \backslash, /)$ (respect. right eloops $(Q, \cdot,
\backslash, f)$ or  left eloops $(Q, \cdot, /, l)$) satisfies the
identities (5) and the congruences of equasigroups (respect. right
eloops or left eloops) are permutable.}\vspace*{0.1cm}

According to (20) - (22) the homomorphic image of a left
equasigroup (respect. right equasigroup either left eloop, or
right eloop) is a left equasigroup (respect. right equasigroup
either left eloop, or right eloop).\vspace*{0.1cm}

\textbf{Corollary 21.} \textit{Let $\varphi: Q \rightarrow
\overline Q$ be a homomorphism of a topological quasigroup
(respect. topological left loop or topological right loop),  let
$\theta$ be the congruence corresponding to the homomorphism
$\varphi$ and let for element $a \in Q$ $K_a = \{x \in Q \vert x
\theta a\}$ denote the class of congruence $\theta$. Then:}

\textit{1) $\varphi$ is an open (closed) mapping;}

\textit{2) for every open (closed) set $A \subseteq Q$ the union
of all  congruence classes intersecting $A$  is open (closed) in
$Q$;}

\textit{3) if $\varphi$ is a homomorphism $''$on$''$, then
$\varphi$ is a quotient homomorphism;}

\textit{4) the isomorphism $Q/\theta \rightarrow \overline Q$ is a
homeomorphism.}

\textit{6) if $a \cdot b = c$ (or $a \backslash b = c$) for
elements $a, b, c$ in topological right quasigroup $(Q, \cdot,
\backslash)$ then $a\cdot K_b = K_c$ (or $a \backslash K_b =
K_c$), $K_a \cdot b \subseteq K_c$ (or $K_a \backslash b \subseteq
K_c$) and the spaces $K_b, K_c$ are homeomorphic.}\vspace*{0.1cm}

\textbf{Proof.} The items 1) - 4) follow from Proposition 4,
Theorem 10 in \cite{Mal2} and Propositions 2.4.3, 2.4.4 in
\cite{Eng}.

5). We prove the assertion only for operation $(\cdot)$ as for
$(\backslash)$ it is analogical. Clearly, if $a \in K_b$ then $K_a
= K_b$ and conversely. Then, obviously, $K_a \cdot b \subseteq
K_c$.

We prove that $a\cdot K_b = K_c$. Indeed, if $y \in K_{ab}$, then
$y \theta ab$. Let $y = az$, then by (21) $az \theta ab$,
$a\backslash (az) \theta a\backslash(ab)$, $z \theta b$ and $y =
az \in aK_b$. Hence $K_{ab} \subseteq aK_b$. Let now $x \in aK_b$.
Then $x = ab_1$, where $b_1 \in K_b$. Hence $b_1 \theta b$, $ab_1
\theta ab$, $x \in K_{ab}$, $aK_b \subseteq K_{ab}$. Consequently,
$aK_b = K_{ab}$, as required.\vspace*{0.1cm}

\textbf{Corollary 22.} \textit{Let $\varphi: Q \rightarrow
\overline Q$ be a continuous homomorphism of a locally compact
left quasigroups (respect. locally compact  right quasigroups)
$Q$, $\overline Q$, satisfying the second axiom of countability,
and let $\theta$ be the congruence corresponding to the
homomorphism $\varphi$. Then the similar assertions to 1) - 5) of
Corollary 21 hold.}\vspace*{0.1cm}

It follows from Theorems 9, 10, Theorem 12 in \cite{Mal2} and
\cite[Propositions 2.4.3, 2.4.4]{Eng}.\vspace*{0.1cm}

\textbf{Proposition 11.} \textit{The underling space of any
topological quasigroup  is rectifiable.}\vspace*{0.1cm}

\textbf{Proof.} Let $(Q, \cdot, /, \backslash)$ be a topological
guasigroup. The right translations $R_{(\cdot)}(x)$, $R_{(/)}(x)$,
$x \in Q$, are homeomorphisms of the space $Q$ according to
Theorem 10 . We use (20). For $x, y \in Q$ we have
$R_{(\cdot)}^{-1}(x) = R_{(/)}(x)$, $R_{(\cdot)}(x\backslash y)(x)
= y$. Let $a$ be a fixed element in $Q$. We put $\Phi(x, y) = (x,
R_{(\cdot)}(x \backslash a)(y))$, $\Phi^{-1}(x, y) \break = (x,
R_{(/)}(x \backslash a)(y))$. It is easy to see that $\Phi$ is a
homeomorphism of $Q \times Q$ onto itself and $\Phi(x, x) = (x,
a)$. Then from the definition it follows that the space $Q$ is
rectifiable, as required.\vspace*{0.1cm}

For compact quasigroups the Proposition 11 is proved in \cite[pag.
1112]{Usp2}.

The definitions of Mal'cev space and of Mal'cev operation can be
found in the assertion $^{(8)}$ below. In \cite{Usp1} (see, also,
\cite{Usp2}) it is proved that compact Mal'cev spaces are Dugundji
(Theorem 1) and if $X$ is a countable Mal'cev space, then the
Stone-$\check{C}$ech compactification $\beta X$ of $X$ is Dugundji
(Theorem 2). It is known that if $\beta X$ is Dugundji, then $X$
is pseudocompact. Conversely, if $X$ is a pseudocompact
rectifiable space and satisfies the Suslin property, then $\beta
X$ is compact Dugundji (Theorem 3). Any countably Mal'cev space
satisfies the Suslin property (Theorem 8 in \cite{Usp2}).

Except enumerated facts the papers \cite{Usp1}, \cite{Usp2},
\cite{Gul} contain many other results relatted to  rectifiable
spaces, Mal'cev spaces which with the help of the Theorem 10 and
the Propositions 10, 12 are  easy transferred  on quasigroup
structures. Particularly, from the listed  results and
\cite[Corollary 2.2, Theorem 3.2]{Gul} it follows\vspace*{0.1cm}

\textbf{Proposition 12.} \textit{Let $Q$ be the topological space
of a topological quasigroup either topological left loop, or
topological right loop. Then:}

\textit{1) if $Q$ is a compact space then $Q$ is Dugundji;}

\textit{2) if $Q$ is a countably  space, then the
Stone-$\check{C}$ech compactification $\beta Q$ of $Q$ is
Dugundji;}

\textit{3)  if $Q$ is a pseudocompact  space and satisfies the
Suslin property, then $\beta X$ is compact Dugundji;}

\textit{4) the space $Q$ is regular;}

\textit{5) if $Q$ is a countable  space, then $Q$ is
metrizable.}\vspace*{0.1cm}

In \cite[Proposition IX.1.15 (Salzmann 1757)]{Smith1} it is proved
that the underlying space of a topological quasigroup is
regular.\vspace*{0.1cm}

\textbf{Remark 6.} In \cite[Corollary 9]{Usp2} it is proved that
any compact quasigroup is Dugundji (see item 1) of Proposition
11). It is also market  that this result is due to Choban,
published in: Choban M. M. The structure of locally compact
algebras//Bacu internat. topol. conf.  (3—9 of October 1787 ă.):
Tez. Ch. 2. Bacu: Communist, 1787. P. 334. 56. In \cite[Corollary
6]{Usp2}  it is proved that any  compact Mal'cev space is
Dugundji. It is marked (pag. 1094, 1109) that this result is  due
to Choban, but is not indicated where are they published.

Further, in \cite[pag. 110]{Gul} it is mentioned  that the Theorem
3.2 (stated above) is a generalization of Choban's theorem, who
proved under the same conditions that $X$ is a Moore space. Again,
it is not mentioned where this result is published. Other cases,
where in detriment to authors results are unmerited assigned to
Choban, will be analyzed  in assertion $^{(9)}$. The item 4) of
Proposition 12 is similar to \cite[Theorem 7]{Mal2}.

According to Theorems 9, 10  the left and right translations in
any topological quasigroup are homeomorphisms. But there are loops
defined on topological spaces in which multiplication is jointly
continuous and translations fail to be homeomorphisms. Hence such
loops are paratopological loops, and  not topological
loops.\vspace*{0.1cm}

\textbf{Example.} (\cite[Example IX.1.5]{Smith1}). Let $L$ be the
real Hilbert space $l^2$ of all square summate real sequences and
define for $x = (x_n)_{n \in N}$ a product $$xy = ((x -n + y_n)/(1
+ (n - 1)|x_1||y_1|))_{n \in N}.$$ Then $L$ is a commutative loop
and $L \times L \rightarrow L; \quad (x,y)\mapsto xy$ is
continuous, but translation with $(1, 0, 0, \ldots )$ is not open.
Hence $L$ is not a topological loop.

If $L_1$ is the subset of all sequences $x = (x_n)_{n \in N}$ with
$x_n = 0$ for $n > 1$ and $L_2$ is the subset of all $x$ with $x_1
= 0$, then $L_2$ is a normal subloop with $L/L_2 \cong L_1 \cong
R$, and $L_2$ is an abelian topological group isomorphic to the
additive group of $l^2$. Furthermore, $L = L_1L_2$ and $L_1
\bigcap L_2= \{0\}$.\vspace*{0.1cm}

\textbf{Remark 7.} In some works, particularly \cite{Usp2}, pag.
1095] the notion of topological quasigroup is defined as a space
$X$ with such continuous operation $(x, y) \rightarrow xy$  that
the equations $ax = b$, $ya = b$ have  unique solutions for any
$a, b \in X$. According to \cite{Mal2} this definition is not
correct. This is  the definition of paratopological quasigroup. In
\cite{Usp2} there are examined only  compact paratopological
quasigroups. According to the following assertion (\cite{Smith1},
Proposition IX.1.6) this does not influence the truth of other
results of the given work. (This fact for compact spaces is
mentioned in \cite{Usp2}). It is hold.

\textit{Let $(Q, \cdot)$ be a paratopological quasigroup. Then $Q$
is a topological quasigroup if at least one of the following
conditions is satisfied: (i) the space $Q$ is compact,  (ii) the
space $Q$ is locally compact and locally connected, and the
mappings $Q \rightarrow Q$; $x \mapsto ax, xa$ are open for all $q
\in Q$.}\vspace*{0.1cm}

\textbf{Lemma 5.} \textit{Let $(Q, \cdot, \backslash, /)$ be a
paratopological quasigroup and let the set $aH$ (respect. $Ha$) be
open for every open set $H$ of $Q$ and every $a \in Q$. Then the
operation $(\backslash)$ (respect. $(/)$) is
continuous.}\vspace*{0.1cm}

\textbf{Proof.}  Let $M$ be an arbitrary open set of $Q$. Then the
set $K = \cup_{m \in M}(m\cdot H)$ is open. Let $k \in K$, $m \in
M$, $h \in H$, i.e. $k = m \cdot h$. According to (20) we have $m
\backslash k = m \backslash (m \cdot h) = h$, $k / h = (m \cdot h)
/ h = m$, i.e. $m \backslash k = h$, $k / h = m$. For the first
equality we have: for all open set $H$ such  an open sets $M, K$
exists  that $M \backslash K = H$. This means that the operation
$(\backslash)$ is continuous. Similarly, the operation $(/)$ is
continuous, as required.\vspace*{0.1cm}

\textbf{Corollary 23.} \textit{If $G$ is a paratopological group,
then from conditions $a \in G$, $H$ is an open set it does not
result that the sets $a\cdot H$, $Ha$ are open.}\vspace*{0.1cm}

\textbf{Proposition 13.} \textit{Let $\mathcal{P}\mathcal{K}$
denote the class of all paratopological  loops. Then:}

\textit{1) the class $\mathcal{P}\mathcal{K}$ contains
paratopological loops which are not  topological algebraic
systems;}

\textit{2) the class $\mathcal{P}\mathcal{K}$ is not a variety: in
such a case it is not possible to consider  free paratopological
loops of $\mathcal{P}\mathcal{K}$;}

\textit{3) the class $\mathcal{P}\mathcal{K}$ contains
paratopological loops with non-permutable congruences;}

\textit{4) every paratopological  loop admits such a structure of
ternary operation  \break $[x,y,z]$  that it satisfies the
identities $[x,x,y] = y$, $[y,x,x] = y$. But  in the class
$\mathcal{P}\mathcal{K}$ there are paratopological  loops $(L,
\cdot)$, which do not contain polynomials $\Psi(x,y,z)$,
satisfying identities  (5), i.e. the $GD$-loops $(L, \cdot)$ do
not satisfy the Theorem 1;}

\textit{5) on  paratopological  loops $(L, \cdot)$ from item 3) it
is not possible to define a structure of biternary
algebra.}\vspace*{0.1cm}

\textbf{Proof.} 1) The definition of $GD$-loop does not satisfy
the definition of algebraic system.

Further, we use $GD$-loops constructed in \cite{Cowell},
\cite{Thurston}. But this $GD$-loops may be  considered as
paratopological  loops (for example, with discrete topology or
anti-discrete topology).

 2) By Birkhoff Theorem, a class of algebraic
systems is a variety  if and only if this class
 is closed with respect to  taking   subsystems,
 direct products and  homomorphic images (see, the Proposition
2).  \cite{Cowell}  proves that for any groupoid with division
(not necessary eloop) $M$ there exists a such $GD$-loop $F$
homomorphic to $M$  that the corresponding congruence is
permutable with all congruences on $F$. Hence, the class
$\mathcal{P}\mathcal{K}$ is not closed with respect to homomorphic
images. Consequently, the class $\mathcal{P}\mathcal{K}$ is not a
variety. Further, it is known that any variety  is characterized
by its free objects.

3)  \cite{Thurston} offers  an example of an $GD$-quasigroup with
a pair of non-permu\-ta\-ble congruences and this method can be
used to construct a similar example for $GD$-loops \cite{Cowell}.

4) Let $(Q, \cdot, e)$ be a right $GD$-loop with the right unit
$e$. Then the equality $a\cdot x = b$ has an unique solution $x_0$
for any $a, b \in Q$. We denote $A(a,b) = x_0$. Then $x \cdot
A(x,y) = y$ for any $x, y \in Q$. Further, from $x \cdot A(x,x) =
x$, $x \cdot e = x$  it follows that $A(x,x) = e$ for all $x \in
Q$. We put $[x,y,z] = x \cdot A(y,z)$. Then $[x,x,y] = y$, $y,x,x]
= y \cdot A(x,x) = y \cdot e = y$.

In the proof of Theorem 1 the free algebras of varieties are used
essentially. But it is impossible to prove the item 4) directly
using Theorem 1 since  the item 2) does not permit it.

Let $(L, \cdot)$ be the $GD$-loop satisfying the item 3). We
assume that the loop $(L, \cdot)$ contains a polynomials
$\Psi(x,y,z)$, satisfying identities $\Psi(x,x,z) = z, \quad
\Psi(x,z,z) = z$.  Let $\theta_1, \theta_2$ be such congruences of
$(L, \cdot)$  that $\theta_1 \theta_2 \neq \theta_2 \theta_1$.
Then such elements $a, b \in L$ exist  that $a \equiv b(\theta_1
\theta_2)$, but $a \not \equiv b (\theta_2 \theta_1)$. From $a
\theta_1 \theta_2 b$ it follows that  $a \theta_1 c$, $c \theta_2
b$ for some $c \in L$. Then $$\Psi(a, c, b) \equiv \Psi(a, a,
b)(\theta_1), \quad \Psi(a, c, b) \equiv \Psi(a, c, c)(\theta_2)$$
or $$\Psi(a, c, b) \equiv b(\theta_1), \quad \Psi(a, c, b) \equiv
a(\theta_2).$$ From here it follows that $\theta_1 \theta_2 =
\theta_2 \theta_1$. We get a contradiction, as $\theta_1 \theta_2
\neq \theta_2 \theta_1$. Hence  item 4) holds.

5) Let us define a ternary system $(L, \alpha, \beta)$ on
$GD$-loop $(L, \cdot)$ from item 3). Then the polynomial
$\Psi(x,y,z) = \beta(\alpha(x,y,a), z, a)$, where $a$ is a fixed
element in $L$, satisfies the identities (5). Contradiction with
item 3). This completes the proof of Proposition
13.\vspace*{0.1cm}

Now, let us denote by $\mathcal{T}\mathcal{K}$ the class of all
topological quasigroups (respect. topological right quasigroups
either topological left quasigroups, or topological loops, or
topological right loops, or topological left loops) and  denote by
$\mathcal{P}\mathcal{K}$ the class of all paratopological
quasigroups (respect. paratopological right quasigroups or
paratopological left quasigroups, or paratopological loops, or
paratopological right loops, or paratopological left loops). The
translations with respect to the basic operations of topological
algebraic systems are continuous. Clearly, the Proposition 13
holds for any class $\mathcal{P}\mathcal{K}$. Then from
Propositions 17, 13 and Example it follows.\vspace*{0.1cm}

\textbf{Corollary 24.} \textit{The  inclusion
$\mathcal{T}\mathcal{K} \subset \mathcal{P}\mathcal{K}$ is
strong.}\vspace*{0.1cm}

 \section{Correlation  of rectifiable spaces and  topological right loops}

We shall use the notion of rectifiable space from \cite{Gul},
given at the beginning of our work. Other information about
rectifiable space are given in the next section in the analyzes of
paper \cite{Chob1}, \cite{Chob2}, \cite{AC1}.

Let $A$ be a topological space. We denote by $\pi_1, \pi_2 : A
\times A \mapsto A$ the projections of first  and second
coordinates of Cartesian product $A \times A$. A mapping $\Phi: A
\times A \mapsto A \times A$ will be called \textit{l-mapping} if
$\pi_2 \circ \Phi = \pi_2$ and will called \textit{r-mapping} if
$\pi_1 \circ \Phi = \pi_1$. Clearly, if $\Phi: A \times A
\rightarrow A \times A$ is a homeomorphism (respect. one-to-one
mapping) of $A \times A$ onto itself, then the inverse mapping
$\Phi^{-1}$ is a homeomorphism (respect. one-to-one mapping) of $A
\times A$ onto itself. We will follow the  scheme of the proof of
Proposition 2.1 from \cite{Gul}.\vspace*{0.1cm}

\textbf{Theorem 11.} \textit{Let $A$ be a topological space and
let $\Phi$,$\Psi: A \times A \mapsto A \times A$ be a surjective
homeomorphisms. We denote $(\cdot) = \pi_2 \circ \Psi$, $(\ast) =
\pi_1\circ \Phi$,  $(\backslash) = \pi_2 \circ \Psi^{-1}$, $(/) =
\pi_1 \circ \Phi^{-1}$. Then:}

\textit{1) a space $A$ is a topological left  quasigroup with
respect to multiplication $(\cdot)$ and right  division $(/)$ iff
$\Phi$ is a $l$-mapping;}

\textit{2) a space $A$ is a topological right quasigroup with
respect to multiplication $(\cdot)$ and left division
$(\backslash)$ iff $\Psi$ is a $r$-mapping;}

\textit{3) a space $A$ is a topological  quasigroup with respect
to multiplication $(\cdot)$, left division $(\backslash)$ and
right division $(/)$  iff $\Phi$ is a $l$-mapping, $\Psi$ is a
$r$-mapping and $\pi_2 \circ \Psi^{-1} =  \pi_1\circ \Phi^{-1}$;}

\textit{4) a topological left quasigroup $(A, /)$, defined by
homeomorphism $\Phi$,  is a topological left loop with a left unit
element $f \in A$ if for any $x \in A$ the equality $\Phi^{-1}(x,
x) = (f, x)$ is fulfilled;}

\textit{5) a topological  right quasigroup $(A, \backslash)$,
defined by homeomorphism $\Psi$, is a topological right loop with
a right unit element $e \in A$ if for any $x \in A$ the equality
$\Psi^{-1}(x, x) = (x, e)$ is fulfilled;}

\textit{6) a topological  quasigroup $(A, /, \backslash)$, defined
by homeomorphisms $\Phi$, $\Psi$, is a topological loop with a
unit element $e \in A$ if for any $x \in A$ the equalities
$\Psi^{-1}(x, x) = (x, e)$ and $\Phi^{-1}(x, x) = (e, x)$ are
fulfilled.}\vspace*{0.1cm}

\textbf{Proof.}  1), 2). The mappings $\Phi, \Phi^{-1}$, $\Psi,
\Psi^{-1}$ on Cartesian product are homeomorphisms. Then they are
continuous. Particularly, according to Remark 1  the projections
to the coordinates $(\cdot) = \pi_2 \circ \Psi$, $(\ast) =
\pi_1\circ \Phi$, $(\backslash) = \pi_2 \circ \Psi^{-1}$, $(/) =
\pi_1 \circ \Phi^{-1}$, are continuous.

Let $\pi_2 \circ \Phi = \pi_2$ for item 1) or $\pi_1 \circ \Psi =
\pi_1$ for item 2). Then $(y \ast x) / x) = \pi_1 \circ
\Phi^{-1}(\pi_1 \circ \Phi(y, x), x) = \pi_1 \circ \Phi^{-1} \circ
\Phi(y, x) = \pi_1(y, x) = y$ and similarly for the other three
identities $(y / x) \ast x = y$, $x \cdot (x \backslash y) = y$,
$x \backslash (x \cdot y) = y$. Hence by definitions, $(A, \ast,
/)$ is a topological left quasigroup and $(A, \cdot, \backslash)$
is a topological right quasigroup.

Conversely, let $(A, \ast, /)$ be a topological left quasigroup,
and $(A, \cdot, \backslash)$ be a topological right quasigroup. We
put $\Phi(y, x) = (y \ast x, x)$ $\Phi^{-1}(y, x) = (y / x, x)$,
for left quasigroup and $\Psi(x, y) = (x, x\cdot y)$,
$\Psi^{-1}(x, y) = (x, x\backslash y)$ for right quasigroup.
Further we shall consider only the  right quasigroup as for the
left quasigroup reasonings are similar.

We have $\Psi \circ \Psi^{-1}(x,y) = \Psi(x, x\backslash y) = (x,
x(x \backslash y)) = (x, y)$ and $\Psi^{-1} \circ \Psi(x, y) = (x,
y)$. Hence $\Psi \circ \Psi^{-1} = \Psi^{-1} \circ \Psi = E$,
where $E$ is the identical mapping. Then from here and the
definition of right quasigroup it follows that $\Psi,  \Psi^{-1}$
are an one-to-one mappings of $A \times A$ onto itself.

For any fixed element $x \in A$ we have $\Psi(\{x\} \times A) =
\{x\} \times A$. By hypothesis the operation $(\cdot)$ is
continuous, hence the mapping $y \rightarrow \pi \circ \Psi(x, y)$
is continuous. Then from the definition of topology of Cartesian
product it follows that the mapping $\Psi: A \times A \rightarrow
A \times A$ is continuous. It is similarly proved  that the
mappings $\Psi^{-1}, \Phi, \Phi^{-1}$ are continuous.

3) From $\pi_2 \circ \Psi =  \pi_1\circ \Phi$ it follows that
$(\cdot) = (\ast)$. Then the item 3) follows from items 1), 2).

4), 5), 6). From $(/) = \pi_1 \circ \Phi^{-1}$ and $\Phi^{-1}(x,
x) = (f, x)$ it follows that $x / x = f$. Then by (3) $(x / x)
\ast x = f \ast x$, $x = f \ast x$. Hence $f$ is a left unit
element of topological left loop $A, \ast, /)$. Similarly it is
proved that from $(\backslash) = \pi_2 \circ \Psi^{-1}$ and
$\Psi^{-1}(x, x) = (x, e)$ it follows that $e$ is a right unit
element of right loop $(A, \cdot, \backslash)$. This completes the
proof of Theorem 11.

The topological structures, considered in   Theorem 11, are
topological algebraic systems. Hence for its investigation  the
powerful methods of theory of topological algebraic systems
\cite{Mal2}, \cite{Mal3} will be applied, as well as of theory of
topological quasigroups and loops  \cite{Smith1}. However, if  in
Theorem 11 the condition that the mappings $\Phi, \Psi$ are
homeomorphisms is weakened to the condition that the mappings
$\Phi, \Psi$ are only  continuous mappings, then the obtained
topological structure is no longer a topological algebraic
system.\vspace*{0.1cm}

\textbf{Proposition 14.} \textit{Let $A$ be a topological space
and let $\Phi$,$\Psi: A \times A \mapsto A \times A$ be continuous
one-to-one surjective mappings. We denote $(\cdot) = \pi_2 \circ
\Psi$, $(\ast) = \pi_1\circ \Phi$, $(\backslash) = \pi_2 \circ
\Psi^{-1}$, $(/) = \pi_1 \circ \Phi^{-1}$. Then:}

\textit{1)  $A$ is a paratopological left  quasigroup with respect
to multiplication $(\cdot)$ and right  division $(/)$ iff $\Phi$
is a $l$-mapping;}

\textit{2)  $A$ is a paratopological right quasigroup with respect
to multiplication $(\cdot)$ and left division $(\backslash)$ iff
$\Psi$ is a $r$-mapping;}

\textit{3)  $A$ is a paratopological  quasigroup with respect to
multiplication $(\cdot)$, left division $(\backslash)$ and right
division $(/)$  iff $\Phi$ is a $l$-mapping, $\Psi$ is a
$r$-mapping and $\pi_2 \circ \Psi^{-1} =  \pi_1\circ \Phi^{-1}$;}

\textit{4) a paratopological left quasigroup $(A, /)$, defined by
homeomorphism $\Phi$,  is a paratopological left loop with a left
unit element $f \in A$ if for any $x \in A$ the equality
$\Phi^{-1}(x, x) = (f, x)$ is fulfilled;}

\textit{5) a paratopological  right quasigroup $(A, \backslash)$,
defined by homeomorphism $\Psi$, is a paratopological right loop
with right unit element $e \in A$ if for any $x \in A$ the
equality $\Psi^{-1}(x, x) = (x, e)$ is fulfilled;}

\textit{6) a paratopological  quasigroup $(A, /, \backslash)$,
defined by homeomorphisms $\Phi$, $\Psi$, is a paratopological
loop with a unit element $e \in A$ if for any $x \in A$ the
equalities $\Psi^{-1}(x, x) = (x, e)$ and $\Phi^{-1}(x, x) = (e,
x)$ are fulfilled.}\vspace*{0.1cm}

\textbf{Proof.} The proof of Proposition 14 is contained in the
proof of Theorem 11.\vspace*{0.1cm}

\textbf{Remark 8.} The topological structures considered in item
5) of Theorem 11 coincide with the notion of rectifiable space,
introduced in \cite{Usp1},  \cite{Gul}. This structure also
coincides with the notion of homogeneous space, introduced in
\cite{Chob1}, \cite{Chob2}.\vspace*{0.1cm}

Let $(X, \mathcal{T})$ be a topological space. A subset of $X$ is
$\mathcal{T}$-open if it is an union of elements of $\mathcal{T}$.
We define on the set Let  $Y$ be a  subset of $(X, \mathcal{T})$.
We define on the set $Y$  the topology $\mathcal{A}$ induced by
the topology $\mathcal{T}$ as the family of intersections of
elements of $\mathcal{T}$ with the set $Y$. The space $(Y,
\mathcal{A})$ is called subspace of space $(X, \mathcal{T})$.
Hence the set $U$ belongs to the induced topology $\mathcal{A}$ if
and only if $U = V \cap Y$ for some $\mathcal{T}$-open set $V$. If
the set $Y$ is $\mathcal{T}$-open then every $\mathcal{A}$-open
set of $Y$ is $\mathcal{T}$-open as an intersection of
$\mathcal{T}$-open set with $Y$.\vspace*{0.1cm}

\textbf{Lemma 6}. \textit{Let $(Y, \mathcal{A})$ be a subspace of
a  topological space $(X, \mathcal{T})$ and let the set $Y$ be
$\mathcal{T}$-open. Then a subset of $Y$ is $\mathcal{T}$-open if
and only if it is $\mathcal{A}$-open.}\vspace*{0.1cm}

\textbf{Corollary 25.} \textit{Let a subspace $A$ of a topological
space $(X, \mathcal{T})$ be a rectifiable space. Then  $A$  is
open.}\vspace*{0.1cm}

\textbf{Proof.} According to item 5) of Theorem 11 the rectifiable
space $A$ admits a structure of topological right loop with right
unit $(A, \cdot, \backslash, e)$ with topology $\mathcal{A}$
induced by topology $\mathcal{T}$.

Let $U(e)$ be a neighborhood of the unit $e$ in $(X,
\mathcal{T})$. By Lemma 6 $U^{\prime}(e) = U(e) \cap A$ is open in
$A$. Let $a \in A$. By Lemma 4 $aU^{\prime}(e)$ is open in $(A,
\mathcal{A}$. Then by Lemma 6 $aU^{\prime}(e)$ is open in $(A,
\mathcal{T}$. As $e \in U^{\prime}(e)$, then $a \in
U^{\prime}(e)$.Hence $a \in aU^{\prime}(e) \subseteq X$. This
means that  $A$  is open, as required.\vspace*{0.1cm}

\textbf{Proposition 15}. \textit{If a non-empty subspace $A$ of a
compact space $X$ is a rectifiable space, then $A$ is a compact
space. Moreover, the set $A$ is open,  regular, Hausdorff and,
consequently, normal.}\vspace*{0.1cm}

\textbf{Proof.} By Corollary 25 the set $A$ is open. Then the set
$B = X \backslash A$ is closed and $A \cap B = \emptyset$. Any
closed subspace of a compact space is a compact space \cite{Eng}.
Then $B$ is a compact space. Let $\{U_i \vert i = 1, \ldots, r\}$
be  a finite open cover of $B$. According to the definition of
induced topology on a subspace, for any $i$ let $U_i = B \cap
U_i^{\prime}$, where $U_i^{\prime}$ is an open set of $X$. Then $A
\cup \{U_i^{\prime} \vert i = 1, \ldots, r\}$ is a open cover of
$X$.

Let $\{V_j \vert j \in J\}$ be a cover of $A$. As $A$ is an open
set then by Lemma 6 $V_j$ is an open set of $X$ for any $j \in J$.
Hence $\{V_j \vert j\in J\} \cup \{U_i^{\prime} \vert i = 1,
\ldots, r\}$ is an open cover of $X$. According to the definition
of compact space, from open cover $\{V_j \vert j\in J\} \cup
\{U_i^{\prime} \vert i = 1, \ldots, r\}$ of $X$ we pull out a
finite open cover $\{V_{j_i} \vert j_i = 1, \ldots, s\} \cup
\{U_i^{\prime} \vert i = 1, \ldots, r\}$. As $A \cap B =
\emptyset$, then $\{V_{j_i} \vert j_i = 1, \ldots, s\}$ is a
finite open cover for $A$. Consequently, $A$ is a compact space.

Now we prove that the space $A$ is regular. Indeed, every
rectifiable space $X$ is regular \cite{Gul} (see, also,
\cite[Theorem 7]{Mal2}. Then, obviously, $X$ is a Hausdorff space.
Every compact (in other terminology  bicompact) Hausdorff space is
normal. This completes the proof of Proposition 15.\vspace*{0.1cm}

From Proposition 15 it follows.\vspace*{0.1cm}

\textbf{Corollary 26.} \textit{Let a subspace $A$ of a compact
space $X$ admit a structure of topological right loop. Then the
set $A$ is compact, normal, regular, Hausdorff and open in $X$.}

\section{Algebraical analysis of some papers}

I have acquired minimal knowledge of topology through my
scientific activity. However I have decided to publish this paper
with the purpose to correct the mistakes found out in papers \cite
{AC1}, \cite {ChobCalm}, published in the journal:
\textsl{Buletinul Academiei de \c{S}tiin\c{t}e a Republicii
Moldova, Matematica.} In paper \cite {AC1} are available reference
on papers \cite{Chob1}, \cite{Chob2}: \cite{AC1} $\Rightarrow$
\cite{Chob1}, \cite{Chob2}. Further, \cite{Chob1}, \cite{Chob2}
$\Rightarrow$ \cite{Chob3}. Initially I planned to analyze papers
\cite{AC1}, \cite{Chob1}, \cite{Chob2} , \cite{Chob3} (the
analysis is presented below). But my plan changed cardinally after
I saw the materials of the \textsl{$20^{th}$ conference of applied
and industrial mathematics, dedicated to Academician Mitrofan M.
Ciobanu, 2012:} www.romai.ro/conferinte.
But my opinion was affected the most by papers \cite{AC1},
\cite{Chob1}, \cite{Chob2} , \cite{Chob3} and other works, for
which references are available. As a result I saw an awful
reality, a catastrophe. It appeared that the works of Academician
Ciobanu M. M. and his followers, specifically monographs
\cite{Dum}, \cite{Chob12}, \cite{Chob6}, \cite{Calm1},
\cite{ChobCalm1}, \cite{Chir1} and dissertation \cite{Calm2},
\cite{Ipate}, \cite{Ciob3}, \cite{Pavel}, \cite{Chir2}, contain
senseless gross errors. These works are anti-scientific, contain a
lot of deceptive statements and lies. Obviously the authors are
unaware of the elementary concepts of algebras and topology.

To stay within the limits of this paper we will analyze only several works,
paying special attention to the algebraic nature of these works. The necessary
theoretical basis is presented in Sections 1 -- 5.

\centerline{\textsf{Continuous signature. Quasivarities. Free
topological algebraic systems}}

We will start with paper  \cite{Chob3}. Ideologically it is connected to papers
\cite{Mal3}, and also \cite{Mal1}. It refers many times to questions related to the
results from\cite{Mal1}, \cite{Mal3}, though no such referrals can be found in \cite{Chob3}.
The fundamental paper \cite {Mal3} contains the basis  of topological algebraic systems,
in particular, of
free topological algebraic systems of given variety (primitive
class) of   such systems. In \cite{Mal1}, \cite{Mal3} any
algebraic system $G$ is defined only with operations of form
$(g_1, g_2, \ldots, g_{n}) \rightarrow g_{n+1}$, $n = 0, 1,
\ldots$, where $g_1, \ldots, g_{n+1} \in G$, i.e. $(n,
1)$-operations, where $n$ is a finite integer.  \cite{Chob3}
considers $(n, m)$-operations with $n < \alpha$, $m < \beta$,
where $\alpha, \beta$ are  any fixed ordinal numbers. However
difference between \cite{Chob3} and \cite{Mal1}, \cite{Mal3} is
really catastrophic. Paper \cite{Chob3} could be compared with a world
of  absurd statements, both in terms of definitions, competence, and logical
and mathematical judgement.  Excuse my lack of modesty, but I think, that
a mentally healthy person cannot even imagine such things, even less publish them.
i

We bring some notions and results from paper \cite[$\S$1]{Chob3}.
\textmd{Fix a non-empty set $G$ and an ordinal numbers $n$, $m$
$q$. An unequivocal  mapping $\omega(n,m): G^n \rightarrow G^m$ is
called \textit{$(n,m)$-operation} on $G$. We put $\Omega^{m}_{n} =
\{\omega(n,m)\}$, $\Omega = \cup \{\Omega^m_n\}$, $P = \cup P_q$,
where $P_q = \{p \vert p \quad \text{is a $q$-ary predicate on}
  G\}$, and let $\alpha_S, \beta_S, \gamma_S$ will be such a
least ordinal numbers that $n < \alpha_S$, $m < \beta_S$, $q <
\gamma_S$.  \textit{The object $<G, \Omega, P>$ is called
algebraic system  of signature $S = \Omega \cup P$} (Definition
1).}

The Definition 1 of algebraic systems has drawbacks, which as will
further result in roughest mistakes. The author does not know that
for rare cases in algebra only  finite-place relations, functions,
words, terms, formulas, operations, predicates are considered
(see, for example, \cite[pag. 42, 138, 141, 167]{Mal1}). For this
look at the comment after  Proposition 1 of present paper.
 Definition 1 for ordinal numbers $n, m, q, \alpha_S, \beta_S,
\gamma_S$  does not impose  even restrictions pag. 108,
 third line above): they can be any ordinal
numbers. It is not clear for what reasons in the last $\S$1 algebraic systems
are investigated, defined by $(1, m)$-operations,
where $m$ is any ordinal number. Let's specify some drawbacks (items
1a), 1b)) of Definition 1.

Remind, that according to \cite{Mal2}, \cite{Mal3} an algebraic
system $\mathcal{A} = <A, \Omega_F, \break \Omega_P>$ is called
topological algebraic system if the basic set $A$ is a topological
space $(A, \mathcal{T})$ and the basic operations of $\Omega_F$
are continuous. From definition of topological algebraic system
follows that the study of  topological algebraic system
$\mathcal{A}$ consists of separate study of the algebraic system
$\mathcal{A}$, or   separate study of  topological space $(A,
\mathcal{T})$, or separate study of  interrelation of algebraic
system $\mathcal{A}$ and  topological space $(A, \mathcal{T})$.

1a). Let $\mathcal{A} = <A, \Omega_F, \Omega_P>$ be an  algebraic
system and let $\omega$ be a basic $(n, 1)$-operation with $n \geq
\infty$. Then by  Proposition 1 of present paper the expression
"subsystem $<H>$ of system $\mathcal{A}$ algebraically generated
by set $\emptyset \neq H \subseteq A$" is not correct, is false.
It is possible only to assert that the subsystem $<H>$ is
generated by set $H$. In general, according to the definition,
$<H>$ is the intersection of all subsystems of $\mathcal{A}$ which
contain $H$ and for $<H>$ it is impossible to apply the
Proposition 1 of present paper. Therefore, if for basic
$(n,m)$-operations of algebraic system $\mathcal{A}$  not to put
any restrictions on ordinal numbers $n, m$ (can be $n \geq \infty$
or $m > 1$), then the statement "the elements in $<H>$ are
expressed through elements in $H$ by means of basic and derived
operations" is without sense.

1b) Let $\mathcal{A} = <A, \Omega_F, \Omega_P>$ be a topological
algebraic system and let $\omega \in \Omega_n^m$. From the definition
of topological algebraic system it follows that a topology $\mathcal{T}$  is
defined on set $A$ and the basic $(n,m)$-operation
$\omega$ is continuous. Let $m > 1$. From last condition it follows
that the given topology $\mathcal{T}$ of $A$ can be a topology on
degree $A^m$. It takes place when $\mathcal{T}$  is  a Tychonoff
topology. Hence, the condition $m > 1$ for basic operations
attracts restriction for given topology  $\mathcal{T}$.

We pass to $\S$2. Let $\mathcal{A} = <A, \Omega_F, \Omega_P>$ be
an algebraic system of signature $S = \Omega_F  \cup \Omega_P$.
Let's present the excerpts.  \textmd{Let some topologies $\tau(n,m)$
be given on set $\Omega^m_n$. In such a case we say that
signature $\Omega$ is \textit{continuous}. If the topologies on
$\Omega^m_n$ are discrete, then  $\Omega$ is called a
\textit{discrete signature.}  Define the mappings $\psi_{(G,n,m)}:
\Omega^m_n \times G^n \rightarrow G^m$, where
$\psi_{(G,n,m)}(\omega, g) = \omega(g)$ for all $\omega \in
\Omega$, $g \in A^n$, and $\varphi_{(G,q)}: P_q \times G^q
\rightarrow L$, where $L$ is a complete latices,
$\varphi_{(G,q)}(p, a) = p(a)$ for all $p \in P_q$, $a \in G^q$.
An algebraical system $G$ together with the topology, given on its
$\mathcal{T}$ and given topologies $\tau(n,m)$ on $\Omega^m_n$ is
called a  \textit{topological system of continuous signature}
$\Omega$, if all mappings $\psi_{(G,n,m)}$, $\varphi_{(G,q)}$ are
continuous (Definition 3.)}

Let's show that  the notion of continuous signature from
Definition 3 is fictitious. It does not have any value and cannot
bring any value added for the research of topological algebraic systems.
As a confirmation hereof, we present some gaps from from the definition of
 continuous signature.

1c). A necessary condition for the topology $\tau(n,m)$ to be
non-discreetly is the infinity of set $\Omega^m_n$. The sum
$\Omega = \cup \{\Omega^m_n\}$ in Definition 3 is discrete (for an
obvious confirmation to this see \cite{Chob2}). Therefore, there
is no connection between topologies $\tau(n,m)$ by the
definition of discrete sum of topological spaces.

1d). Without loss of generality we shall consider that $P =
\emptyset$.  We denote by $f$ the mapping $\psi_{(G,n,m)}:
\Omega^m_n \times G^n \rightarrow G^m$ and assume that $f$ is
continuous. We consider $\Omega$ as a set of symbols of basic
operations.  By Remark 1 of present paper for a fixed operation
$\omega \in \Omega_n^m$ the mapping $\mu: \Omega_n^m \rightarrow
G^m$: $\omega \rightarrow h = \omega(g)$ for a fixed $g \in G^n$
is continuous. Similarly, for a fixed $\omega \in \Omega_n^m$ the
mapping $\nu: g \rightarrow h = \omega(g)$, $g \in G^n, h \in G^m$
is continuous. It means that the topology $\tau(n,m)$ on
$\Omega_n^m$ does not bear any  topological and algebraic
information of topological system  $<G, \Omega>$, i.e. does not
change the topological and algebraical structure of $<G, \Omega>$,
does not change the topological and algebraical structure of
objects of any class $\mathcal{K}$ of topological algebraic
systems of signature $\Omega$. Consequently, the topology
$\mathcal{T}$ on $G$ does not depend on topologies $\tau(n,m)$ on
$\Omega_n^m$. Finally, if the mappings $\mu, \nu$ are continuous
then by  Remark 1 of present paper the mapping $f$ should not be
necessarily continuous. Hence the definition of topology
$\tau(n,m)$ on set $\Omega_n^m$ does not depend  on the  topology
of algebraic system $<G, \Omega>$ and on the requirement for the
continuity of mapping $\psi_{(G,n,m)}: \Omega^m_n \times G^n
\rightarrow G^m$,

1e). Consequently, without prejudice to research of classes of
topological algebraic systems $\mathcal{K}$, also separate systems
of $\mathcal{K}$ expressed in the form $''$\ldots system of
continuous signature . . .$''$ the words $''$of continuous
signature$''$ should be omitted. In such  a case almost all
results of paper (for example, the sections 4, 7, 8 and others)
become trivial or without sense.  The reasoning for introducing
notion of "continuous signature" is not clear. But if necessary,
it is possible to define the topology $\tau(n,m)$ considering
$\Omega_n^m$ as a set of symbols of basic operations.

$1e_1)$. Further, we quote by \cite[pag. 110, 111]{Chob3}.
\textmd{Fix a class of topological systems of continuous
signature.} \textmd{Let's allocate the following properties, which
can have the class $\mathcal{K}$.}

\textmd{Condition 1ě. Closed with respect to subsystems.}

\textmd{Condition 2ě. Closed with respect to Tychonoff products.}

\textmd{Condition 3ě. All objects from $\mathcal{K}$ satisfy a
given  set of  topological, algebraical properties $Q$ (the set
$Q$ can be $\{\emptyset\}$ or the separation axioms $T_0, T_1,
T_2, \break T_3$, or the requirement of a regularity, or the
requirement of a completely regularity and others).}

\textmd{Condition 4ě. If $(G, \tau) \in \mathcal{K}$, $\tau_d$ is
the discrete topology and  $(G, \tau_d)$ is a topological system
of continuous signature, then $(G, \tau_d) \in \mathcal{K}$.}

\textmd{Condition 5ě. Closed with respect to continuous
homomorphic images with property $Q$.}

\textmd{Condition 6ě. If $(G, \tau_1) \in \frak K$ and with
respect to topology $\tau$ the pair $(G, \tau)$ is a topological
system with properties $Q$ of continuous signature $S$, then $(G,
\tau) \in \frak K$.}

\textmd{Definition 4. A class $\mathcal{K}$ is called:}

\textmd{1. Topological $Q$-quasivariety, if  the conditions  1m -
4m are satisfied.}

\textmd{2. Topological $Q$-prequasivariety, if  the conditions 1m
- 5m are satisfied.}

\textmd{3. Topological $Q$-variety, if  the conditions  1m - 6m
are satisfied.}

\textmd{4. Topological complete $Q$-quasivariety, if  the
conditions 1m - 4m and 6m are satisfied.}

\textmd{Definition 5. A class $\mathcal{K}$ of algebraical systems
of signature $S = \Omega \cup P$ is called quasivariety, if}

\textmd{1. The class $\mathcal{K}$ is closed with respect to
subsystems;}

\textmd{2. The class $\mathcal{K}$ is closed with respect to
Cartesian product.}

The notions defined in the following excerpt  are  a basic object of
research not only in paper \cite{Chob3}, but also in many
other papers, therefore we analyze this excerpt in more details.
 But at first remind that according
to \cite{Mal2}, \cite{Mal3} an algebraic system for which the
basic set of  elements is a topological space (one and only
one space) and the basic operations are continuous is called
\textit{topological algebraic system}. The totality of topological
algebraic systems (defined above by one and only by one topological
space) which make a variety (respect. quasivariety) in algebraic
sense is called \textit{a variety  of topological algebraic
system} (respect. \textit{a quasivariety  of topological algebraic
system}.

1f). Let's show that the defined notions are not correct, without
sense, a total chaos, and do not correspond to the classical
 notions. We give reason for it.

$1f_1)$. The conditions 1m, 2m coincide with conditions $K_1)$,
$K_2)$ (Section 2) by  \cite[Proposition 2.3.1]{Eng}. In condition
4m the relation $(G, \tau_d) \in \frak K$ is false. As
 marked in the comment below   conditions $K_1)$, $K_2)$
(Section 2) the discrete topology does not satisfy  conditions
$K_1)$, $K_2)$. It generates topological spaces which are
regular and completely discontinuous in the sense that for any
finite set of points $x_1, \ldots, x_n $ of such space $X$
there exists a splitting $X$ on $n$ in pairs not crossing closed sets,
containing  one of the specified points. These topological
spaces satisfy  conditions $K_1)$, $K_2)$.

$1f_2)$. Even not considering item 1d) the topologies $\tau_1$ and
$\tau$ from the conditions are "almost" independent by Remark 1of
present paper. Hence the topology $\tau$ from relation $(G, \tau)
\in \frak K$ is "almost" arbitrary. To what extent such topologies
$\tau$ are available and how are are connected with conditions 1m,
2m?

1g). By  item 1f) the notions from Definition 4 are not correct,
false. Let's show that a similar situation takes place and the
definition of quasivariety (Definition 5).

$1g_1)$. In the literature, \cite{Con}, \cite{Mal1} and others
(see Section 2) the notion of quasivariety (or quasiprimitive
class) is defined  by means  of quasiidentities. An quasiidentity
is a certain set
 of quasiatomic formulas, of
$\forall$-formulas. Any quasiatomic formula is a finite totality
of free variables of basic operations of algebraic system. Any
basic operation of algebraic system is a $(n, 1)$-operation, where
$n$ is a finite integer. Definition 4 examines $(n,
m)$-operations with arbitrary numbers $n, m$, that is a roughest
mistake.

$1g_2)$. The definition of quasivariety  from Proposition 5 is
incorrect. Besides conditions 1g, 2g it is necessary to add the
conditions $Q_1)$, $Q_4)$ of Proposition 3 from Section 2. The
condition $Q_4$): the class  $\frak{K}$ contains an unitary system
is a characteristic property for quasivarieties by Corollary 8.
The condition $Q_1$): the class  $\frak{K}$  is ultraclosed is
also necessary in the definition of quasivariety. There are
classes $\frak{K}$ with conditions $Q_2)$, $Q_3)$, $Q_4)$ which
are not quasivarieties (\cite[pag. 272]{Mal1}).

$1g_3)$. Further, the notion of Cartesian product does not characterize
fully the quasivariety (\cite[pag. 193]{Mal1}. The notion of filter
product corrects this drawback. For example, by
\cite[Corollary V.11.3]{Mal1} a class $\frak K$ of algebraic
systems is a quasivariety if and only if $\frak K$: (1) is closed
with respect to filter products; (ii) is closed with respect to
captures of subsystems; (iii) contains an unitary system.

1h). The section 3 $''$Free algebraical systems$''$ of \cite[pag.
112 - 118]{Chob3} begins with the following definitions, which are
similar to the Definitions in \cite{Mal3}.

\textmd{Definition 7 (resp. Definition 8). We fix  a non-trivial
class $\mathcal{K}$ of topological systems of continuous signature
$\Omega$ and non-empty set (respect. space) $X$. The pair
$(F^a_{\mathcal{K}}(X), i_X)$ (resp. $(F_{\mathcal{K}}(X),
\delta_X)$), when $F^a_{\mathcal{K}}(X)$, $F_{\mathcal{K}}(X) \in
\mathcal{K}$ and $i_X: X \rightarrow F^a_{\mathcal{K}}(X)$ (resp.
$\delta_X: X \rightarrow F_{\mathcal{K}}(X)$) is unequivocal
mapping (resp. continuous mapping), is called algebraically free
system of set (respect. space) $X$, if it satisfies the
conditions:}

\textmd{1a (resp. 1t). The set $i_X(X)$ (resp. $\delta_X(X)$)
algebraically generate the system $F^a_{\mathcal{K}}(X)$ (resp.
$F_{\mathcal{K}}(X)$);}

\textmd{2a (resp. 2t). For every unequivocal mapping  $f: X
\rightarrow G$, when $G \in \mathcal{K}$, it is exists a
homomorphism $\hat{f}: F^a_{\mathcal{K}}(X) \rightarrow G$
(respect. continuous homomorphism  $\hat{f}: F_{\mathcal{K}}(X)
\rightarrow G$ ) such that $f(x) = \hat{f}(i_X(x))$ (resp. $f(x) =
\hat{f}(\delta_X(x))$) for all $x \in X$.}

At the first sight the introduced notions (and their properties
proved below) are an essential generalization of classical
notion of topological algebra of a given variety of topological
algebras with given generating  topological space and with given
generating relations, and of  free objects of varieties of
topological algebras (and their properties), which were introduced and
investigated in detail by A. I. Mal'cev in \cite{Mal3}.
Actually the introduced notions and their properties from
\cite{Chob3} are a senseless imitation of \cite{Mal3} (which are not
mentioned in \cite{Chob3}), are a complete absurdity. We will
give arguments to support this. But first we note that a complete
generalization of Mal'cev's notions and results is presented in
Sections 2, 3.

$1h_1)$. The condition "of continuous signature"  should be ignored
by virtue of item 1e).

1i). With such a wording the item 1a (resp. 1t) is incorrect: it
is impossible to use the expression "algebraically generated" by
item 1a) and Proposition 1. It is meaningful if and only if the
basic operations of topological systems in $\mathcal{K}$ are $(n,
1)$-operations with finite integer $n$. With such restrictions on
basic operations the Definition 7 (resp. Definition 8) is a
special case (together with Lemma 2) of Definition 1 (see, also
\cite{Mal3}) provided that   the item 1t is replaced by the
condition: the set $\delta_X$ topologically generate the system
$F_{\mathcal{K}}(X)$.

Below Definition 7 (resp. Definition 8) one may find it follows something
improbable, a reasoning without sense, which does not follow any
mathematical logic. To present Properties 1 - 15, their corollaries,
and definitions that describe the systems from Definitions 7, 8 in a case
when the class of topological spaces $\mathcal{K}$ satisfies only the
conditions 1m, 2m, let us present some of them. Let us have a non-trivial class
$\mathcal{K}$ of topological systems of continuous signature $S =
\Omega \cup P$ which satisfies the conditions 1m, 2m. Then the following hold.

\textmd{Property 1. For each non-empty space $X$ exists an unique
algebraically free topological system  $(F_{\mathcal{K}}^a(X),
i_X)$  of space $X$ in the class  $\mathcal{K}$ and the mapping
$i_X$ which to apply one-to-one  $X$ in $(F_{\mathcal{K}}^a(X)$.}

\textmd{Property 2. For each non-empty space $X$ exists an unique
topologically free topological system $(F_{\mathcal{K}}(X),
\delta_X)$ of space $X$ in class  $\mathcal{K}$.}

\textmd{Property 3. For each non-empty space $X$ exists an unique
continuous homomorphism $g_X: F^a_{\mathcal{K}}(X) \rightarrow
F_{\mathcal{K}}(X)$ such that $i_X \circ g_X = \delta_X$.}

\textmd{Property 4. If the space $X$ is discrete then the
homomorphism $g_X$ is a topological isomorphism.}

\textmd{The following fact shows, that the conditions 1m, 2m are
not only sufficient for existence of free systems, but also almost
necessary.}

\textmd{Property 13. Let $\mathcal{K}$ be a quasivariety of
algebraic systems of signature $S$. Then for  each non-empty space
$X$ exists an unique algebraically free topological system
$(F_{\mathcal{K}}^a(X)^a, i_X)$. If $\mathcal{K}$ contains a not
one-element system, then  $i_X$ one-to-one apply  $X$ in
$(F_{\mathcal{K}}^a(X)$.}

\textmd{Proof. Similarly of proof of Property 1.}

To prove Properties 1 - 4, 13 the following Proposition 2 is used essentially.

\textmd{Proposition 2. Let $G$ be an algebraic system of signature
$\Omega = \cup \Omega_n^M$ with $n < \alpha_S$, $m < \beta_S$.
Let, further, $A(Y)$ be a subsystem of system $G$ generated by
non-empty set $Y \subseteq G$. Then $|A(Y)| \leq |\Omega \cup
Y|^{\tau} + 2^{\tau}$, where $\tau = \max\{|\alpha_S|, |\beta_S|,
\aleph_0\}$.}

\textmd{Proof. Denote by $\tau^{+}$ the first ordinal number more
than $\tau$. If $\omega \in \Omega_n^m$ and $g \in G^n$ then
$\omega(g) = \{g_\xi \vert \xi < m\} \subseteq G^m$ and put
$B(\omega(g)) = \}g_\xi \vert \xi < m\} \subseteq G$. For any
non-empty set $Z \subseteq G$ denote $\Omega_n^m(Z) =
\cup\{B(\omega(g)) \vert g \in Z^n \subseteq G^n$ and $\omega in
\Omega_n^m \}$ and $\Omega(Z) = \cup\{\Omega_n^m(Z) \vert n <
\alpha_S, m < \beta_S\}$. Let's observe, that $|\Omega_n^m(Z)|
\leq |\Omega_n^m \times (Z^n)| \cdot |m| = |\Omega_n^m | \cdot
|Z^n| \cdot |m|.$ Hence $|\Omega(Z)| \leq |\Omega \cdot |Z^{\tau}|
\cdot \tau$. Put $A_0(Y) = Y$ and $A_{\alpha}(Y) =
\Omega(\cup\{A_{\xi}(Y) \vert \xi < \alpha\})$ for all $\alpha <
\tau^{+}$. Then for every $\alpha < \tau^{+}$ we have
$|A_{\alpha}(Y)| \leq |\Omega \cup Z|^{\tau} + 2^{\tau}$. As $A(Y)
= \cup \{A_{\alpha}(Y) \vert \alpha < \tau^{+}\}$, that proof is
completed.}

1j). With such wordings the Proposition 2 and its proof are a
roughest absurdity, which cannot be explained by any means. The
nature of $ \aleph_0 $ is not clear. It is known that
$2^{\aleph_0}$ is the continuum. Then what kind of restrictions
$\leq$ are mentioned in Proposition 2? The role of number $\tau^+$
is also not clear. Statement: "As $A(Y) = \cup \{A_{\alpha}(Y)
\vert \alpha < \tau^{+}\}$, that proof is completed." contains
severe errors. Remind that the subsystem $A(Y)$ is the
intersection of all subsystem of $G$ which contain the set $Y$. If
$m > 1$, for example $m = 5$, for any $n,m$-operation $\omega \in
\Omega$, then it is not clear how to express every element of
$A(Y)$ through elements of $Y $ with the help of operations
$\omega$, i.e. that equality $A(Y) = \cup \{A_{\alpha}(Y) \vert
\alpha < \tau^{+}\}$ is correct. According to Proposition 1 of
present paper this equality holds if and only if all operations in
$\Omega$ have a form $(n, 1)$-operation with finite integer $n$.
In such case the statement of Proposition 2 is known: $|A(L)| \leq
\max \{|\Omega|, |Y|, \aleph_0\}$ \cite[Example 2, pag.
163]{Mal1}.

On proof of Properties 1 - 4, 13. Schematically the proof is
without sense, imitating the proof of Theorem 2 from \cite{Mal3}(see
also the Theorem 6), though no reference is made in the paper about it.
We present and will analyze  the  proof of Properties 1 - 4. First,
let's note that the proof is erroneous as it used the
incorrect Proposition 2.

Proof of Properties 1 - 4. \textmd{Uniqueness to within natural
isomorphism  easily follows from Definition 7. Let's prove
existence of algebraically free system. Put $\tau = \max\{|X|,
|\Omega|, |\alpha_S|, |\beta_S|, \aleph_0\}$. By Proposition 2  we
can consider( only a systems of cardinality $\leq 2^{\tau}$. The
totality  $\mathcal{K}_{\tau} = \{G \in \mathcal{K} \vert |G| \leq
2^{\tau}\}$ is a set. Then is set and the totality $\{f_{\xi}: X
\rightarrow G_{\xi} \vert \xi \in A\}$ all mappings of set  $X$ in
systems  $G_{\xi}$ from $\mathcal{K}_{\tau}$. We consider the
diagonal product $i_X \rightarrow G = \prod G_{\xi} \vert \xi \in
A\}$. Let $F^a_{\mathcal{K}}(X)$ is the subsystem of system $G$
generated by set $i_X(X)$. Using projections $\hat{f}_{\xi} =
\pi_{\xi}: F^a_{\mathcal{K}}(X) \rightarrow G_{\xi}$ and
Proposition 2 easily to show that $(F^a_{\mathcal{K}}(X), i_X)$ is
an algebraically free system of set $X$ in class $\mathcal{K}$.
One-to-one unambiguity of mapping  $i_X$ follows from a condition
of  not trivialities of class $\mathcal{K}$. To it the property 1
is proved. Let now  $X$ is a  topological space. Put $B = \{\alpha
\in A \vert f_{\alpha}: X \rightarrow G_{\alpha} \quad \text{is
 continuous}\}$.  Considering diagonal product $\delta_X =
\prod \{f_{\alpha} \vert \alpha \in B\}: X \rightarrow H = \prod
\{G_{\alpha} \vert \alpha \in B\}$, we similarly let's construct a
topologically  free system $F_{\mathcal{K}}(X), \delta_X)$ of
space $X$ in class $\mathcal{K}$. If the  space $X$ is discrete,
then $A = B$, hence $F^a_{\mathcal{K}}(X) = F_{\mathcal{K}}(X)$.
It proves the following three properties 2 - 4.}

By item 1j) the proof of Properties 1 - 4 is erroneous as it used
the Proposition 2. Assume that the Proposition 2 holds. It is
known that $2^{\aleph_0}$ is continuum. We can assume that
$\tau > \aleph_0$. It is not clear why totality
$\mathcal{K}_{\tau}$ is a set. Let's assume that this is correct. We
consider the Cartesian product $G = \prod G_{\xi} \vert \xi \in
A\}$, where $A = \{\xi \vert G_{\xi} \in \mathcal{K}_{\tau} \}$.
Assume that all systems $G_{\xi}$ are non-empty. There is a
question: will the Cartesian product be non-empty? If the
totality of indexes is finite, then the answer to this question is
positive. Non-emptiness of set $G_{\xi}$ means that in each set
$G_{\xi}$ it is possible "to choose" at least one element. Having made
such a choice for every $\xi$ we will receive an element of
Cartesian product. If the totality $A$ is infinite, then it is necessary
to make an infinite number of choices. But it is an axiom of choice. Hence,
for $\prod G_{\xi} \vert \xi \in A\} \neq
\emptyset$ it is necessary to require that it should be possible
to index  the totality $A$. Only in such a case it is possible to
speak about existence of system $F^a_{\mathcal{K}}(X)$. To have a free
system $F^a_{\mathcal{K}}(X)$ it is necessary to
require that the totality of all projections $\\pi_{\xi}:
F^a_{\mathcal{K}}(X) \rightarrow G_{\xi}$ is complete
\cite[Theorem I.2.2]{Mal1}.

According to the above stated we conclude.

1k). The entire material from Section 3 of \cite{Chob3}, where
the notions of algebraically (or topologically) free
topological systems of various  classes of topological systems are
defined and the various properties of such systems  are investigated,
particularly their existence, are erroneous, incorrect, and without
any sense. Similar questions are investigated in detail in
Section 3 and Section 2 of this paper.

1l). Properties 1 - 4 prove the existence and uniqueness of the
free systems from Definitions 7, 8. These properties are the basic tool
for all further researches of this paper. By item 1k), and 1e), 1f), 1g), 1h), 1i)
we conclude that the entire content of paper \cite{Chob3}, including numerous
generalizations of known results, is erroneous and senseless, consisting of
statements assembled mechanically as proofs. The given paper \cite{Chob3} is
not mathematical. I even do not know how to name it. To confirm these statements
we will analyze some more excerpts from the given paper.

In $''$Introduction$''$ it is mentioned that the main purpose of paper \cite{Chob3}
is \textmd{the decision
of problem on representation of  topological systems in form of
quotient homomorphic images of zero dimensional systems from same
class. This problem is completely solved for discrete signature.
(As corollary it follows that  amplifies the result from
\cite{Arh4} and given problem is positively solved for topological
universal algebras of discrete signature, for topological
semigroups, for topological quasigroups and loops, for topological
rings, and also for topological groups with multioperators, where
on set of multioperators to consider the discrete topology. It is
possible to consider $\Phi$-algebras also with set of
differentiated $\Delta$, where $\Phi$, $\Delta$ are discrete
(see, pag. 132). If the signature is not discrete, that follows
from the Theorem 5, then the given problem has a negative
decision even in rather well known varieties (see, item 1q). For
topological groups this problem is solved completely in
\cite{Arh4}.} The proof scheme from \cite{Chob3} is similar to the
proof scheme from \cite{Arh4}. Except for the concepts and some results
taken from \cite {Arh4}, it is absolutely impossible to regard as authentic
the solution of the problem stated in \cite {Chob3}. Let's specify only a
few rough mistakes in comparison with \cite{Arh4}.

1m). The expressions "continuous signature", "discrete signature"
should be ignored by item 1e).

1n). The problem is solved for classes of topological systems for
which the conditions 1m, 2m are correct. Thus these are supposed
mistakes of items 1f), 1g).

1o). The following concept is used to solve the problem in \cite{Arh4}.
It is known (see Corollary 2 of present paper),that any  continuous
homomorphism of topological group  on
topological group  is open. For factor-system of topological
system the property  of open homomorphism is not always correct as
the conditions of a Consequence 2 are not always fair.

1p). In \cite{Arh4} in quality of group $G$ from item 1f is
Considered (the free topological group described below Definition
1 of present paper. For classes of topological systems in
\cite{Chob3} is considered the topologically free topological
system from Definitions 7, 8. The admitted mistakes are specified in
the analysis of validity   of Properties 1 - 4.

1q). Really, \textmd{the class  $\mathcal{K}_l$ of all  linear
topological spaces form a topological variety of topological
algebras with continuous signature. But in $\mathcal{K}_l$ only
one-pointed spaces are null-dimensional.} This is incorrect. The
author confuses elementary notions: the dimension of linear space
and the dimension of topological space. Example: the subset of
$n$-dimensional Euclidean  space, consisting from all points with
real coordinates, is a null-dimensional topological space as this
subset does not contain intervals.

Now we pass to analysis of $\S$4 and let's show complete
incompetence of the author with respect to theory of identities
and quasiidentities of algebraic systems. Like other sections of
analyzed paper, the $\S$4 consists only of senseless notions
and  statements. To back this we will use the following excerpt.
Fix a continuous signature and an algebraic system $G$ of signature $S$.
Let $\frak K(S, i)$ denote the totality of all algebraic systems of signature
$S$ which satisfies the separated axiom  $T_i$. Fix $i \in \{-1, 0,1,
2, 3, 3\frac{1}{2}, r, \rho\}$. The class $\frak{K}(\emptyset, i)$
consists of all non-empty $T_i$-spaces. For any non-empty space
$X$ we denote by $(X/i, g_{iX})$ such a $T_i$-spaces $X/i$ and
continuous mapping $g_{iX}: X \rightarrow X/i$ that for every
continuous mapping $f: X \rightarrow Y$ in $T_i$-space $Y$ a continuous
mapping $f_i: X/i \rightarrow Y$ exists that $f = f_i \circ
g_{iX}$. Actually,  the pair $(X/i, g_{iX})$ is a free topological
system of space $X$ in class $\frak{K}(\emptyset, i)$. For class
$\frak{K}(\emptyset, i)$ the signature is empty. We consider the
signature $S = \Omega \cup P$. If $\omega, \mu \in
\Omega_{\alpha}^{\beta}$, then the identity of form $\omega = \mu$
is called a simple identity. (Absolutely unclear: the
identity from what? Actually, $\omega, \mu \in
\Omega_{\alpha}^{\beta}$ are fixed elements of the set of symbols of
basic operations of fixed signature $S = \Omega \cup P$. The
equality $\omega = \mu$ means that the basic operations $\omega,
\mu$ coincide). Actually, the identities can have a very difficult
form. Really, let $\omega \in \Omega_{\alpha}^{\beta}, \mu \in
\Omega_{\lambda}^{\beta}$ and let the set $A$ and the
mappings $f: A \rightarrow \prod(\alpha), g: A \rightarrow
\prod(\lambda)$ be given, then for any $x = \{x_i\} \in G^{\alpha}$ and $y
= \{y_{\eta}\} \in G^{\lambda}$ we have $\omega(x) = \mu(y)$ as
only $x_{f(a)} = y_{f(a)}$ for all $a \in A$. Instead of basic
operations it is possible to take a derived operations, then under
$\beta_S \leq 2$ and $P = \emptyset$ thus all identities can be
received. Let a quasivariety $\frak K$ of algebraic
systems of signature $S$ be given. Through $J(\frak K)$ and $qJ(\frak K)$
denote the totality of all identities and quasiidentities which
define respectively the class $\frak K$ (uttermost absurdity:
according to definition of quasivarieties they are not defined
with the help of identities). Through  $Q\frak K$ denote the
totality of all systems from $\frak K$, for which topologies with
properties $Q$ can be defined. Let's say that $Q\frak K$ is topologically
defined  by identities $J(\frak K)$ and
quasiidentities $qJ(\frak K)$. The identities $J(Q\frak K)
\backslash J(\frak K)$ and the quasiidentities $qJ(Q\frak K)
\backslash qJ(\frak K)$ are called relatives identities and
quasiidentities of topological type. The identities $J(Q\frak K)
\backslash J(\frak K)$ and the quasiidentities $qJ(Q\frak K(S,i))
\backslash qJ(\frak K(S,i))$ are called absolute identities and
quasiidentities of topological type. Let $J$ denote the set of all
simple identities. Put $J_{\alpha}^{\beta} = \{(\omega = \mu) \in
J \vert \omega, \mu \in \Omega_{\alpha}^{\beta}\}$. On space
$\Omega_{\alpha}^{\beta}$ the identities $J_{\alpha}^{\beta}$
generate a partition of classes of equivalence. The space of
partitions in factor topology we denote by
$\Omega_{\alpha}^{\beta}/J$ and by $\pi_{\alpha\beta}:
\Omega_{\alpha}^{\beta} \rightarrow \Omega_{\alpha}^{\beta}/J$ we
denote the projection. By $\frak K(S, J,  i)$ denote all
topological systems from $\frak K(S,i)$ on which the identities
$J$ are true.

Further these senseless notions are used to describe the free
topological system $F_{\frak K(S, J, i)}$ (statements 1 - 5).
Moreover, in Example 7 even is given
constructive construction such free systems and a maximal (in some
sense) topologies on such free systems. We cannot even speak about
its validity.  According to item 1l), such free systems do not exist.

To raise the paper‘s prestige the author uses the following methods
(at least in the author’s opinion).

1r). In the beginning of paper there is a statement: I use this
opportunity to express my gratitude to professor A. V. Arhangel'skii,
the conversations with whom have inspired the author to write the paper.

1s). The author claims that the $X_{\lambda}$-topologies construction
procedure, offered by A. I. Mal'cev \cite[pag. 448]{Mal3} served as a
basis for $ \S$4.  This is far from reality.

1t). In some of its papers \cite{Chob5}, \cite{Chob51},
\cite{Chob6}, \cite{Chob8}, \cite{ChobChob}, \cite{ChobChir44} the
author changes the name of the journal, where the given paper
\cite{Chob3} is published. Specifically: M. M. Cho\-ban, On the
theory of topological algebraic systems, Trans. Amer. Math. Soc.,
48, 1986, 115 - 159. But this is not true. Paper \cite{Chob3} was
published in \textsl{$''$Trans. Moscov. Math. Soc.$''$}

$1t_1)$. Moreover, various monographs\cite{Calm}, \cite{CalChob2},
\cite{Chir1},  dissertations \cite{Ciob3}, \cite{Pavel},
\cite{Chir2}, \cite{Ipate} have references to non-existent works,
which were not published in the indicated prestigious journals
\cite{Chir3}, \cite{ChobCiob2}, \cite{ChobKir4}, \cite{ChobKir5},
\cite{CalmChob} and others.

1u). We analyzed only some excerpts from paper \cite{Chob3}. The situation
of the remaining part (Sections 4 - 13) is much worse: absurd notions and results.
Thus, on the basis of the above-mentioned it is possible to ascertain that paper
\cite{Chob3}, which consists only of unclear proofs and severe errors, does not
represent any scientific interest. Unfortunately, such statements are not enough
to describe the paper of Academician Choban M. M. \cite{Chob3}. The absurd and
senseless notions and
results, as well as absurd and senseless proofs, are still used as basic
objects and tools of research for many papers, monographes and
dissertations.
We listed some of them: \cite{Chob12}, \cite{Chob6},
\cite{Chob7}, \cite{Chob8}, \cite{Chob19}, \cite{Dum},
\cite{Chir1}, \cite{Chir2}, \cite{DC}, \cite{ChobChob},
\cite{Calm1}, \cite{Ipate}. Naturally, they are erroneous. Let us analyze
briefly some of them.

With some "peculiarities" paper \cite{DC} introduced and studied the same
notions and results as paper \cite{Chob3}.

 \textmd{Let $G^A$ denote  the totality of indexed
 sequence $\{x_{\alpha} \vert \alpha \in A\}$ of elements of $G$. We
consider a fixed set of indexes and a non-empty set $G$. Every
mapping $\omega: G^A \rightarrow G$ is called $A$-ary operation on
$G$. Often the set $A$ is identified with set of ordinal numbers
which are less than some given number $\alpha$. In such case
$A$-ary operation is called $\alpha$-ary operation.}

\textmd{Let is given the set $B \supset A$. Then the $A$-operation
$\omega$ it is possible to considered as  $B$-operation, if to
identify $\omega$ with operation $\omega_1: G^B \rightarrow G$
such that for each point $\{x_{\beta} \vert \beta \in B\} \in G^B$
$\omega_1(\{x_{\beta} \vert \beta \in B\}) = \omega(\{x_{\beta}
\vert \beta \in A\})$. It allows to consider any given set of
operations $\Omega = \{\omega_{\gamma} \vert \gamma \in \Gamma\}$
above the same set of indexes $A$. If $\omega_{\gamma}$ is an
$A_{\gamma}$-ary operation, then let's assume $A =
\cup\{A_{\gamma} \vert \gamma \in \Gamma\}$ and  we shall consider
that all operations $\omega_{\gamma}$ are $A$-ary} (pag. 29).

2a). Severe error. By equating the finite arity of operation with
non-finite arity one of the basic properties of operation with
finite arity, specified in Proposition 1 of that paper, is lost.
Moreover, according to Proposition 1 it is necessary to consider
only operations of finite arity.

\textmd{If on $\Omega$ is given a topology, then $\Omega$ is
called continuous system of operations. If is given a some set $G$
together with mapping $\Psi: \Omega \times G^A \rightarrow G$,
that speak, that is given an universal algebra $G$ with system of
$A$-ary operations $\Omega$ or is given an $\Omega$-algebra. If
thus $\Omega$ is a continuous system of operations, $G$ is a
topological space it the mapping $\Psi: \Omega \times G^A
\rightarrow G$ is continuous with respect to topology of Tychonoff
product in $\Omega \times G^A$, then $G$ is called continuous
$\Omega$-algebra.}

2b). According to item 1e) the notion of algebra of continuous
signature is fictitious, it can be ignored.

Besides the concept of free topological algebra with continuous
Signature, paper \cite{DC} investigates the following basic notions.
A non-empty class $\mathcal{K}$ of
$\Omega$-algebras is called prevariety if:

A1) $\mathcal{K}$ is closed with respect to direct product;

A2) $\mathcal{K}$ is closed with respect to taking of
subalgebras.

A non-empty class $\mathcal{K}$ of continuous $\Omega$-algebras is
called $T_i$-prevariety, if $\mathcal{K}$ satisfies the conditions
1m, 2m, 3m, 6m, considered above for paper \cite{Chob3}. Similarly
to Definition 7, it defines the notion of free algebra for
algebraical prevariety. The notion of free topological algebra for
$T_i$-prevariety coincides with Definition 8. As mentioned in
item 1u) the notions of prevariety and it free  algebras, of
$T_i$-prevariety and its free topological algebras are a basic
object of researches of many papers. Therefore, we analyze the
evidence of  existence and uniqueness of such free algebras. We
mention only that these proofs are an unsuccessful imitation of
Mal'cev's proofs \cite[Theorem 1]{Mal3}, \cite[Theorem
V.11.4]{Mal1}, and also of Theorem 6 from that paper.

We present it literally.

\textmd{LEMMA 1. Let $G$ be a  $\Omega$-algebra, $X \subset G$ $X
\neq \emptyset$. Then $X$ generate in $G$ a subalgebra $G_1$, for
 which $|G_1| \leq (|X \cup \Omega| + 2)^{2|\Delta|} +
\aleph_0$.}

\textmd{PROOF. Let us note  previously $\Omega^{\star} = \Omega
\cup \{\omega^{\star}\}$, where $\omega^{\star}$ is the identical
operation. For every subset $M \subset G$ let's assume
$\Omega^{\star} = \omega^{\star} \cup \cup (\omega(M^A) \vert
\omega \in \Omega\}$. Then we have $|\omega^{\star}(M)| \leq
(|\Omega \cup M| + 2)^{|\Delta| + 1}$. Let $\lambda_0$ be the
first ordinal of cardinality $2^{|A|} + \aleph_0$. We put  $X_0 =
X$ and $X_{\lambda} = \Omega^{\star} \times (X_{\alpha} \vert
\alpha < \lambda\})$ for all $\lambda < \lambda_0$. Let $G_1 =
\cup\{X_{\lambda} \vert \lambda < \lambda_0\}$. Clearly that the
cardinality of $G_1$ does not exceed $(|X \cup \Omega| +
2)^{2|\Delta|} + \aleph_0$. It is necessary to show, that $G_1$ is
a subalgebra of algebra $G$. Really, if $\{x_{\alpha} \vert \alpha
\in A\} \in G_1^A$, that will be such an element $\lambda_1 \leq
\lambda_0$ that $\{x_{\alpha} \vert \alpha \in A\} \in
X_{{\lambda}_1}$. Then we have  $\Omega(\{x_{\alpha} \vert \alpha
\in A\}) \subset X_{\lambda_1 + 1} \subset G_1$. The lemma  is
proved.}

2c). Gross mistake by Proposition 1 from that paper: $G_1$ is not
a subalgebra of algebra $G$. It is not clear what $\Delta$ means.
For more details see item 1j).

Lemma 1 is used essentially in the proof of Theorem 1.

\textmd{THEOREM 1. What were not $T_i$-prevariety $\mathcal{K}$
with continuous system of $A$-ary operations $\Omega$ and the
topological space $X$ the free topological algebra
$(F_{\mathcal{K}}(X), \sigma)$ exists and is  unique.}

\textmd{PROOF. The  uniqueness is proved similarly of
\cite[Property 1]{Chob3} to expound above. Let us prove now the
existence of $F_{\mathcal{K}}(X)$. We consider the totality
$\mathcal{B}$ of all $G \in \mathcal{K}$, for which $|G| \leq (|X
\cup \Omega| + 2)^{2|\Delta|} + \aleph_0$. Obviously,
$\mathcal{B}$ forms a not empty set, if to consider identical the
topologically isomorphic algebras from  $\mathcal{K}$. Therefore
the family all continuous mappings $\mathcal{F} = \{g: X
\rightarrow G \vert G \in \mathcal{B}\}$ is also non-empty set. We
consider the diagonal product $\sigma = \Delta \mathcal{F}: X
\rightarrow \Pi \{G \vert G \in \mathcal{B}\}$. The set $\sigma X$
generate in $\Pi \{G \vert G \in \mathcal{B}\}$ a some subalgebra
$F$. Then $(F. \sigma) = (F_{\mathcal{K}}(X), \sigma)$ is the free
topological algebra. The Theorem is proved.}

\textmd{THEOREM 2. Let $\mathcal{K}$ be an algebraic prevariety of
$\Omega$-algebras. Then for any non-empty set $X$ the free algebra
$F_{\mathcal{K}}^a(X)$ it is exists and is  unique( (to within to
isomorphism above $\sigma_a X$.}

\textmd{PROOF. The existence and the uniqueness of algebra
$(F_{\mathcal{K}}^a(X, \sigma_a)$ is established similarly of
existence and of uniqueness of algebra $(F_{\mathcal{K}}(X),
\sigma)$ from  Theorem 1.}

2d). It is not clear how to estimate the authenticity of proofs of
Lemma 1 and Theorems 1, 2. It is possible to use items 1j), 1k),
1l). I think that the reader will agree with me, that the proofs
of Lemma 1 and Theorems 1, 2 are a real non-sense, a delirium.
Hence, the existence and uniqueness of free algebras
$(F_{\mathcal{K}}(X), \sigma)$, $(F_{\mathcal{K}}(X), \sigma)$ are
not proved. It is impossible to prove the existence of such free
algebras operating with the definitions of prevariety and
operations from Proposition 1 and Theorems 6, 7 of that paper.

On page 33 is defined the notion of point finite arity operation.
Actually, it is the notion of operation $\omega$ of finite arity
$I_{\omega}$.

Further, let us present the excerpt, page 33.

\textmd{Let $\mathcal{K}$ be a $T_i$-prevariety. If for any
$\omega \in \Omega$ it is exist a finite set $B$ and a
neighborhood  $O\omega$ of point $\omega \in \Omega$ such that
$I_{\omega^{\prime}} \subset B$ for all $\omega^{\prime} \in
O\omega$, then $\mathcal{K}$ is called locally finite arity. On
analogy shall speak about point finite arity and locally finite
arity systems of operations $\Omega$.}

Let $G \in \mathcal{K}$ and let the mapping $\Psi: \Omega \times
G^A$ be continuous by the definition of algebra with continuous
signature. According to item 1d) the topology $\mathcal{T}$ on
$\Omega$ does not change the topological and algebraical structure
of topological algebra. If to consider $\mathcal{T}$ as
anti-discrete topology, then from the definition of locally finite
arity systems it follows that this notion coincides with the
notion of system with basic operations of finite arity.

2e). The notion of point finite arity $\Omega$-operation coincides
with notion of $\Omega$-operation with finite arity of basic
operation from $\Omega$, and the notion of locally finite arity
topological $\Omega$-algebra coincides with notion of topological
$\Omega$-algebra with finite arity of basic operation from
$\Omega$.

The authors claim on page 33 that bringing all operations to the
same arity plays a positive role for the proof of many facts. The
notion of  locally finite arity of prevariety is a basis for this paper’s
theorems and results. Is necessary to add to this the notions of $T_i$-prevariety,
the notions of algebraic and topological free algebras of $T_i$-prevariety. The
paper generalizes many known results, corrects the incorrect results of other
authors (see, items 2h), 2i)) and also gives valuable indications how to
investigate the topological algebras (see item 2j)). We are limited to
the volume of this paper, therefore we will not analyze the above-mentioned.
However, I think that it is enough to ascertain items 2a) - 2e).

2f). The work \cite{DC} consists only of erroneous and senseless
statements. It is not a scientific, but rather an anti-scientific
paper. To confirm the last statement, we will present some
excerpts (pag. 36).

 \textmd{Let is given a $T_i$-prevariety $\mathcal{K}$ with
continuous system of $A$-ary operations $\Omega$. For every
topological space $X$ it is exist a free topological algebra
$(F_{\mathcal{K}}(X), \sigma)$ (Theorem 1) and a free algebra
$(F_{\mathcal{K}}^a(X), \sigma)$ (Theorem 2). The topology on
$F_{\mathcal{K}}(X)$ denote by $\mathcal{T}_F$. There are
following questions.}

\textmd{Question  1. How depends the topology $\mathcal{T}_F$ from
topology of space $X$?}

\textmd{Question 2. Under what conditions the mapping $\sigma$ is
a homeomorphic embedding?}

\textmd{Question  3. Whether is the mapping $i_X:
F_{\mathcal{K}}^a(X)\rightarrow F_{\mathcal{K}}(X)$ an algebraic
isomorphism?}

\textmd{These questions were formulated A. I. Mal'cev \cite {Mal3}
for case, when  $\Omega$ is a finite discrete space.}

2g). The statement is senseless.

2h). \textmd{As will be shown further, answers to questions 2 and
3 are positive for a wide class of prevarieties and spaces, in
particular, and case of completely regular space $X$. Meanwhile \cite{Bur2}
presents a negative answer to question 3 even for the case $X$ and $T_2$-prevariety
$\mathcal{K}$ (the completely regular space and $T_2$-space are the same.) However a
careful analysis of the statements from \cite {Bur1} shows that the author interprets
incorrectly the problem...  The author \cite {Bur2} has a wrong approach to the A. I.
Mal'cev problem, which resulted in other mistakes. Some of the author’s statements about
$\beta$-classes are in particular wrong. It follows from the results obtained in paper (\S
6)that any non-trivial prevariety is a $\beta$-class in the sense of A. I.
Mal'cev \cite{Mal3}. Therefore this concept, from our point of
view, to examine is inexpedient, as a case non-trivial
prevarieties does not represent interest.}

2i). (page 51). REMARK. Theorems 7 and 8 generate the lemma and its consequence,
formulated in \cite[page 186]{Mal3}.
The basic result of M. S. Burghin  \cite{Bur2} is false.

2j). (page 39). REMARK. It is possible, following directly A. I.
Mal'cev \cite {Mal3}, to construct sequence enclosed topologies
$X_{\lambda}$ using polynomials, though the last definition, when all operations
from $ \Omega $ are $A $-ary, is very cumbersome. We define the initial topology  $\mathcal{T}_0$
otherwise, not using polynomials. It is possible to show that the sequence of topology,
constructed by us, contains the sequence of A. I. Mal'cev as a co-final sequence.
In particular, we receive the initial topology $X_0 $ at a certain step $\lambda_0$.

2k). Excuse my immodesty, but such ingenious ideas are impossible
to comment.

2l). The authors express their gratitude to the seminar leader,
professor V. I. Arnautov, and all its participants for
valuable remarks made during the discussion of this work.

Item 2l) is similar to item 1f). The authenticity of 2l) follows
from items 2a) - 2k). Item 2l) also refers to other papers, published by Academician
Choban M. M. in journal
"Matematicheskie issledovania" \cite{Chob31}, \cite{Chob32},
\cite{Chob32}. We present literally the excerptions a) -- h) from
\cite{Chob31} and analyze them briefly.

a). \textmd{A space $X$ is zero-dimensional if $\dim \beta X =
0$.} (pag. 121)

3a). This statement is not correct. \cite[example 6.2.20]{Eng}
presents a zero-dimensional space $X$ with $\dim \beta X \neq
0$. This excerpt is copied mechanically from another paper. This
statement is correct for non-empty class of locally compact
paracompacts, for non-empty class of  compacts \cite[Theorem
6.2.9, Corollary 6.2.10]{Eng}.

3b). A topological space $X$ is called zero-dimensional, $\dim X =
0$, if $X$ is a $T_1$-space and posed a basis from open-closed
sets. Any discrete space is zero-dimensional. A topological space
$X$ is zero-dimensional if and only if $\dim X = \text{ind} X = 0$
\cite[pag. 562]{Eng}.

b). \textmd{Let us consider a non-empty class $\mathcal{K}$ of
topological groups and  emphasize some properties, which can have
such class:}

\textmd{I. Closed with respect to direct topological products.}

\textmd{II. Closed with respect to subgroups.}

\textmd{III. Closed with respect to imagines of continuous
homomorphisms "on".}

\textmd{IV. If $G \in \mathcal{K}$ and with respect to topology
$\tau$ the pair $(G, \tau)$ is a topological group, then $(G,
\tau) \in \mathcal{K}$.}

ń). \textmd{A class $\mathcal{K}$ with properties I - III will
called variety. A class $\mathcal{K}$ with properties I - IV will
called primitive class according to \cite[pag. 425]{Mal2},
\cite[pag. 23]{Mal3}.}

The reference to Mal'cev's papers is inaccurate. In these
papers a primitive class of algebraic systems is defined as
totality of all algebraic systems having same basic operations of
finite arity and all elements which satisfy some system of
identities of form $f = g$, where $f, g$ are a polynomials.
Clearly, that the identities of form $f = g$ are quasiatomic
formulas.

In \cite[pag. 269]{Mal1} a variety of algebraic systems is defined
as totality of all algebraic systems having the same basic operations
of finite arity  and  satisfying the same system of quasiatomic
formulas of form $f = g$. Further, by \cite{Mal2}, \cite{Mal3} a
topological primitive class is called a totality of topological
algebraic systems, which make a primitive class in algebraically
sense.

3c). The notion of topological primitive class defined in item a)
is incorrect. The reference from item b) to Mal'cev's
papers is inaccurate. The notions of topological (algebraical)
primitive class and topological (algebraical) variety of algebraic
systems, defined by A. I. Mal'cev, coincide. By \cite[Theorem VI.13.1]{Mal1} or
Proposition 2 of that paper, these  notions for topological groups coincide
with a class of topological groups that satisfy conditions I, II, III of item a).

d). \textmd{Let $\mathcal{K}$ be a non-empty class of topological
groups. By $\mathcal{M}(\mathcal{K})$ we denote least variety,
containing the class $\mathcal{K}$ and by
$\mathcal{P}(\mathcal{K})$ we denote least primitive class,
containing  the class $\mathcal{K}$. Always
$\mathcal{M}(\mathcal{K})\subset \mathcal{P}(\mathcal{K})$. The
classes $\mathcal{M}(\mathcal{K})$, $\mathcal{P}(\mathcal{K})$ are
supplied with same "set" of abstract groups, but each group in
class  $\mathcal{P}(\mathcal{K})$ is supplied with the large
number of topologies.}

\textmd{LEMMA 2. Let $\mathcal{K}$ be a non-empty class of
bicompact groups. Then: 1. $\mathcal{P}(\mathcal{K}) \neq
\mathcal{M}(\mathcal{K})$. ... }

\textmd{PROOF. ... In $\mathcal{M}(\mathcal{K})$ it exists a
non-finite group $G$. Let $\tau$ denote the discrete topology on
$G$. Then $(G, \tau) \in \mathcal{P}(\mathcal{K})$... Lemma is
proved.}

We presented a fragment of Lemma 2 with proof to show the sense of
condition IV for introduced notion of primitive class. Takes place

3d). In definition of primitive class the condition IV is
senseless. It contradicts conditions I, II, III. Does not
exist neither  class $\mathcal{K}$ of topological groups such that
the primitive class $\mathcal{P}(\mathcal{K})$ exists.

Really, we assume that $\mathcal{M}(\mathcal{K}) \subset
\mathcal{P}(\mathcal{K})$. Let $(G, \tau_1) \in
\mathcal{M}(\mathcal{K})$ be a non-discrete topological group and
let $\tau$ denote the discrete topology on $G$. Similarly to Lemma
2 $(G, \tau) \in \mathcal{P}(\mathcal{K})$ by condition
IV. Let $(G, \tau)_F$ denote the free discrete group with free
generators of variety of groups in Mal'cev sense, generated by
$(G, \tau)$ \cite[Corollary VI.13.4]{Mal1}. The congruences of
groups are permutable. Then by \cite[Theorem 3]{Mal3}
the identical mapping of set $G$ induces a homeomorphism between
discrete space $(G, \tau)_F$ and non-discrete space $(G, \tau_1)$.
We get a contradiction. Hence the primitive class
$\mathcal{P}(\mathcal{K})$ does not exist.

3e). The basic object of researches in Lemmas 4, 7, Theorems 1, 2,
Propositions 1, 2, Corollaries 5, 6, 7 is the primitive class of
topological groups. Hence, by item 3d) the listed
statements are erroneous, without sense.

e) \textmd{DEFINITION 2. We consider a class of topological groups
$\mathcal{K}$ and a set $X$. The group $F \in \mathcal{K}$ free
topological group of space $X$ with respect to class
$\mathcal{K}$, if:}

\textmd{1. $X$ is a subspace of space $F$.}

\textmd{2. $X$ topologically generate $F$.}

\textmd{3. For any continuous mapping "on" $\varphi: X \rightarrow
G$, where $G \in \mathcal{K}$, it is exists a continuous mapping
$\tilde{\varphi}: F \rightarrow G$ such that $\tilde{\varphi}/X =
\varphi$.}

\textmd{If for a space $X$ in a class $\mathcal{K}$ a free
topological group exists, then this group we denote by
$F_{\mathcal{K}}(X)$.}

\textmd{From Lemma 3 it follows}

\textmd{COROLLARY 1. Let $\mathcal{K}$ be a non-empty class of
zero-dimensional bicompact groups and $\mathcal{M} =
\mathcal{M}(\mathcal{K})$. For a space $X$ the group
$F_{\mathcal{M}}(X)$ exists if and only if $\text{ind} X = 0$.}

\textmd{Thus, the object $F_{\mathcal{K}}(X)$ does not exists in
any class $\mathcal{K}$, even in any variety. With other party,
the Swierczkowsti's Theorem \cite{Sw} assert that object
$F_{\mathcal{K}}(X)$ exists for any  completely regular  set $X$,
if $\mathcal{K}$ is a primitive class. Hence, the topologically
notion of primitive class more flexibly, than notion of variety.
In some works on topological algebras these notions are
identified, what it is no correct (see, \cite{Bur1}). It can to
result and to erroneous results. For example, the Theorems 1 - 4
and them corollaries from \cite{Bur2} are erroneous.}

3f). The reference of Swierczkowsti's Theorem \cite{Sw} is inaccurate
by item 3d).

According to Theorem 8 of this paper and item 3b), the Corollary
1 is also senseless. Hence, the assertion that the Theorems 1 - 4 and
their corollaries from \cite{Bur2} are erroneous is not correct (see, also, items 2i), 2k)).

3g). Similarly to Corollary 1, according to Theorem 8 of this
paper and item 3b), the following excerpt (page 133) is also a gross mistake.
\textmd{Let a class $\mathcal{K}$ satisfies to
conditions I, II of item b), i.e. is a prevariety. In such case
for any space $X$ exists $F_{\mathcal{K}}(X)$ as soon as
$\text{ind} X = 0$.} Theorems 5, 6, 7 are presented and they claim that they
can be proved with the help of the presented excerpt and methods from $\S 2$,
consisting only of mistakes. Hence, the Theorems 5, 6, 7 are not proved,
actually they are erroneous.

d). (pages  125, 126). \textmd{The Theorem 1 is more general, than
results from Bel'nov's work \cite{Bel'nov}. In particular, from
Theorem 1 it follows}

\textmd{COROLLARY 2. Every metrizable topological group is a
quotient group of some zero-dimensional metrizable topological
group.}

e). \textmd{The Corollary 2 gives the complete answer to one A. V.
Arhangel'skii question \cite{Arh6} and independently is proved by
Bel'nov V. K \cite{Bel'nov}.}

f). Corollary 7 gives an answer to other question from Arhangel'skii's
work \cite{Arh6} (page 127).

3i). As noted in item 3e), Theorems 1, 2 and Propositions 1, 2 are
inaccurate. Therefore, their Corollaries 2 -- 7 are also inaccurate. Hence,
the Arhangel'skii questions from items e), f) are not solved.

3j). On a deceit, on senseless of item e) has got the Professor A.
V. Arhangel'skii. Paper \cite[page 1037]{Arh4} quotes the false statement e).

g). Items a) -- f), 3a) -- 3i) contain the analysis of the entire
work \cite{Chob31}, except for  the $\S 3$. $\S 3$ uses the
researches conducted in Choban's work \cite{Chob32}, \cite{Chob33}.
We will not analyze the authenticity of the last works’ results.
We will only note that the material and the proofs presented are not mathematical.
This confirms  (see, for example, \cite[page 118]{Chob33}),
that these works  generalize many known results, proved earlier by various authors.
The same is true about paper \cite{Chob31}. We present some excerpts from this paper
(pag. 130, 131). \textmd{... From above mentioned reasonings (27 lines) and results of
work \cite{Mich}, \cite{Chob32} it follows. THEOREM 4. .... The statement 4 of Theorem 4
supplements many statements from \cite{Chob32}.} It is followed by COROLLARIES 8 --  13
without proof. Let us present another excerpt.

h).  \textmd{The Corollaries 8, 9, 11, 12 generalize many
statements of Michael's paper \cite{Mich}, of Gleason's paper
\cite{Glea}, of Bartle's and Grawes's work \cite{Bart}, of
Kenderov's paper \cite{Kend} ... The Corollary 13 generalizes many
statements of Hofman's and Mostert's  monograph \cite{Hofm} and
known Iwasawa Theorem \cite{Iwas}.}

I think that the comments to items g), h) are clear.

3k). To result from items a) -- h), 3a) -- 3j) to ascertain, that
 the item 2f) concern and to papers  \cite{Chob31},
\cite{Chob32}. \cite{Chob33}. The   paper \cite{Chob31}   consist only from
erroneous nonsense statements. \cite{Chob31} is an anti-scientific paper.

4a). The senseless and incorrect notions and results from the
above analyzed works \cite{Chob3}, \cite{DC} about operations and
arity of operations (items 1a), 1b), 2a)), about continuous
signature (items 1m), 1e), 1h), 2b)), about various definitions of
quasivarieties (items $1e_1)$ -- $1g_3)$), about existence and the
interrelations of various definitions of free topological
algebraic systems (items 1h) -- 1l), 2b) -- 2d)) and many others are
used in other works, which implies incorrect and many times senseless
results of these works. We present only some examples. But before
that let us present some well-known facts.

$4a_1)$. We present the following literally.

THEOREM 3 (\cite{Mal3}).  Let $\mathcal{A}$ be a topological algebra with defining space
$X$ and with defining relations $S$ in primitive class (variety), defined by identities
$V$ and let  $\mathcal{B}$ be a topological algebra with defining space
$Y$ and with defining relations $T$ in primitive class,  defined by identities
$W$. We suppose that such a continuous mapping $\tau$ of space $X$ on space $Y$ is
given that every equality from $S$, written for the corresponding elements $Y$,
is contained in Ň, and let every identity of $S$ be contained in $T$.
Then the mapping $\tau$ can be continued to continuous mapping $\tau^{\star}$
of algebra $\mathcal{A}$ on algebra $\mathcal{B}$. This homomorphism  will be
open if in the class, defined by identities
$V$ the congruences are permutable, $\mathcal{B}$ topologically containing
$Y$ the mapping $\tau$ of $X$ on $Y$ is homeomorphism.

$4a_2)$. Every topological algebra of primitive class $\frak K$ is an image of
continuous homomorphism of free algebra of class $\frak K$ according to Theorem 3.
But to construct a topological quotient algebra the homomorphism should be
not only continuous, but also open.
One such sufficient condition is presented in Theorem 3. But in general
this is not correct. \cite[pag. 177]{Mal3} contains an example of a
continuous homomorphism of free semigroup which is not open.

4b). We analyze briefly papers \cite{Chob9}, \cite{Chob19}. These
papers consider the topological algebras of continuous signature
and no limitations are imposed on the arity of  basic operations.
In such a case Proposition 1 and Lemma 2 of the given paper cannot
be applied. How could we speak about the soundness of these
papers’ results if the author Academician Choban M. M. is unaware
of the basic notions of algebras, which can be found in any manual
and monograph on algebra \cite{Mal1}, \cite{Con}, \cite{Kurosh}:
terms, polynomials, identities, varieties, replica (\cite[pag.
185, Definition 5.5]{Chob9}, \cite[paragraph 2.4; 3.1]{Chob19},
signature, arity (\cite[Problems 4, 4]{Chob9}, direct and
Cartesian products (\cite[pag. 14]{Chob19} and many others. We
will present in form of excerpts from \cite{Chob9}.

$4c_1)$. \textmd{The discrete sum $E = \oplus\{E_n: n \in N = \{0, 1,
2, \ldots\}\}$ of topological spaces $\{E_n: n \in N\}$ is called
the continuous signature.}

\textmd{\textbf{Definition  1.1} An $E$-algebra or a universal
algebra of the signature $E$ is a family $\{G, e_{nG}:  n \in N\}$
for which:}

\textmd{(1) $G$ is a nonempty set.}

\textmd{(2) $e_{nG}: E_n \times G^n \rightarrow G$ is a mapping
for every $n \in N$.}

\textmd{Let $i \in \{-1, 0, 1, 2, 3, 3.5\}$. If a topological
$E$-algebra $G$ is a $T_i$-space, then $G$ is called
$T_i$-$E$-algebra.}

\textmd{If $p \in E_0$ and $G$ is an $E$-algebra, then
$e_{0G}(\{p\} \times G^0) = 1_p$. If $n \in N$ and$ q \in E_n$,
then $q: G^n \rightarrow  G$, where $q(x_1, \ldots,  x_n) =
e_{nG}(q, (x_1, \ldots,  x_n))$,  is an operation of type $n$ on
$G$.} (Correct: operation of arity $n$).

$4c_2)$. \textmd{The set of terms $T(E)$ is the smallest class of the
operations on the $E$-algebras such that:}

\textmd{(1) $E \subseteq T(E)$ and $e_G \in T(E)$, where $e_G(x)
= x$ for every $x \in G$.}

\textmd{(2) If $n > 0$, $e \in  E_n$ and $u_i,\ldots , u_n \in
T(E)$, then $e(u_i, \ldots, u_n) \in  T(E)$.}

$4c_3)$. \textmd{Let $1 \leq m < n$, $N_m = \{1, 2,\ldots , m\}$
and $h: N_n \rightarrow,  N_m$ (onto) be a mapping. The operation
$\omega: G_n \rightarrow  G$ of type $n$ and $h$ generates the
operation $\psi: G_m \rightarrow G$, where $\psi(x_1, \ldots, x_m)
= \omega(x_{h(1)}, \ldots, x_{h(n)})$,  and $\psi$  is called an
$h$-permutation of the operation $\omega$.} (it is not allowed to
change the arity of basic operation, i.e. the signature of given
algebra).

$4c_4)$. \textmd{The set of the polynomials $P(E)$ or of the derived
operations is the smallest class of operations on $E$-algebras
such that:}

\textmd{(1) $T(E) \subseteq P(E)$.}

\textmd{(2) If $f \in P(E)$  and $g$ is an $h$-permutation of $f$,
then $g \in P(E)$.}

$4c_5)$. \textmd{If $\omega$ and $\psi$ are polynomials of types
$n$ and $m$, then the form \break $\omega(x_1, \ldots ,x_n) =
\psi(y_1, \ldots, y_m)$ is called an identity on the class of
$E$-algebras.}

$4c_6)$.  Any comments are irrelevant.  Similarly to item 1q), the
presented definitions are a delirium, senseless statements.

$4d_1)$. Sections 1, 2 contain the definitions of $T_i$-quasivariety, of
$T_i$-variety, of complete $T_i$-quasivariety, of complete
$T_i$-variety of topological $T_i$-$E$-algeb\-ras similarly to item $1e_1$.

$4d_2)$. For example, a class $\frak K$ of topological $T_i$-$E$-algebras
is called a $T_i$-quasivariety if $\frak K$ is closed with respect
to subalgebras and Tychonoff product.

$4d_3)$.  It contains  the definitions of topologically free
algebra $(F(X, \frak K), i_X)$ of a space $X$ in a
$T_i$-quasivariety $\frak K$ and of algebraically free algebra
\break $(F^a(X, \frak K), j_X)$ of a space $X$ in a
$T_i$-quasivariety $\frak K$ as in item $1g_3)$, which are
presented in item 6c). Almost all results of paper \cite{Chob9}
are related to free algebras or the free algebras are a crucial
instrument to prove the statements of the analyzed paper. Item
$6d_1)$ shows that this notions are not corrects, senseless.
Hence, all results from the analyzed paper \cite{Chob9} are false,
many of them are absurd.

 We present the following literally.

$4e_1)$. \textbf{Definition 5.5.}  Let $E$ be a continuous signature.
The class $K$ of $E$-algebras is
called a Mal’cev class if there exists a polynomial $p(x, y, z)$ such that the
equations $x = p(y, y, x) = p(x, y, y)$
hold identically in $K$. The polynomial $p(x, y, z)$ is called a Mal’cev polynomial.

$4e_2)$. Lemma 5.8 (Mal’cev \cite{Mal2} for finite $E$). Let
$g: A \rightarrow B$ be a homomorphism of a
topological $E$-algebra $A$ with the Mal’cev polynomial onto an $E$-algebra $B$. Then the
$E$-algebra $B$ equipped with the quotient topology is a topological $E$-algebra and the
quotient mapping $g$ is open.

Proof. Let $p(x, y, z)$ be a Mal’cev polynomial on the algebra $A$
and $V$ be an open set in $A$. We prove that $U = g^{-1}(g(V))$ is
open in $A$ (this statement was proved by Mal’cev \cite{Mal2}. Fix
a point $x \in U$. Then $g(x) = g(y)$ for some $y \in V$. . .

$4e_3)$. The Definition 5.5 is unclear, it does not correspond to the generally
accepted notion of Mal'cev algebra (see Theorem  1 of this paper) according to
item $4c_4)$ about the senselessness of notion of polynomial.

The expression $x = p(y, y, x) = p(x, y, y)$ are identities, rather than equations.

The referral to \cite{Mal2} is the Theorem 10 from \cite{Mal2}. It cannot be
used in the proof of Lemma 5.8 according to item $4e_3)$. Hence, Lemma 5.8 is not
proved. From here it follows that all results of
Sections 5, 6 are false.

$4e_4)$. Items 4a), 4b) -- $4e_3)$ confirm only partially the following statement:
paper \cite{Chob9}
consist only of erroneous,  senseless  statements and introduced notions.
\cite{Chob9}  is an anti-scientific paper.

$5a)$. Moreover, the results from \cite{Dm}, \cite{Chir3}, presented in
Sections  5, 6 from \cite{Chob9} (see item $4c_8)$) are false. Consequently,
the results of works \cite{Dm}, \cite{Chir3} are false.

Now let us analyze paper \cite{Chob19}.

$6a$. Section 2 contains the same senseless notions as in
\cite{Chob9}: algebra of continuous signature ($E$-algebra),
various notions of quasivariery, polynomials, terms, identities,
free algebras. The list for the analyzed paper is complemented
with the following notions: constants (see \cite[pag. 335]{Mal1},
supports of elements, replica, saturated classes of algebras and
others. We will present them literally.

$6b$. \textmd{Let $D \subseteq E$ and $D_n = D \cap E_n$ for any $n \in N$. Then $D$ is a continuous
signature and every topological $E$-algebra is a topological $D$-algebra, where
$d_{nA} = e_{nA} | D_n \times A^n$ for any $n \in N$. Moreover, every $E$-homomorphism $f:
A \rightarrow B$ of $E$-algebras $A$ and $B$ is a $D$-homomorphism, too.}

$6c$. \textmd{\textbf{Definition 2.2.1.} \textit{Fix $i \in \{-1;
0; 1; 2; 3; 3,5\}$, a $T_i$-quasivariety $K$ of topological
$E$-algebras and a topological $D$-algebra $X$. Then:}}

\textit{(T) A couple $(F(X, K, D), i_X)$ is called a topological $T_i$-free $E$-algebra of the
$D$-algebra $X$ in the class $K$ if the following conditions hold:}

\textit{((1) $(F(X, K, D) \in K$ and $i_X: X \rightarrow (F(X, K, D)$ is a continuous
$D$-homomorphism;}

\textit{(2)  the set $i_X(X)$ generates $(F(X, K, D)$.}

\textit{(3) for each  $D$-homomorphism $f: X \rightarrow G \in K$
there exists a continuous $E$-homomorphism $\overline{f}: (F(X, K,
D) \rightarrow G$ such that $f(x) = \overline{f}(i_X(x))$ or every
$x \in X$.}

\textit{(A) A couple $(F^a(X, K, D), j_X)$ is called a  $T_i$-free
$E$-algebra of the $D$-algebra $X$ in the class $K$ if the
following conditions hold:}

\textit{((4) $(F(X, K, D) \in K$ and $i_X: X \rightarrow (F(X, K, D)$ is a
$D$-homomorphism;}

\textit{(5)  the set $j_X(X)$ generates $(F^a(X, K, D)$.}

\textit{(6) for each  $D$-homomorphism $f: X \rightarrow G \in K$ there exists a continuous
$E$-homomorphism $\overline{f}: (F^a(X, K, D) \rightarrow G$ such that $f(x) = \overline{f}(i_X(x))$
for every $x \in X$.}

$6d$. \textbf{Theorem 2.2.3} \textit{Let $K$ be a
$T_i$-quasivariety of topological $E$-algebras. Then for every
topological $D$-algebra $X$ there exist:}

\textit{- a unique $D$-free algebra $(F^a(X,  K, D), j_X)$;}

\textit{- a unique topologically  $D$-free algebra $(F^a(X,  K, D), i_X)$}

\textit{- a unique continuous homomorphism $k_X: F^a(X, K, D) \rightarrow
F(X, K, D)$ and $i_X = k_X \cdot j_X$.}

\textbf{Proof.} Let $tau$ be an infinite cardinal and $tau \geq
\max\{|X|, |E|\}$. Then the collection $\{f_{\xi}: X \rightarrow
G_{\xi} :\xi \in L\}$
 of all $D$-homomorphism of $X$ into the $E$-algebra $G \in  K$ of the
 cardinality $\leq \tau$ is a non-empty set. We consider the
diagonal product $j_X : X \rightarrow G = \prod \{G_{xi} : \xi \in
L\}$, where $j_X(x) = (f_{\xi}):x) : \xi \in L)$, and denote by
$F^a(X, K,D)$ the $E$-subalgebra of $G$ generated by the set
$j_X(X)$. If $f : X \rightarrow A \in K$ is a $D$-homomorphism,
then for some $\xi \in L$ we have $G_{\xi} = d(E; f(X))$ ($d(E,
f(X))$ denote the $E$-subalgebra generated by $f(X)$) and
$f_{\xi}$. The natural projection $\overline{f} : F^a(X, K, D)
\rightarrow G_{\xi} \subseteq A$ is the homomorphism generated by
$f$.

If $M = \{\xi \in L: f_{\xi} \quad \text{is continuous}\}$, then
we consider the diagonal product $i_X : X \rightarrow G^{\prime} =
\prod\{G_{\xi} : \xi \in L\}$ and $(X, K, D)$ is the
$E$-subalgebra of $G^{\prime}$ generated by the set $i_X(X)$. For
some $\xi \in  M$ we have $F(X, K, D) = G_{\xi}$ and $i_X =
f_{\xi}$. The proof is complete.

$6d_1)$. The uniqueness of $D$-free algebra $(F^a(X, K, D),  j_X)$
and topologically $D$-free algebra $(F^a(X, K, D),  i_X)$ is not
proved at all. Proofs of similar statements can be found in
\cite{Mal3} and Corollary 3 of this paper. The proof of the third
statement is obviously wrong. The very first sentence of Proof is
not correct. The Proof is meaningful only in case when
$|\prod\{G_{\xi} : \xi \in L\}| \leq \tau$.

$6e)$. From items $6a)$ – $6d_1)$ it follows that paper
\cite{Chob19} can be appraised in the same manner as item $4e_4)$
for \cite{Chob9}, i.e. \cite{Chob19} it is an anti-scientific
paper.

7a). Moreover, in \cite{Chob19} is mentioned that \cite{Chob191},
\cite{Chob192} were examinees the same questions about free
objects, various types of equivalence that [?]. From here we
conclude that the results from \cite{Chob191} and \cite{Chob192}
are false.

\centerline{\textsf{Homogeneous algebra, rectifiable space,
Mal'cev operation}}

In terms of scientific importance the paper \cite {Chob2}, published in
journal: \textsl{Serbica, Bulgaricae mathematicae publications},
does not differ from the previous paper \cite{Chob3}. It, in
general, consists of easily noticed gross errors and incorrect or senseless  statements.

In the beginning of the paper, following \cite{Chob3}, presents the
definitions of topological algebra of continuous signature (the
basic operations have a finite arity), calls  $E$-algebra, of
quasivariety $\mathcal{K}$ of $E$-algebras as class of all
$E$-algebras of given continuous signature with conditions 1m. 2m.
For quasivariety $\mathcal{K}$ presents the Definition 8 and the
Property 2.

Further, for a quasivariety $\mathcal{K}$ of topological
$E$-algebras present Theorems 01, 02 as analogical to Theorems 2,
3 from \cite{Mal2}, proved for primitive classes (variaties) of
algebras. Claims that the proofs of Theorems 2, 3 for primitive
classes of algebras are valid for Theorems 01, 02. But it is not
true. Really, the existence of free algebras for primitive classes
is used essentially to prove Theorems 2 and 3. However, for
quasivarieties it is not true. It follows from Theorem 7 of this
paper, if to ignore the requirement of continuous signature by
item 1e).

The paper contains only three statements with proof, specifically Theorems
1.1, 2.1, and 3.1. Theorem 1.1 investigates such regular $\sigma$-compact
$E$-algebras for which the signature $E$ is a
$\sigma$-compact space. By item 1e) the requirement of
signature can be ignored. Then, Theorem 1.1 and its proof literally coincides
with the appropriate result from \cite{Prod}.

Further, let us consider the excerpt from \cite{Chob2}.

\textmd{Theorem 2.1. \textit{Let $\mathcal{K}$ be a quasivariety
of topological $E$-algebras and let the free algebra $F(D,
\mathcal{K})$ be discrete. Then the following conditions are
equivalent:}}

\textmd{1. \textit{Any algebra $A \in \mathcal{K}$ has permutable
congruences;}}

\textmd{2. \textit{There exists a term $p(x,y,z)$ such that the
equations $p(x,y,y) = \break p(y,y,x) = x$ hold identically in
$\mathcal{K}$;}}

\textmd{3. \textit{Any congruence in an algebra of $\mathcal{K}$
is open.}}

We ignore the requirement about quasivariety and about $E$-algebra, because,
as mentioned above, they are fictitious. Let's assume that $\mathcal{K}$
is a primitive class of algebraic systems. Then the equivalence 1)
$\leftrightarrow$ 2) is the famous Mal'cev's Theorem, the Theorem 1.
The implication 1)$\rightarrow$ 3) is not less remarkable Mal'cev's result, the
Theorem 10 from \cite{Mal2} (the congruence, considered in Theorem 10 are
called complete congruences, but not open congruences). See also Theorem 4
and Corollary 2 of this paper.

Let the implication 3) $\rightarrow$ 1) hold. Then by
\cite[Proposition 2.4.9]{Eng} to equivalent conditions 1), 2), 3)
of Theorem 2.1 it is possible to add and equivalent conditions 1),
2), 3) of Corollary 3, and also the items 4), 5) of Corollary 3.
In result the Mal'cev's results mentioned above is very much
amplified. But, naturally, it not so, it is utopia. From here, it
is even possible to see by intuition that the implication
3)$\rightarrow$ 1) is false. It is also possible to see that the
author of paper  \cite{Chob2} has not understood the algebraical
and topological essences of Mal'cev Theorems.

Actually, without going into too many details we will note that the proof of
implication 3) $\rightarrow$ 1) of Theorem 2.1 is false.
Particularly, to prove this implication the author uses essentially Properties 1, 2
from \cite{Chob3}, i.e. the existence of free topological algebras in any quasivariety.
As shown during the analysis of paper \cite{Chob3}, Properties 1 and 2 cannot be regarded
as proved. More specific, their proofs consist of a set of gross errors, see items 1j) - 1l).
Further, after Theorem 2.1 one may find erroneous statement presented without any proof as
Corollaries 2.1, 2.2, 2.3 and Remarks 2.2,
2.4. It is claimed that they can be proved with the help of the false Theorem 2.1.

Before starting the analysis of Theorem 3.1 we will present an excerpt from \cite{Chob2}.

$\ll$\textmd{If for the $E$-algebra $A$ there exist two
polynomials $p(x, y)$ and $q(x, y)$ with equations $p(x, x) = p(y,
y)$, $p(x, q(x, y)) = y$, then $A$ is called a homogeneous
algebra}$\gg^{(5)}$. $\ll$\textmd{We have $q(x, p(x, y)) =
y$.}$\gg^{(2)}$ $\ll$\textmd{Let \ldots $P_a(x) = p(a,x)$, $Q_a(x)
= q(a,x)$ \ldots Then \ldots $P_a$ and $Q_a$ are homeomorphisms
and $Q^{-1}_a = P_a$}$\gg^{(3)}$. $\ll$\textmd{An algebra $A$ is
homogenous iff $A$ is a biternary algebra}$\gg^{(4)}$. $\dots$

$\ll$\textmd{The mapping $g: X \times X \rightarrow X \times X$ is
a centralizing mapping with center $a \in X$ if $g$ is a
homeomorphism}$\gg^{(5)}$, \textmd{$g(x,x) = (x, a)$ and $\{x\}
\times X = g(\{x\} \times X)$ for every $x \in X$. The space $X$
is called a centralizable space if there exists a centralizable
mapping for $X$} (see \cite{Shap}).

\textmd{Theorem 3.1 The space $X$ is a centralizable space iff for
some two continuous binary operations $X$ is a homogeneous
algebra.}

\textmd{Proof. Let for the continuous binary operations $p$ and
$q$ the space $X$ be a homogeneous algebra. Then $h(x,y) = (x,
p(x,y))$ is a centralizing mapping with centre $p(x,x)$.
$\ll$\textmd{Let $h$ be a centralizing mapping with centre $a \in
X$. We put $g(x,y) = y$, $p(x,y) = g(h(x,y))$ and $q(x,y) =
g(h^{-1}(x,y))$}$\gg^{(6)}$. Then $p(x,x) = a$, and
$\ll$\textmd{$q(x, p(x,y)) = y$}$\gg^{(7)}$. $\Box$}

$^{(1)}$. The $E$-algebras are considered only here. By item 1e) a question arises:
what is the role of the continuous signature of $E$ in the definition of
homogeneous algebras? Further, an elementary
remark: the notions of equation and identity are totally different.
The expressions $p(x, x) = p(y, y)$, $p(x, q(x, y)) = y$ are not equations, but identities.

According to item 5) of Theorem 11 the identities $p(x,x) =
p(y,y)$,\break $p(x, g(x,y)) = y$,  $q(x, p(x,y)) = y$ transform a
non-empty space $A$ into a topological right loop $(A, \cdot,
\backslash, e)$, where $e = p(x,x)$, $x \cdot y = q(x,y)$, $x
\backslash y = p(x,y)$. The solution of equation $a \cdot x = b$
follows from identity $p(x, g(x,y)) = y$ and its uniqueness
follows from identity  $q(x, p(x, y))  = y$. For a groupoid $(A,
\cdot)$ the existence and the  uniqueness of solution are
independents (see, \cite{Bruk}, \cite{Bel}, \cite{Smith1}). Then
on space $A$ the identity $q(x, p(x,y)) = y$ does not follow from
identity $p(x, g(x,y)) = y$ and the identity $p(x, q(x,y)) = y$
does not follow from identity $q(x, p(x,y)) = y$. Consequently,
these identities defined on topological space are independents.
Hence the assertion $^{(2)}$ is false. In this case the definition
of homogeneous algebra from this paper is not correct with respect
to the corresponding definitions from \cite{Usp1}, \cite{Usp2},
\cite{Gul}.

Under definition $^{(1)}$ of homogeneous algebra the assertion
$^{(3)}$ is false. But under the corrected definition the assertion
$^{(3)}$ is not an obvious fact. For the proof of $^{(3)}$ under such a
condition see the comment on assertion $^{(19)}$, presented below.

The assertion $^{(4)}$ is false. For more details see the comment on assertion $^{(11)}$.

$^{(5)}$, $^{(6)}$. As is mentioned in the beginning of this paper
the unique distinction between classes of centralizable spaces
and rectifiable spaces is that for rectifiable spaces in the
definition of centralizable spaces it is necessary to require that
the homeomorphism $g: X \times X \rightarrow X \times X$  be surjective.
Otherwise the operation $q(x,y) = g(h^{-1}(x,y))$ is not determined for all $x, y \in X$.

The assertion  $^{(7)}$ is not an obvious fact, but not proofs are given for
$^{(7)}$. Moreover, $^{(7)}$  does not correspond to definition $^{(1)}$.

Further, the following are included without any proof

\textmd{Lemma 3.1. Let $f: A \rightarrow B$ be a homomorphism of
the homogeneous topological $E$-algebra $A$ onto the $E$-algebra
$B$. Then for every $x, y \in B$ the spaces $f^{-1}(x)$ and
$f^{-1}(y)$ are homeomorphic.}

Lemma 3.1 holds when $A, B$ are topological groups.
According to Theorems 9 and 10, it holds when $A, B$ are topological loops.
But Lemma 3.1 is false when $A, B$ are homogeneous topological $E$-algebras.
We ignore the senseless statement that $A, B$ are $E$-algebras, i.e. we ignore
the condition of continuous signature. According to $^{(1)}$ and $^{(9)}$
below let's consider that the notion of homogeneous algebra is
defined correctly. In this case by item 5) of Theorem 11 $A, B$ are topological right loops.

Recall that a congruence is called correct, if it is defined
unequivocally by any adjacent class (\cite{Mal2}). Let $x, y$ be
arbitrary elements in topological right algebra $(A, \cdot,
\backslash, f)$. If $(A, \cdot, \backslash, f)$ is not a left loop
then from the definition of right loop it follows that there exist
such elements $x, y$ in $A$ that $T(x) \neq y$ for any translation
$T$ $(A, \cdot, \backslash, f)$. It means that $(A, \cdot,
\backslash, f)$ does not satisfy  a necessary condition for a
congruence of algebraic system to be correct \cite[pag. 10]{Mal2}.
Hence Lemma 3.1 is false. This statement also follows from item 6)
of Corollary 21.

From the above-stated it follows that the false  Lemma 3.1 was written
as follows. The expression "topological group" was replaced mechanically
with expression "homogeneous topological $E$-algebra" in true version
of Lemma 3.1 for topological groups.

The following can be noticed from the above-mentioned comments.

a). The listed properties (without proofs) for homogeneous
algebras defined by $^{(1)}$ are false. Theorem 3 is false for
centralizable spaces. As result, the proof of Theorem 3 represents a
chaotic set of senseless and unrelated statements.

b). The listed properties, except for ${(4)}$, and Theorem 3 hold for
the modified definitions of homogeneous algebras and of centralizable
spaces (called rectifiable spaces) from \cite{Usp1}, \cite{Usp2}, \cite{Gul}.
The first two papers were published before the analyzed paper was submitted for publication.

c). Unless knowing beforehand the results of papers \cite{Usp1},
\cite{Usp2}, \cite{Gul}, it is not possible to understand and verify the
soundness of the statements from the excerpt, in particular Theorem 3 and its proof.

d). In the beginning of the article it was mentioned that, according to [22],
the notion of rectifiable space was introduced at the seminar of Prof.
Arhangel’skii. I believe that the above-mentioned excerpt is an unsuccessful
compilation of statements, collected during the mentioned seminar but not
understood by the author himself.

e). The above excerpt, in particular Theorem 3 and its proof, cannot
be regarded as a mathematical research.

The content of the remaining part of the given paper covers the entire
spectrum of topological science, and probably, the entire set of topological
notions. It presents plenty of remarks and open questions, which are compiled
by the following principles.  In text of the
original, mainly on topological groups, the expression
$''$topological group$''$ is replaced by $''$homogeneous topological
$E$-algebra$''$ (for example, the Lemma 3.1) and it is stated that they
are proved by same methods as the original, using not only the incorrect,
but as shown above, the senseless  Theorem 2.1 and false Lemma 3.1.

Other principle underlying the listed remarks and open questions is
the incompetence of the author and groundlessness introduction of new
concepts. Let's confirm this by excerpts.

$\ll$\textmd{Let $E = E_3 = \{p\}$. We denote by $V$ the class of
all topological $E$-algebras with Mal'cev equations $p(x,y,y) =
p(y,y,x) = x$. By $W_1$ we denote the class of all algebras $A \in
V$ with the equation $p(y,x,y) = x$. Further, let $W_2$ designate
the class of all algebras $A \in V$ with the equation $p(x,y,x) =
x$ and $W = W_1 \cup W_2$.}$\gg$

There is not any connection between the definitions of classes $V,
W_1, W_2,  \break W$ and topology of topological $E$-algebras.
Further, Mal'cev identities is the correct form, and not Mal'cev
equations. The author has not understood the meaning of Theorem 1
(or Theorem 3), the importance of which is shown in Theorem 10 of
\cite{Mal2} (or Theorem 4, Corollary 2). There is not any
connection between Mal'cev equations of $E$-algebras and ternary
operation $\Psi(x,y,z)$ of Theorems 1, 3. It follows from item 4)
of Proposition 13. Therefore it is inadmissible to name the
identities $p(x,y,y) = p(y,y,x) = x$ as Mal'cev equations.
Consequently, if one follows this principle, then it is possible
to introduce and investigate classes of $E$-algebras, defined by
arbitrary operations and arbitrary identities, which obviously is
not correct.

As a consequence of the above-stated, we list the remarks and open questions
that are incorrect. Let us present several examples.

$\ll$\textmd{Proposition 4.2. If $\mathcal{K} \subset W$ is a
quasivariety, then $\mathcal{K} \subset W_1$ or $\mathcal{K}
\subset W_2$ }$\gg$. Without sense  and trivial assertion, at than
here requirement about quasivariety. The proposition is false.
This   follows from following example. Let $A$ be an $E$-algebra
consisting from one element. Then $A \in \mathcal{K}$, $A \in W$,
$A \in W_1$, $A \in W_2$.

$\ll$\textmd{Remark 4.3. If $X$ is a zero-dimensional metrizable
space, then: 1. $X \in W_1$ for some $E_{3X}$; 2. $X \in W_2$ for
some $E_{3X}$.}$\gg$. The incompetence follows from not knowledge
of definition of zero-dimensional  space. No any connection
between it notion and definitions of $W_1, W_2$.

$\ll$\textmd{Question 4.1. Let $A \in V$ be a pseudocompact space.
Is $A^2$ a pseudocompact space?}$\gg$. According to Proposition 2
$A \in A^2$. The requirements $A \in V$ and the pseudocompactly of
space are independents. The answer to this question is negative,
as the family of all pseudocompact spaces is not closed with
respect to product of spaces \cite[pag. 746]{Eng}.

$\ll$\textmd{Question 4.2. Let $A$ be a pseudocompact homogeneous
algebra. Is $A^2$ a pseudocompact space?}$\gg$.

$\ll$\textmd{Question 4.3. Let $A \in V$ be a compact space. Is it
true that \text{dim}$A$ = \text{ind}$A$?}$\gg$ Rash question.
Conditions $A \in V$ and compact space are independents. An
example of compact space $X$ for which \text{dim}$X$ $\neq$
\text{ind}$X$, see \cite[pag. 612]{Eng}. For every metrizable
compact spaces $X$ \text{dim}$X$ = \text{ind}$X$, see \cite[pag.
600]{Eng}.

 Now we pass to analysis the paper
\cite{AC1}. Theorems 2, 3, 4, 5, 9, Proposition 9, and
Corollary 6 have the following form:

\textmd{Suppose that $B$ is a compact Hausdorff space \ldots
Suppose further that $B = X \cup Y$, where $X$, $Y$ are
non-locally compact rectifiable spaces. Then \ldots}

\textmd{Proof. Clearly, $Y$ and $X$ are non-empty, since they are
not locally compact. Hence, $B$ is non-empty. \ldots}

According to Proposition 15 the spaces $X$, $Y$ are compact. Then,
obviously, they are locally compact (see, \cite[pag. 196]{Kelli}.
Consequently, Theorems 2, 3, 4, 5, 9,  Proposition 9 and
Corollary 6  are trivial, they hold only for $B = \emptyset$.
But even without this remark the proofs of these theorems are
erroneous. Let's note only that the proofs use essentially the Dichotomy
Theorem for remainders of rectifiable spaces from \cite{AC2}. When analyzing
paper \cite{AC2} we will show that in \cite{AC2} this theorem is not proved.

The expression $\ll$\textmd{It is obvious that if a $k$-gentle paratopological
group is a $k$-space, then this paratopological group is a topological
group}$\gg$ from the proof of Proposition 2 is not clear (and is incorrect).

Based on the above-mentioned we conclude that all statements from sections 2, 3
(Theorems 2, 3, 4, 5, Propositions 1, 2, Example 7) are either trivial, or
senseless and their proofs are erroneous. Let's also
note that in the beginning of section 2 a definition of rectifiable space is
presented from papers \cite{Chob1}, \cite{Chob2}. But it is not correct.
These papers do not define the rectifiable space; moreover, the term
$''$rectifiable space$''$ cannot be even found in these papers(see,
the analysis of these papers above). The circumstances how this definition
appeared are presented in the beginning of this article.

Further, section 4 begins with the following excerpt.

 $\ll$\textmd{\textit{A Mal’cev operation} on a space $X$ is a
continuous mapping $µ : X^3 \rightarrow X$ such that $µ(x, x, z) =
z$ and $µ(x, y, y) = x$, for all $x, y, z \in X$. A space is
called a \textit{Mal’cev space} if it admits a Mal’cev operation
(see \cite{Chob1}, \cite{Chob2},  \cite{Chob9}, \cite{Mal2},
\cite{Usp2}} $\gg^{(8)}$.

$\ll$\textmd{\textit{A homogeneous algebra} on a space $G$ is a
pair of binary continuous operations $p, q : G \times G
\rightarrow G$ such that $p(x, x) = p(y, y)$, and $p(x, q(x, y)) =
y$, $q(x, p(x, y)) = y$ for all $x, y \in G$}$\gg^{(9)}$. If the
above conditions are satisfied, then $\ll$\textmd{the ternary
operation $µ(x, x, z) = q(x, p(y, z))$}$\gg^{(10)}$ is a Mal’cev
operation (see [9, 10]).

$\ll$\textmd{(In [9,10] (see also [19]) it was proved that for an
arbitrary space $G$ the following conditions are equivalent:}

\textmd{1) $G$ is a rectifiable space;}

\textmd{2) $G$ is homeomorphic to a homogeneous algebra;}

\textmd{3) There exists a structure of a biternary algebra on
$G$.}$\gg^{(11)}$

A structure of a \textit{topological quasigroup} on a space $G$ is
a triplet of $\ll$\textmd{binary continuous operations $p, l, r :
G \times G \times G \rightarrow G$  such that $p(x, l(x, y)) =
p(r(y, x), x) = l(x, p(x, y)) = l(r(x, y), x) = r(p(y, x), x) =
r(x, l(y, x)) = y$, for all $x, y \in G$}$\gg^{(12)}$. If there
exists an element $e \in G$ such that $p(e, x) = p(x, e) = x$ for
any $x \in G$, then we say that $G$ is a topological loop and $e$
is the identity of $G$. $\ll$\textmd{Any topological quasigroup
admits the structure of a topological loop (see
\cite{Mal2})}$\gg^{(13)}$. $\ll$\textmd{If $e \in G$ and $p(e, x)
= x$ for any $x \in G$, then $x + y = p(y, x)$ and $x · y = r(y,
x)$ is a structure of a homogeneous algebra}$\gg^{(14)}$.
$\ll$\textmd{If $(G, ·)$ is a topological group with the neutral
element $e$, then the mapping $\varphi(x, y) = (x, x^{-1} · y)$ is
a rectification on the space $G$ with the neutral element $e$, and
the mappings $p(x, y) = x^{-1}·y$ and $q(x, y) = x·y$ form a
structure of homogeneous algebra on $G$. Therefore,
$\ll$\textmd{every topological quasigroup is a rectifiable space
}$\gg^{(16)}$ }$\gg^{(15)}$.

$^{(8)}$. The need to study the Mal'cev operation and Mal'cev space
follows from Mal'cev's paper \cite{Mal2}, which puts the foundation
of topological algebraic systems. The cited literature specifies
several papers \cite{Chob1}, \cite{Chob2}, \cite{Chob9}, \cite{Mal2},
\cite{Usp2} that study (use) these notions. We can add papers \cite{Usp1},
\cite{Usp3}, \cite{Gul}, as well as \cite{AC1}, \cite{ACM1}, \cite{ACM2},
\cite{ACM3}, where various generalizations are introduced. However, some
peculiarities arise during their consideration. The definitions of Mal'cev
operation and Mal'cev space are too general and it is impossible to apply to
them the fundamental results from \cite{Mal2}.

Theorem 1 characterizes an algebraic systems with permutable
congruences. For this purpose on such algebraic systems a ternary
operation \break $\Psi(x,y,z)$ is defined with identities (5),
i.e. ( a Mal'cev operation $\Psi(x,y,z)$ is defined. But in item
4) of Proposition 13 a loop is indicated, which admits a Mal'cev
operation  with non-permutable congruences. Such a situation
occurs because in Theorem 1 the ternary operation  should be a
polynomial (from the basic operations) of the given algebraic
system.

To avoid such misunderstandings some authors (see, for example,
\cite{Smith2}) do the following. Define a Mal'cev algebraic system as
algebraic system with permutable congruences. Then for such
Mal'cev algebraic systems it is possible to apply Theorem 1.
Moreover, for the Mal'cev algebraic systems defined in such a way it is
possible to apply the important  Theorem 10 from \cite{Mal2} (see
also the Theorem 4 of this paper): every quotient homomorphism
of Mal'cev algebraic system is open.

Some authors, \cite{Usp1}, \cite{Usp2}, apply Theorem 10 in the
following way: for Mal'cev operation $[ , ,]$, defined generally, they
consider only those homomorphisms $\varphi$ of Mal'cev space $(A,
[ , , ])$ which are consistent with Mal'cev operation $[ , , ]$:
$\varphi[x,y,x] = [\varphi x, \varphi y, \varphi z]$ for all $x,
y, z \in A$. Then from here (or Theorem 1 and \cite[Theorem
10]{Mal2}), Proposition 10  and the equivalence $1)
\leftrightarrow 2)$ from $^{(11)}$ below we get.

Every homomorphism of a topological right loop (respect.
homogeneous algebra) is open.

$^{(9)}$. The notion of homogeneous algebra, presented here, differs
from the notion of homogeneous algebra from \cite{Chob2}. I believe
that the term "homogeneous algebra" is not successful: it
should be concretized. Though obviously there is no connection between
the notion of homogeneous algebra and the notion of homogeneous space, I
think that the term "homogeneous algebra" appeared under the influence of
the notion of strongly homogeneous space from \cite{Usp2}
(or something similar). Recall, that for a compact $K$ by
\text{Aut} $K$ we denote the group of all autohomeomorphisms of
space $K$, provided with compact open topology. A compact $K$ is
called strongly homogeneous if there exists a continuous mapping $(x,
y) \rightarrow h_{x,y}$ from $K^2$ into $\text{Aut}$ $K$ such that
$h_{x,y}(x) = y$ for all $x, y \in K$.

According to (23) any homogeneous algebra coincides with a
topological right loop. By Corollary 1 the underling space of any
topological right loop if homogeneous. Consequently, the underling
space of any homogeneous algebra is topologically homogeneous.
 Moreover, from item 5) of Theorem 11 and
Proposition 15 from \cite{Usp2} it follows: a compact $K$ is a
homogeneous algebra if $K$ is a strongly homogeneous space.

$^{(10)}$. At the first glance it seems that the senseless
expression $µ(x, x, z) = \break q(x, p(y, z))$ is a typo. But
similar senseless statements can be found and in other works of
the author (see, for
 example, \cite{AC2},\cite{DC}).

$^{(11)}$. Strange things occur in connection with the proof of the
equivalence of conditions 1) and 2). The authors of this
paper \cite{AC1} specify that this equivalence is proved in
\cite{Chob2}. But according to assertions $^{(6)}$ - $^{(7)}$ the
proof of equivalence of conditions 1), 2) is and absurdity, it cannot
be called a mathematical proof. It is a set of basically unproved and
unrelated statements. As is mentioned in $^{(5)}$ even the
notion of rectifiable space is not correct.

The equivalence of conditions 1), 2) is proved in \cite{Usp1},
\cite{Usp2}, \cite{Gul}. These works give various correct
definition of rectifiable space. Clearly, the authors of paper
\cite{AC1} (at least M. M. Choban) were aware of these papers, but
do not make any references to them and assign this result to M. M.
Choban. On the other hand, it is weird and surprising that the
authors of \cite{Usp1}, \cite{Usp2}, \cite{Gul} assign this result
to M. M. Choban, but fail to indicate where is was published.

The term $''$rectifiable space$''$ and the equivalence of items
1), 2) are copied from \cite {Gul}, but this is not mentioned in the article.
To hide this fact item 2) is written in such a
$"$confusing$"$ form. The following would be more correct: On  underling
space $G$ of rectifiable space $G$ there exist a structure of a
homogeneous algebra.

The authors do not understand (or are not capable to understand) the
meaning and importance of Mal'cev's Theorem 2. The notion of
biternary algebra is introduced in Theorem 2 as a characteristic
notion (if and only if) for primitive classes of algebraic
systems, is introduced as derived operation of algebraic systems.
The implication $3) \rightarrow 2)$ is false. This follows from
item 5) of Proposition 5.

$^{(12)}$. According to (22) any guasigroup $(G, p, r, l)$ is
defined by identities $p(x, l(x, y)) = y$,  $l(x, p(x, y)) = y$,
$p(r(y, x), x) = y$, $r(p(y, x), x) = y$ with respect to
multiplication $p$, left division $l$ and right division $r$. The
remaining identities $l(r(x, y), x) = r(x, l(y, x)) =
y$ follow from the first and it is not clear why they are given.
There exists six parastrophs with respect to the basic operation
$(\cdot)$ for any quasigroup $(Q, \cdot, \backslash, /)$ (by our
terminology - equasigroup) \cite{Bruk}, \cite{Bel}. Then the given
list can be complemented with many identities. It seems that the authors are
not familiar with definition of equasigroup.

$^{(13)}$.  This expression is too general, it does not contain any
information. Unlike \cite {Mal2}, it does not indicate the relation
between quasigroup operation and loop operation. We concretize
$^{(13)}$ for topological quasigroups.
Recall that if $Q_1$ and $Q_2$ are topological quasigroups then we
say that $Q_1$ and $Q_2$ are \textit{isotopic} if there are
homeomorphisms $f, g, h : Q_1 \rightarrow Q_2$ such that $f(xy) =
g(x)h(y)$. The triple $(f, g, h)$ is called an \textit{isotopy
from $Q_1$ to $Q_2$}. The concept of isotopy is one of basic tools
in theory of topological projective planes.  Then from [16, pag.
5] ( it follows one from forms of expression $^{(7)}$. Let $Q$ be
a topological quasigroup and $e$ any fixed element of $Q$. If we
take any homeomorphism $f : Q \rightarrow Q$ with $f(e) = ee$ (e.
g. the one defined by $f(x) = L((ee)/e)x = ((ee)/e)x$), then the
underlying space of $Q$ together with the multiplication $x \circ
y = f^{-1}((f(e)/e)(e\backslash f(y)))$ define a topological loop
$L$ with unit $e$; and the triple $(f, R(e)^{-1}f, L(e)^{-1}f)$ is
an isotopy from $Q$ to $L$.

$^{(14)}$. Unclear and senseless assertion. According to
item 5) of Theorem 11 any homogeneous algebra is a topological
right loop with right unit. Then instead of the erroneous $p (e, x) =
x $ more correct would be $p(x,e) = x$. It is not clear why operations
$(+)$, $(\cdot)$ were introduced. Probably the authors wanted
to show the affinity of homogeneous algebras with classical notions of
algebras, rings.

$^{(15)}$. Ingenious logic. From the absolutely illogical and unrelated
statements (copied from different sources) the authors guessed and
actually generalized the following result from \cite{Usp2}: every
compact topological quasigroup is rectifiable. The statement $^{(16)}$
holds. This follows from Proposition 11.

Propositions 4, 5, and Corollary 8, which use the notion of
rectifiable space and the connected to it homogeneous algebra,
Mal'cev space, are not on the above-mentioned list of senseless
and false statements: Theorems 2, 3, 4, 5, Propositions 1, 2, and
Example 7. But Propositions 4, 5, and Corollary 8 could be added on
that list. As a confirmation of that we will analyze Proposition 13 and its proof.

\textmd{Proposition 5. $\ll$Let $X$ be a subalgebra of a
homogeneous algebra $G$. If the space $G$ is regular$\gg^{(17)}$
and Lindel$\ddot{o}$f, and the space $X$ is of pointwise countable
type and is dense in $G$, then there exist a separable metrizable
homogeneous algebra $G'$ and a homomorphism $g : G \rightarrow G'$
such that $X = g^{-1}(g(X))$ and the mapping $g$ is open and
perfect. In particular, it follows that $X$ is a Lindel$\ddot{o}$f
$p$-space.}

\textmd{Proof. $\ll$By the assumptions, there is a pair of binary
continuous operations $p, q : G \times G \rightarrow G$ on the
space $G$ such that:}

\textmd{$p(x, x) = p(y, y)$, and $p(x, q(x, y)) = y$, $q(x, p(x,
y)) = y$ for all $x, y \in G$;}

\textmd{$p(x, y) \in X$ and $q(x, y) \in X$ for all $x, y \in X$.}

\textmd{We put $e = p(x, x)$. If $a \in G$, then $p_a(x) = p(a,
x)$ and $q_a(x) = q(a, x)$ for any $x \in G$. We have $q^{(-1)}_ a
= p_a$ and $q_a(e) = a$. $\ll$Thus, $p_a$ and $q_a$ are
homeomorphisms. Moreover, $p_a(X) = q_a(X) = X$ for each $a \in
X$.$\gg^{(19)}$$\gg^{(18)}$}

$\ll$\textmd{$\ll$Let $F$ be a non-empty compact subspace of $X$
with a countable base of open neighborhoods in $X$.$\gg^{(21)}$ We
can assume that $e \in F$. $\ll$ Since $X$ is dense in $G$, the
set $F$ also has a countable base of open neighborhoods in the
space $G$. Therefore, $X$ and $G$ are $p$-spaces$\gg^{(9)}$ (see
[6], Proposition 2.1).}$\gg^{(20)}$

$\ll$\textmd{Since $F$ is a compact $G_{\delta}$-subset of the
Lindel$\ddot{o}$f algebra $G$, there exist a separable metrizable
homogeneous algebra $G'$ and a homomorphism $g: G \rightarrow G'$
such that $F = g^{-1}(g(F))$ and $g$ is a perfect mapping
\cite{Chob2}.$\gg^{(23)}$ $\ll$The quotient homomorphism of a
Mal’cev algebra is an open mapping \cite{Chob2}.$\gg^{(24)}$
$\ll$Thus, the mapping $g$ is open.$\gg^{(25)}$ We can assume that
$e' = g(e)$ and $p(z, z) = e'$ for any $z \in G'$.}

$\ll$\textmd{Let $b \in g(X) \subseteq G'$. Fix $a \in X \cap
g^{-1}(b)$. If $H = g^{-1}(e'$, then $H \subseteq F \subseteq X$
and $q_a(H) = g^{-1}(qb(e') \subseteq X$. Thus, $g^{-1}(b) =
g^{-1}(q_b(e') \subseteq X$. Therefore, $g^{-1}(g(X)) = X$. The
proof is complete.}$\gg^{(26)}$

$^{(17)}$. The contents of Proposition 5 is ill-considered. From the
equivalence $1) \leftrightarrow 2)$ of $^{(11)}$ and
\cite[Corollary 2.2]{Gul} it follows that the underlying space of
every homogeneous algebra is regular. Then the expression
"regular" from $^{(17)}$ is unnecessary.

$^{(18)}$, $^{(19)}$. The text $^{(18)}$ is copied from
\cite{Chob2}, but with corrected definition of homogeneous
algebra. Unlike \cite{Chob2} the statements from $^{(19)}$
hold, but they are not proved, only facts are presented. This can be
explained by item d) above. Let us present evidence. We use item 5) of
Theorem 11, i.e. let's be assume that the homogeneous algebra
$(G, q, p, e)$ is a topological right loop with right unit $e$.
From $q(a, p(a,x)) =
x$ it follows that $q_a p_a = E$, where $E$ denote the identical
mapping. But from $q_a p_a = E$ it does not follow that $q_a^{(-1)} =
p_a$. This equality and the assertion "$p_a$ and $q_a$ are
homeomorphisms" follow from Theorem 9. The equalities $p_a(X) =
q_a(X) = X$ follow from the property that $(X, q, p, e)$ is a right
loop. We mentioned that $^{(19)}$ is not used further in the proof of
Proposition 5.

$^{(20)}$. The implication \cite[Proposition 2.1]{AC2}
$\Rightarrow$ $^{(20)}$ is false, does not correspond to reality.

$^{(21)}$. It is necessary to show that space $X$ contains a
non-empty subspace $F$, at least from the following reasons.
Every topological group that contains a non-empty compact
subspace with a countable base of open neighborhoods is a
paracompact $p$-space \cite[pag.
122]{Arh5}.

$^{(9)}$. The excerpt is copied from \cite[end of proof of
Theorem 2.18]{Arh5}.

$^{(23)}$. Another deceptive statement. It is not clear
how $^{(23)}$ follows from the results of \cite{Chob2}.
Moreover, we have shown above that the results of \cite{Chob2}
are either false, or not proved.  Indeed, $^{(23)}$ follows from the
theorem of V. V. Filippov on
preservation of the class of $p$-space by perfect mappings
\cite{Fil} (see, \cite[pag. 120]{Arh5}).

$^{(24)}$. In \cite{Chob2} are not considered Mal'cev algebras as
defined in $^{(8)}$. Hence, the assertion $^{(24)}$ is not
proved, either. In such a case, assertion $^{(25)}$ is also not proved.
For more details see $^{(8)}$.

$^{(26)}$. In this excerpt only the words "Let, fix,
then, $\ldots$, the proof is complete." are correct. All the rest are
delirium, senseless assertions. For example, "Fix $a \in X \cap g^{-1}(b)$.".
Obviously the case $a \in g^{-1}(b) \backslash X$ should be considered, as well.
Further, "If $H = g^{-1}(e'$, then $H \subseteq F \subseteq X$ and
$q_a(H) = g^{-1}(qb(e') \subseteq X$. Thus, $g^{-1}(b) =
g^{-1}(q_b(e') \subseteq X$. Therefore, $g^{-1}(g(X)) = X$." The role
of inclusions $H \subseteq F \subseteq X$ is not clear in the proof of
$^{(26)}$. In general, $q_a(H) \subseteq
g^{-1}(qb(e')$. It is not clear, why $g^{-1}(q_b(e') \subseteq X$.
It is easier to understand all this if to consider the classes $K_a =
g^{-1}$, $a \in G$, of congruence, induced by  homomorphism $g$, and to
take into account that $K_a$ is a complete class of congruence, as well
as item 6) of Corollary 3 and the analysis of \cite[Lemma 3.1]{Chob2},
presented above.

$^{(27)}$. Paper \cite{ChobCalm} is published in the same issue of
the scientific journal "Buletinul Academiei de \c{S}tiin\c{t}e a
Republicii Moldova"  as \cite{AC1}. With the accuracy of
translation the paper \cite{ChobCalm} is copied word by word from
monograph \cite{CalChob2} (for example, sections 4, 5 of
\cite{ChobCalm} Coincide literarily with $\S\S$ 4.14, 4,15
respectively. This is not mentioned in \cite{ChobCalm}. It is also
not shown that the results of monograph \cite{CalChob2} coincide
exactly, including the paragraph numbering, with the results of
monograph \cite{Calm1} (with some minor differences in Chapter 7).

I got interested in other papers that consider the rectifiable
spaces. Let us start with \cite{AC2}. The main result of this
paper is the Theorem 3.1, known as Dichotomy Theorem.

\textbf{Theorem 3.1} \textit{For every rectifiable space $G$, any
remainder of $G$ in a compactification $bG$ is either
pseudocompact or Lindel$\ddot{o}$f.}

A crucial moment in proof of this theorem is the Corollary 2.8,
which is a direct corollary of Proposition 2.7. The latter is a
generalization of Proposition 2.6. Let us present it with its proof.

\textbf{Proposition 2.6.} \textit{Let $p, q : G \times G
\rightarrow G$ be a structure of homogeneous algebra on a space
$G$, $F$ and $\Phi$ be two non-empty compact subsets of $G$, $H =
p(F \times \Phi)$ and $L = q(F \times \Phi)$. Suppose further that
$\frak{B}_1$ is a base (a $\pi$-base) of $G$ at $F$ and that
$\frak{B}_2$ is a base (a $\pi$-base) of $G$ at $F$. Then:}

\textit{– $\{p(F \times U): U \in \frak{B}_2\}$ and $\{p(V \times
U): V \in \frak{B}_1, U \in \frak{B}_1\}$ are bases ($\pi$-bases,
respectively) of $G$ at $H$; and}

\textit{– $\{q(F \times U): U \in \frak{B}_2\}$ and $\{q(V \times
U): V \in \frak{B}_1, U \in \frak{B}_1\}$ are bases ($\pi$-bases,
respectively) of $G$ at $LH$.}

\textbf{Proof.} Let $F_1$ and $F_2$ be two compact subsets of $G$,
$U$ be an open subset of $G$ and $m(F_1 \times F_2) \subseteq U$
(where $m$ is either $p$ or $q$). Then, by Wallace theorem [11,
Theorem 3.2.10], there exist open subsets $U_1$ and $U_2$ of $G$
such that $F_1 \subseteq U_1$, $F_2 \subseteq U_2$ and $m(U_1
\times U_2) \subseteq U.$ $\square$

Let us now analyze the next Section $''$4. On paratopological
groups and dissentive spaces$''$ from paper \cite{AC2}, which
starts with the definition: a group $G$  T is called: – a
paratopological group if the multiplication $(x, y) \rightarrow x
\cdot y$ is a continuous mapping of $G \times G$ onto $G$. Recall
that every topological group $(G, \cdot, ^{-1}, e)$ is an
associative topological loop $(G, \cdot, \backslash, /, e)$. From
here it follows that

\textit{* if $G$ is a paratopological group (respect. topological
right loop) then from conditions $a \in G$, $H$ is an open set does
not result that the set $a\cdot H$ is open. On the contrary $G$
is a topological group (respect. topological loop).}

Further one may find a genuine nightmare (Proposition 4.1,
Corollaries 4.2, 4.3, 4.5, and Theorem 4.4). To justify this we
present the proof of Proposition 4.1, which has a crucial importance
for the proof of other statements.

\textbf{Proposition 4.1.} \textit{Let $G$ be a paratopological
group. Then:}

1. \textit{If there exists a non-empty compact subset of $G$ of
countable character in $G$, then $G$ is a space of countable
type.}

2. \textit{If there exists a non-empty compact subset of $G$ of
countable $\delta$-charac\-ter, then $G$ is a space of countable
$\delta$-type.}

\textbf{Proof.} Let $\{U_n: n \in \omega\}$ be a sequence of open
subsets of $G$, $F$ be a non-empty compact subset of $G$, $F
\subseteq cl_XU_n$ $U_{n+1} \subseteq U_n$ for each $n \in
\omega$. We may assume that $e \in F$.

Fix a compact subset $L$ of $G$. Put $\Phi = L \cdot H$ and $V_n =
L \cdot H_n$ for any $n \in \omega$. Then $\Phi$ is a compact
subset of $G$, $L \subseteq \Phi \subseteq  cl_X V_n$, and $V_n =
\cup\{Q_a(U_n) =a \cdot U_n: a \in L\}$ is open in $G$ for each $n
\in \omega$.

Let $\{U_n: n \in \omega \}$ be a $\pi$-base of $G$ at $F$. Then,
by Proposition 2.6, $\{V_n: n \in \omega\}$ is a $pi$-base of $G$
at $\Phi$. Thus, $\Phi$ is a compact subset of countable
$\delta$-character in $G$. Statement 2 is proved.

Let $\{U_n: n \in \omega \}$ be a base of $G$ at $F$. Then $\{V_n:
n \in \omega\}$ is a base of $G$ at $\Phi$. Thus, $\Phi$ is a
compact subset of countable character in $G$. Statement 1 is
proved.

This is not a mathematical proof. Actually, this could be regarded as
a compilation of unclear and unrelated denotations and topological
notions. The meaning of denotations $cl_XU_n$, $H$, $H_n$,  $Q_a(U_n)$
is not explained. It is not explained why $\Phi$ is a compact set of $G$.
 It is stated that $V_n = \cup\{Q_a(U_n) = a \cdot U_n: a \in L\}$ is open
 in $G$. But this contradicts the statement presented before Proposition 4.1.

Moreover, according to item 5) of Theorem 11 the Proposition 2.6
is proved for topological right loops. For such loops the
operation $(\backslash)$ is continuous. Then Proposition 2.6
cannot be used for paratopological groups.

In the beginning of the analyzed paper \cite{AC2} it is mentioned
that it will use the terminology from \cite{Eng}. Let us present
the following notions and results from \cite{Eng}. Let $\{X_s\}_{s
\in S}$ be a family of topological spaces such that
$$X_s \cap X_{s^{\prime}} =
\emptyset \quad \text{for} \quad s \neq s^{\prime}. \eqno{(14)}$$.
Let $X = \cup_{s \in S}X_s$. We define the set $U \subseteq X$
open in $X$ if and only if $U \cap X_s$ is an open set in $X_s$
for each $s \in S$. The family of all such defined open sets is a
topology on $X$.  The set $X$ together with this topology is
called \textit{sum of spaces $\{X_s\}_{s \in S}$} and is
designated $\bigoplus_{s \in S}X_s$ or $X_1 \bigoplus X_2
\bigoplus \ldots \bigoplus X_n$ if $S = \{1, 2, \ldots, n\}$.

It is also possible to define the sum $\bigoplus_{s \in S}X_s$ for
the family of topological spaces $\{X_s\}$ without the condition
(14), but in \cite{Eng} it is proved that the sum defined in such
a way is homeomorphic to the sum of topological spaces with
condition (14). An analogical notion is considered in \cite{Mal2},
the sum $\bigoplus_{s \in S}X_s$ without condition (14) is called
\textit{topological union of spaces $X_s$} and the sum
$\bigoplus_{s \in S}X_s$ with condition (14) is called
\textit{free topological union of spaces $X_s$}. We mentioned also
that the notion of \textit{singleton set}  is used to define the
topology of topological algebras with the help of limit of partial
algebraical operations with infinite of arguments.

We present three excerpts from Section 4 of \cite{AC2} for
consideration.

(*) A \textit{Mal’cev operation} on a space $X$ is a continuous
mapping $\mu: X^3 \rightarrow X$ such that $\mu(x, x, z) = z$ and
$\mu(x, y, y) = x$, for all $x, y, z \in X$. A space is called a
\textit{Mal’cev space} if it admits a Mal’cev operation (see
[8–10,14,15,19]). Let $\{X_{\alpha}: \alpha \in A\}$ be a
non-empty family of non-empty Mal’cev spaces. Assume that the set
$A$ is well ordered. Denote by $\mu_{\alpha}$ the Mal’cev
operation on $X_{\alpha}$. Let $X$ be the discrete sum of the
spaces $\{X_{\alpha}: \alpha \in A\}$. If $\alpha \in  A$ and $x,
y, z \in X_{\alpha}$, then $\mu(x, y, z) = \mu_{\alpha}(x, y, z)$.
Let $x, y, z \in X$ and $\{x, y, z\} \ X_{\alpha} \neq \emptyset$
for each $\alpha \in A$. Then $\{\mu(x, y, z)\} = \{x, y, z\} \cap
X_{\alpha}$, where $\alpha$ is the first element $\beta$ of $A$
for which $\{x, y, z\} \cap X_{\beta}$ is a singleton set. Then
$\mu$ is a Mal’cev operation on $X$. Thus, the discrete sum of
Mal’cev spaces is a Mal’cev space.

(**) A \textit{dissentive} operation (a \textit{dissentor}) on a
space $X$ is a continuous mapping $\mu: X^3 \rightarrow X$
satisfying the following conditions:

– $\mu(x, x, y) = y$ for all $x, y \in X$;

– for every open set $U$ of $X$ and all $b, c \in X$ the set
$\mu(U, b, c) = \{\mu(x, b, c): x \in U\}$ is open in $X$.

A space is \textit{dissentive} if it admits a dissentive
operation. Every rectifiable space is dissentive. If $p, q : G
\times G \rightarrow G$ is a structure of homogeneous algebra on a
space $G$, then $\mu(x, y, z) = p(x, q(y, z))$ is a Mal’cev
dissentive operation.

(***) \textbf{Example 5.18.} Let $G_1$ and $G_2$ be two disjoint
topological groups with the following properties:

1. $G_1$ is metrizable and is not locally compact;

2. $G_2$ is pseudocompact and is not compact.

Then the free topological sum $G = G_1 \oplus G_2$ is a space with
the following properties.

\textbf{Property 1.} $G$ is a Mal’cev non-dissentive space. In
particular, $G$ is not rectifiable.

\textbf{Property 2.} Every remainder of $G$ is not pseudocompact.

\textbf{Property 3.} Every remainder of $G$ is not
Lindel$\ddot{o}$f.

Therefore, the Dichotomy Theorem cannot be extended to all Mal’cev
spaces.

The proof of paragraph (*)is a mockery of the readers. It is not
clear what the discrete sum of topological spaces means.

Assume that the Mal'cev operation $\mu$ on discrete sum $X$ is
defined, though it is not so. Indeed, let $x_{\alpha} \in
X_{\alpha}$, $x_{\beta} \in X_{\beta}$, $x_{\gamma} \in
X_{\gamma}$ and let $x_{\alpha} < x_{\beta} < x_{\gamma}$. From
$\{x_{\alpha}, x_{\beta}, x_{\gamma}\} \cap X_{\alpha} =
\{x_{\alpha}\}$ it follows that $\mu(x_{\beta}, x_{\alpha},
x_{\gamma}) = x_{\alpha}$. From here it follows that $\mu(x, y,
z)$ is not a Mal'cev operation. Moreover, unlike the sum of
topological spaces, defined above, there isn’t any connection
between the topology of discrete sum $X$ and the topologies of
spaces $\{X_{\alpha}: \alpha \in A\}$. Now let us present excerpts
from \cite{AC2}.

Senseless properties and incorrect conclusion. The main topological object
of research of this paper is the remainder of a space $X$ (for example, Theorem 3.1).
This is the subspace $bX \backslash X$ of a Hausdorff compactification $bX$ of $X$.
 However, taking into account the Example 5.18 it seems that the authors are unaware
 of the definition of Hausdorff compactification.
Example 5.13 also seems weird, stating that Theorems 3.1, 3.4 do
not generalize to homogeneous spaces, though in excerpt $^{(4)}$
it is stated that the rectifiable spaces and homogeneous spaces
coincide topologically.

Further, let us move on to the next excerpt.

It is clear that statement $''$Every rectifiable space is
dissentive$''$ aims at showing the importance of introducing the
notion of $''$dissentive operation$''$. But this statement,
included in the paper without any proof, is incorrect. Indeed,
according to Yheorem 11  every homogeneous algebra $(G, p, q)$ is
a topological right loop $(G, \cdot, \backslash, e)$, which is not
a left loop. If $p = (\cdot)$, $q = (\backslash)$ then $\mu(x, y,
z) = x \cdot (y \backslash z)$. Let $H$ be an open set of $G$ and
fix the element $y, z \in G$.. According to definition of
dissentive operation the set $\mu(H, y, z)$ should be open, i.e.
the set $H \cdot (y \backslash z)$ should be open. This
contradicts the statement *. In a similar way we get a
contradiction for $q = (\cdot)$, $p = (\backslash)$.

The following question arises: was it necessary to introduce a new term
$''$dissentive operation$''$? For a space $X$ this definition is
equivalent to assertion: for any $a, b \in X$ the continuous mapping
$\varphi_{ab}(x) = \mu(x, a, b)$ is open and $\varphi_{ab}(a) = b$.

The last assertion without condition $\varphi_{ab}(a) = b$ is
changed and new terms, such as o-homogeneous, io-homogeneous,
co-homogeneous, do-homo\-geneous are introduced in paper \cite{ACM1}.
Besides, other new terms as c-sequence, $A$-sieve, densely sieve-complete,
densely q-complete, densely fan-complete and others are introduced.
These new notions, in tandem with plenty of topological terms,
used both with and without sense, make it impossible to verify
the authenticity of the paper’s results. Papers \cite{ACM2},
\cite{ACM3} are written in the same manner.

\textbf{Note:}\texttt{} Part 2 of this work will follow soon, where we will examine
the other papers, in particular: \cite{ChobChob}, \cite{Chir1},
\cite{ChobChir44}.

Nicolae I. Sandu,

Tiraspol State University of Moldova,

Chisin\u{a}u, R. Moldova

sandumn@yahoo.com

\end{document}